\documentclass[11pt]{article}

\usepackage{amsmath}
\usepackage{amssymb}
\usepackage{amsthm}
\usepackage{latexsym}
\usepackage{color}
\usepackage{enumerate}
\usepackage{graphicx}
\usepackage{color} 

\usepackage[T1]{fontenc}
\usepackage[english]{babel}
\usepackage{geometry}
\DeclareSymbolFont{calletters}{OMS}{cmsy}{m}{n}
\DeclareSymbolFontAlphabet{\mathcal}{calletters}

\def\be{\begin{eqnarray}}
\def\ee{\end{eqnarray}}

\def\b*{\begin{eqnarray*}}
\def\e*{\end{eqnarray*}}

\newtheorem{Theorem}{Theorem}[part]
\newtheorem{Definition}{Definition}[part]
\newtheorem{Proposition}{Proposition}[part]

\newtheorem{Assumption}{Assumption}[part]
\newtheorem{Lemma}{Lemma}[part]

\newtheorem{Remark}{Remark}[part]

\makeatletter \@addtoreset{equation}{section}

\@addtoreset{Definition}{section}

\@addtoreset{Theorem}{section}

\@addtoreset{Proposition}{section}

\@addtoreset{Property}{section}

\@addtoreset{Assumption}{section}

\@addtoreset{Corollary}{section}

\@addtoreset{Lemma}{section}

\@addtoreset{Remark}{section}

\@addtoreset{Example}{section}

\@addtoreset{Properties}{section}

\@addtoreset{footnote}{page}

\newcommand{\No}[1]{\left\|#1\right\|}     
\newcommand{\abs}[1]{\left|#1\right|}     

\def \D{\mathbb{D}}
\def \E{\mathbb{E}}
\def \F{\mathbb{F}}
\def \H{\mathbb{H}}
\def \L{\mathbb{L}}

\def \N{\mathbb{N}}
\def \P{\mathbb{P}}

\def \R{\mathbb{R}}

\def\Ac{{\cal A}}

\def\Dc{{\cal D}}

\def\Fc{{\cal F}}

\def\Kc{{\cal K}}
\def\Lc{{\cal L}}

\def\Nc{{\cal N}}

\def\Pc{{\cal P}}

\def\Rc{{\cal R}}

\def\Uc{{\cal U}}
\def\Vc{{\cal V}}

\def\Yc{{\cal Y}}
\def\Zc{{\cal Z}}

\def\Eh{{\hat E}}
\def\Fh{{\hat F}}

\def\Tr#1{{\rm Tr}\left[#1\right]}

\def \Sup{\displaystyle\sup}

\def\einf{{\rm ess \, inf}}
\def\esup{{\rm ess \, sup}}
\def\trace{{\rm Tr}}

\def\={\;=\;}
\def\.{\;.}

\def\eps{\varepsilon}

\def\reff#1{{\rm(\ref{#1})}}

\def\1{{\bf 1}}
\def \ep{\hbox{ }\hfill{ ${\cal t}$~\hspace{-5.1mm}~${\cal u}$   } }

\def \proof{{\noindent \bf Proof. }}
\def \ep{\hbox{ }\hfill$\Box$}

 \def\normeL2#1{\left\|{#1}\right\|_{L^2}}

 \title{Second Order BSDEs with Jumps: Existence and probabilistic representation for fully-nonlinear PIDEs}

\author{ Nabil {\sc Kazi-Tani}\footnote{CMAP, Ecole Polytechnique, Paris, mohamed-nabil.kazi-tani@polytechnique.edu.} \and Dylan {\sc Possama\"{i}}\footnote{CEREMADE, Universit\'e Paris-Dauphine, Paris, possamai@ceremade.dauphine.fr. Part of this work was carried out while the author was working at CMAP, Ecole Polytechnique,  whose financial support is kindly acknowledged. }      \and Chao {\sc Zhou}\footnote{Department of Mathematics, National University of Singapore, Singapore, matzc@nus.edu.sg. Part of this work was carried out while the author was working at CMAP, Ecole Polytechnique,  whose financial support is kindly acknowledged.} }          
 \date{\today}

\begin{document}

 \maketitle

\vspace{3mm}

 \begin{abstract}
In this paper, we pursue the study of second order BSDEs with jumps (2BSDEJs for short) started in our accompanying paper \cite{kpz3}. We prove existence of these equations by a direct method, thus providing complete wellposedness for 2BSDEJs. These equations are a natural candidate for the probabilistic interpretation of some fully non-linear partial integro-differential equations, which is the point of the second part of this work. We prove a non-linear Feynman-Kac formula and show that solutions to 2BSDEJs provide viscosity solutions of the associated PIDEs.
\vspace{10mm}

\vspace{5mm}

\end{abstract}
\newpage

\section{Introduction}

Motivated by duality methods and maximum principles for optimal stochastic control, Bismut studied in \cite{bis} a linear backward stochastic differential equation (BSDE). In their seminal paper \cite{pardpeng}, Pardoux and Peng generalized such equations to the non-linear Lipschitz case and proved existence and uniqueness results in a Brownian framework. Since then, a lot of attention has been given to BSDEs and their applications, not only in stochastic control, but also in theoretical economics, stochastic differential games and financial mathematics. Given a filtered probability space $(\Omega,\mathcal F,\left\{\mathcal F_t\right\}_{0\leq t\leq T},\mathbb P)$ generated by an $\mathbb R^d$-valued Brownian motion $B$, solving a BSDE with generator $g$, and terminal condition $\xi$ consists in finding a pair of progressively measurable processes $(Y,Z)$ such that
\begin{align}
 Y_t=\xi +\int_t^T g_s(Y_s,Z_s)ds-\int_t^T Z_s dB_s,\text{ }\mathbb P-a.s, \text{ } t\in [0,T]. \label{def_bsde}
\end{align}

The process $Y$ defined this way is a possible generalization of the conditional expectation of $\xi$, since when $g$ is the null function, we have $Y_t=\E^{\P}\left[\xi | \Fc_t\right]$, and in that case, $Z$ is the process appearing in the $(\Fc_t)$-martingale representation property of $(\E^{\P}\left[\xi | \Fc_t\right])_{t\geq 0}$. In the case of a filtered probability space generated by both a Brownian motion $B$ and a Poisson random measure $\mu$ with compensator $\nu$, the martingale representation for $(\E^{\P}\left[\xi | \Fc_t\right])_{t\geq 0}$ becomes
\begin{align*}
 \E^{\P}[\xi | \Fc_t] =\xi+ \int_0^t Z_s dB_s + \int_0^t \int_{\R^d\backslash \{0\}} U_s(x)(\mu-\nu)(dx,ds),\ \mathbb P-a.s.,
\end{align*}
where $U$ is a predictable function. This leads to the following natural generalization of equation (\ref{def_bsde}) to the case with jumps. We will say that $(Y,Z,U)$ is a solution of the BSDE with jumps (BSDEJ in the sequel) with generator $g$ and terminal condition $\xi$ if for all $t \in [0,T]$, we have $\mathbb P-a.s.$
\begin{align}
 Y_t=\xi +\int_t^T g_s(Y_s,Z_s,U_s)ds-\int_t^T Z_s dB_s -\int_t^T \int_{\R^d\backslash \{0\}} U_s(x)(\mu-\nu)(dx,ds). \label{def_bsdej}
\end{align}

Tang and Li \cite{tangli} were the first to prove existence and uniqueness of a solution for (\ref{def_bsdej}) in the case where $g$ is Lipschitz in $(y,z,u)$. In the continuous framework, Soner, Touzi and Zhang \cite{stz} generalized the BSDE \reff{def_bsde} to the second order case. Their key idea in the definition of the second order BSDEs (2BSDEs) is that the equation has to hold $\P$-almost surely, for every $\P$ in a class of non-dominated probability measures. Furthermore, they proved a uniqueness result using a representation result of the 2BSDEs as essential supremum of standard BSDEs. 

%

\vspace{0.4em}
Our aim in this paper is to pursue the study undertaken in \cite{kpz3}. More precisely, we prove existence of a solution to equation \reff{2bsdej} by a direct approach. Inspired by the representation obtained in Theorem $4.1$ of \cite{kpz3}, we construct a solution by using the tool of regular conditional probability distributions. This gives a complete wellposedness theory for 2BSDEJs. 


\vspace{0.4em}
The last part of our study is to establish a connection with partial integro-differential equations (PIDEs for short). Indeed, Soner, Touzi and Zhang proved in \cite{stz} that Markovian 2BSDEs, are connected in the continuous case to a class of parabolic fully non-linear PDEs. On the other hand, we know that solutions to standard Markovian BSDEJs provide viscosity solutions to some parabolic partial integro-differential equations whose non-local operator is given by a quantity similar to $\langle \widetilde{v},\nu\rangle$ defined in \reff{def_v_nu} (see \cite{bbp} for more details). Then in the Markovian case, 2BSDEJs are the natural candidates for the probabilistic interpretation of fully non-linear PIDEs. This is the purpose of the second part of this article. During the revision of this paper, in two beautiful articles, Neufeld and Nutz \cite{nn2,nn3} constructed so-called non-linear L\'evy processes, and showed that they provided probabilistic representations for viscosity solutions to a certain class of fully non-linear PIDEs. These objects are related to 2BSDEJs in the sense that they roughly correspond to the case of generator equal to $0$. However, the method they used for their construction (which is actually and extension of Nutz and van Handel \cite{nvh} to the Skorohod space of c\`adl\`ag functions) allows them to do not assume any strong pathwise regularity, unlike in our approach. Nonetheless, an extension of their method to the case of a non-zero generator is far from trivial, as it would require to study measurability of fully non-linear (and not only sub-linear) stochastic kernels. 

\vspace{0.4em}
The rest of the paper is organized as follows. In Section \ref{section.1}, in order to introduce our readers to the theory, we provide several definitions and results on the set of probability measures on the Skorohod space $\D$ that we will work with. In Section \ref{sec.2BSDE}, we introduce the generator of our 2BSDEJs and the assumptions under which we will be working, we recall from \cite{kpz3} the natural spaces and norms for the solution of a 2BSDEJ, and give the formulation of the 2BSDEJs. Section \ref{section.3} is devoted to the proof of our existence result. Finally, in Section \ref{sec.PIDE}, we study the links between solutions to some fully-nonlinear PIDEs and 2BSDEJs. The Appendix is dedicated to the proof of some important technical results needed throughout the paper.

\section{Preliminaries on probability measures} \label{section.1}
\subsection{The stochastic basis}
Let $\Omega:= \D([0,T],\mathbb R^d)$ be the space of c\`adl\`ag paths defined on $[0,T]$ with values in $\R^d$ and such that $w(0)=0$, equipped with the Skorohod topology, so that it is a complete, separable metric space (see \cite{bil} for instance). 

\vspace{0.5em}
We denote $B$ the canonical process, $\mathbb F:=\left\{\mathcal F_t\right\}_{0\leq t\leq T}$ the filtration generated by $B$, $\mathbb F^+:=\left\{\mathcal F_t^+\right\}_{0\leq t\leq T}$ the right limit of $\mathbb F$ and for any $\mathbb P$, $\overline{\mathcal F}_t^\mathbb P:=\mathcal F_t^+\vee\mathcal N^\mathbb P(\mathcal F_t^+)$ where
$$\mathcal N^\mathbb P(\mathcal G):=\left\{E\in\Omega,\text{ there exists $\widetilde E\in\mathcal G$ such that $E\subset\widetilde E$ and $\mathbb P(\widetilde E)=0$}\right\}.$$

We then define as in \cite{stz} a local martingale measure $\mathbb P$ as a probability measure such that $B$ is a $\mathbb P$-local martingale. 
We then associate to the jumps of $B$ a counting measure $\mu_{B}$, which is a random measure on $\mathbb R^+\times E$ equipped with its Borel $\sigma$-field $\mathcal B(\R^+)\times\mathcal B(E)$ (where $E:=\mathbb R^d\backslash \{0\}$), defined pathwise by
\begin{equation}
\mu_{B}(A,[0,t]) := \sum_{0<s\leq t} \mathbf{1}_{\{\Delta B_s \in A\}}, \; \forall t \geq 0, \; \forall A \in\mathcal B(E).
\end{equation}


We recall that (see for instance Theorem I.4.18 in \cite{jac}) under any local martingale measure $\P$, we can decompose $B$ uniquely into the sum of a continuous local martingale, denoted by $B^{\P,c}$, and a purely discontinuous local martingale, denoted by $B^{\P,d}$. Then, we define $\overline{\mathcal P}_W$ as the set of all local martingale measures $\mathbb P$, such that $\mathbb P$-a.s.
\begin{itemize}
	  \item[$\rm{(i)}$] The quadratic variation of $B^{\P,c}$ is absolutely continuous with respect to the Lebesgue measure $dt$ and its density takes values in $\mathbb S^{>0}_d$, which is the space of all $d\times d$ real valued positive definite matrices.
	\item[$\rm{(ii)}$] The compensator $\lambda^\mathbb P_t(dx,dt)$ of the jump measure $\mu_B$ exists under $\mathbb P$ and can be decomposed as follows
	$$\lambda^\P_t(dx,dt)=\nu^\P_t(dx)dt,$$
	for some $\F$-predictable random measure $\nu^\P$ on $E$.
\end{itemize}

We will denote by $\widetilde\mu_{B}^\mathbb P(dx,dt)$ the corresponding compensated measure, and for simplicity, we will often call $\nu^\P$ the compensator of the jump measure associated to $B$. Finally, as shown in \cite{kpz3}, it is possible, using results of Bichteler \cite{bich} to give a pathwise definition of the density (with respect to the Lebesgue measure) of the continuous part of $[B,B]$, which we denote by $\widehat a$.

\subsection{Martingale problems and probability measures}
In this section, we recall the families of probability measures introduced in \cite{kpz3}. Let $\mathcal{N}$ be the set of $\mathbb F$-predictable random measures $\nu$ on $\mathcal{B}(E)$ satisfying 
\begin{equation}
\int_0^t \int_{E}(1\wedge \abs{x}^2)\nu_s(\omega,dx)ds <+\infty \text{  and  } \int_0^t\int_{ \abs{x}>1 } x \nu_s(\omega,dx)ds <+\infty,\ \text{for all } \omega \in \Omega, \label{hyp_nu}
\end{equation}
and let $\Dc$ be the set of $\mathbb F$-predictable processes $\alpha$ taking values in $\mathbb S_d^{>0}$ with $\int_0^T|\alpha_t(\omega)|dt<+\infty$, for all $\omega \in \Omega$. We define a martingale problem as follows

\begin{Definition}
For $\mathbb F$-stopping times $\tau_1\leq\tau_2$, for $(\alpha,\nu)\in\mathcal D\times\mathcal N$ and for a probability measure $\mathbb P_1$ on $\mathcal F_{\tau_1}$, we say that $\P$ is a solution of the \textit{martingale problem} $(\P_1,\tau_1,\tau_2,\alpha,\nu)$ if 
	\begin{itemize}
		\item[$\rm{(i)}$] $\P = \P_1$ on $\Fc_{\tau_1}$.
		\item[$\rm{(ii)}$] The canonical process $B$ on $[\tau_1,\tau_2]$ is a semimartingale under $\P$ with characteristics 
			\begin{align*}
				\left( -\int_{\tau_1}^{\cdot} \int_E x \mathbf{1}_{\abs{x}>1}\nu_s(dx)ds, \int_{\tau_1}^{\cdot} \alpha_s ds,\ \nu_s(dx)ds \right).
			\end{align*}
	\end{itemize}
\end{Definition}

%
%
%
%
We say that the martingale problem associated to $(\alpha,\nu)$ has a unique solution if, for every stopping times $\tau_1, \tau_2$ and for every probability measure $\P_1$, the martingale problem $(\P_1,\tau_1,\tau_2,\alpha,\nu)$ has a unique solution.

\vspace{0.5em}
Let now $\overline{\Ac}_W$ be the set of $(\alpha, \nu) \in \Dc \times \Nc$, such that there exists a solution to the martingale problem $(\P^0,0,+\infty,\alpha,\nu)$, where $\P^0$ is such that $\P^0(B_0=0)=1$. We also denote by $\Ac_W$ the set of $(\alpha, \nu) \in \overline{\Ac}_W$ such that there exists a unique solution to the martingale problem $(\P^0,0,+\infty,\alpha,\nu)$. We denote $\P^{\alpha}_{\nu}$ this unique solution and set 
\begin{align*}
\Pc_W := \left\{ \P^{\alpha}_{\nu}, \ (\alpha, \nu) \in \Ac_W \right\}.
\end{align*}

Our main interest in this paper will be a particular sub-class of $\Pc_W$, which can be defined as follows. First, we define
$$\mathcal A_{\rm det}:=\left\{(I_d,F),\ F\in\mathcal N \text{ and $F$ is deterministic}\right\},$$
 and we let $\widetilde{\Ac}_{\rm det}$ be the separable class of coefficients generated by $\Ac_{\rm det}$ (to avoid unnecessary technicalities, we will refrain from giving the precise definition here, and we refer instead the reader to Definition $A.2$ in \cite{kpz3} for more details). We emphasize that thanks to Proposition $A.1$ in \cite{kpz3}, we have $\widetilde{\Ac}_{\rm det}\subset \Ac_W$. For simplicity, we let $\mathcal V$ designate the measure $F\in\mathcal N$ such that $(I_d,F)\in\widetilde{\mathcal A}_{\rm det}$. Moreover, we will still denote $\mathbb P_{0,F}:=\mathbb P_F^{I_d}$, for any $F\in\mathcal V$.

\vspace{0.5em}
Next, we introduce the following set $\mathcal R_F$ of $\F$-predictable functions $\beta:E\longmapsto \R$ such that for Lebesgue almost every $s\in[0,T]$
$$\abs{\beta_s}(\omega,x)\leq C(1\wedge\abs{x}),\ F_s(\omega,dx)-a.e.,\text{ for every }\omega\in\Omega,$$
and such that for every $\omega\in\Omega$, $x\longmapsto\beta_s(\omega,x)$ is strictly monotone on the support of the law of $\Delta B_s$ under $\P_{0,F}$.

%
\vspace{0.3em} 
Next, for each $F\in\mathcal V$ and for each $(\alpha,\beta)\in\mathcal D\times\mathcal R_F$, we define
$$\mathbb P^{\alpha,\beta}_F:=\mathbb P_{0,F}\circ\left(X^{\alpha,\beta}_.\right)^{-1},$$
where
\begin{equation}\label{X_alphabeta}
X^{\alpha,\beta}_t:=\int_0^t\alpha_s^{1/2}dB_s^{\P_{0,F},c}+\int_0^t\int_E\beta_s(x)\left(\mu_B(dx,ds)-F_s(dx)ds\right),\ \mathbb P_{0,F}-a.s.
\end{equation}

Finally, we let
$$\overline{\mathcal P}_S:=\underset{F\in\mathcal V}{\bigcup}\left\{\mathbb P^{\alpha,\beta}_F, \ (\alpha,\beta)\in\mathcal D\times\mathcal R_F\right\}.$$

We recall the following results from \cite{kpz3}
\begin{Lemma}
Every probability measure in $\overline{\mathcal P}_S$ satisfies the predictable martingale representation property and the Blumenthal $0-1$ law.
\end{Lemma}

\section{Preliminaries on 2BSDEJs}\label{sec.2BSDE}
\subsection{The Non-linear Generator}
In this subsection we will introduce the function which will serve as the generator of our 2BSDEJ. Let us define the following  spaces for $p\geq 1$ 
$$\hat{L}^p:= \left\{\xi,\ \mathcal F_T\text{-measurable, s.t. }\xi\in L^p(\nu),\text{ for every $\nu\in\Nc$}\right\}.$$

We then consider a map 
\begin{align*}
H_t(\omega,y,z,u,\gamma,\tilde{v}):[0,T]\times\Omega\times\mathbb{R}\times\mathbb{R}^d\times \hat{L}^2 \times D_1\times D_2\rightarrow \mathbb{R},
\end{align*}
where $D_1 \subset \mathbb{R}^{d\times d}$ is a given subset containing $0$ and $D_2 \subset \hat{L}^1$ is the domain of $H$ in $\tilde{v}$.

\vspace{0.5em}
Define the following conjugate of $H$ with respect to $\gamma$ and $\tilde{v}$ by
\begin{align}
F_t(\omega,y,z,u,a,\nu):=\underset{\{\gamma,\tilde{v}\} \in D_1 \times D_2}{\Sup}\Big\{\frac12\trace(a\gamma)+\int_{E} \tilde{v}(e) \nu(de)-H_t\big(\omega,y,z,u,\gamma,\tilde{v}\big)\Big\}, \label{F_fenchel}
\end{align}
for  $a \in \mathbb S_d^{>0}$ and $\nu \in \mathcal{N}$.

\vspace{0.5em}
In the remainder of this paper, we formulate the needed hypothesis for the generator directly on the function $F$, and the BSDEs we consider also include the case where $F$ does not take the form \reff{F_fenchel}. Nonetheless, this particular form allows to retrieve easily the framework of the standard BSDEs or of the $G$-stochastic analysis on the one hand (see sections 3.4 and 3.5 in \cite{kpz3}), and to establish the link with the associated PDEs on the other hand. In the latter cases, $H$ is evaluated at $\tilde{v}(\cdot) = Av(\cdot)$, where $A$ is the following non local operator, defined for any $\mathcal{C}^2$ function $v$ on $\R^d$ with bounded gradient and Hessian, and $y \in \R^d$ by: 
\begin{equation}\label{def_v_nu}
(Av)(y,e):= v(e+y) - v(y)- e. (\nabla v)(y), \text{  for } e \in E.
\end{equation}

The assumptions on $v$ ensure that $(Av)(y,\cdot)$ is an element of $\hat{L}^1$.

\vspace{0.5em}
The operator $A$ applied to $v$ will only appear again in Section \ref{sec.PIDE}, when we explore the links between 2BSDEJs and solutions to fully-nonlinear PIDEs. For the time being, we only want to insist on the fact that this particular non local form comes from the intuition that the 2BSDEJs is an essential supremum of standard BSDEJs. Indeed, solutions to Markovian BSDEJs provide viscosity solutions to some parabolic partial integro-differential equations with similar non-local operators (see \cite{bbp} for more details). 

\vspace{0.5em}
We define next $\widehat{F}^{\mathbb P}_t(y,z,u):=F_t(y,z,u,\widehat{a}_t,\nu^{\mathbb P}_t) \text{ and } \widehat{F}_t^{\mathbb P,0}:=\widehat{F}^{\mathbb P}_t(0,0,0).$ We also denote by $D^1_{F_t(y,z,u)}$ the domain of $F$ in $a$ and by $D^2_{F_t(y,z,u)}$ the domain of $F$ in $\nu$, for a fixed $(t,\omega,y,z,u)$. As in \cite{stz} we fix a constant $\kappa \in (1,2]$ and restrict the probability measures in $\mathcal{P}_H^\kappa\subset \mathcal{P}_{\widetilde{\mathcal A}}$

\begin{Definition}\label{def}
$\mathcal{P}_H^\kappa$ consists of all $\mathbb P \in \overline{\mathcal{P}}_{S}$ such that
\begin{itemize}
\item[$\rm{(i)}$] $\displaystyle\mathbb E^\mathbb P\left[\int_0^T\int_E \abs{x}^2{\nu}^{\mathbb P}_t(dx)dt\right]<+\infty.$
	\item[$\rm{(ii)}$] $\underline{a}^\mathbb P \leq \widehat{a}\leq \bar{a}^\mathbb P, \text{ } dt\times d\mathbb P-as \text{ for some } \underline{a}^\mathbb P, \bar{a}^\mathbb P \in \mathbb{S}_d^{>0}, \text{ and } \mathbb{E}^{\mathbb{P}}\left[\left(\int_0^T\abs{\widehat{F}_t^{\mathbb P,0}}^\kappa dt\right)^{\frac2\kappa}\right]<+\infty.$
\end{itemize}

\end{Definition}

\begin{Remark}
The above conditions assumed on the probability measures in $\mathcal P^\kappa_H$ ensure that under any $\mathbb P\in\mathcal P^\kappa_H$, the canonical process $B$ is actually a true square integrable c\`adl\`ag martingale. This will be important when we will define standard BSDEJs under each of these probability measures.
\end{Remark}

We now state our main assumptions on the function $F$
\begin{Assumption} \label{assump.href}
\rm{(i)} The domains $D^1_{F_t(y,z,u)}=D^1_{F_t}$, $D^2_{F_t(y,z,u)}=D^2_{F_t}$ are independent of $(\omega,y,z,u)$.

\vspace{0.4em}
\rm{(ii)} For fixed $(y,z,u,a,\nu)$, $F$ is $\mathbb{F}$-progressively measurable in $D^1_{F_t} \times D^2_{F_t} $.

\vspace{0.4em}
\rm{(iii)} The following uniform Lipschitz-type property holds. For all $(y,y',z,z',u,t,a,\nu,\omega)$
\begin{align*}
 &\abs{ F_t(\omega,y,z,u,a,\nu)- F_t(\omega,y',z',u,a,\nu)}\leq C\left(\abs{y-y'}+\abs{ a^{1/2}\left(z-z'\right)}\right).
\end{align*}

\vspace{0.4em}
\rm{(iv)} For all $(t,\omega,y,z,u^1,u^{2},a,\nu)$, there exist two processes $\gamma$ and $\gamma'$ such that
\begin{align*}
 \int_{E}\delta^{1,2} u(x)\gamma'_t(x)\nu(dx)\leq F_t(\omega,y,z,u^1,a,\nu)- F_t(\omega,y,z,u^2,a,\nu)&\leq \int_{E}\delta^{1,2} u(x)\gamma_t(x)\nu(dx),
\end{align*}
where $\delta^{1,2} u:=u^1-u^2$ and $c_1(1\wedge \abs{x}) \leq \gamma_t(x) \leq c_2(1\wedge \abs{x})$ with $-1+\delta\leq c_1\leq0, \; c_2\geq 0,$
and $c_1'(1\wedge \abs{x}) \leq \gamma'_t(x) \leq c_2'(1\wedge \abs{x})$ with $-1+\delta\leq c_1'\leq0, \; c_2'\geq 0,$
for some $\delta >0$.

\vspace{0.4em}
\rm{(v)} $F$ is uniformly continuous in $\omega$ for the Skorohod topology, that is to say that there exists some modulus of continuity $\rho$ such that for all $(t,\omega,\omega',y,z,u,a,\nu)$
$$\abs{F_t(\omega,y,z,u,a,\nu)-F_t(\omega',y,z,u,a,\nu)}\leq \rho\left(d_S(\omega_{.\wedge t},\omega'_{.\wedge t})\right),$$
where $d_S$ is the Skorohod metric and where $\omega_{.\wedge t}(s):=\omega(s\wedge t)$.

\end{Assumption}

\subsection{The Spaces and Norms}

We now define as in \cite{stz}, the spaces and norms which will be needed for the formulation of the 2BSDEJs.

\vspace{0.4em}
For $p\geq 1$, $L^{p,\kappa}_H$ denotes the space of all $\mathcal F_T$-measurable scalar r.v. $\xi$ with
$$\No{\xi}_{L^{p,\kappa}_H}^p:=\underset{\mathbb{P} \in \mathcal{P}_H^\kappa}{\sup}\mathbb E^{\mathbb P}\left[|\xi|^p\right]<+\infty.$$

$\mathbb H^{p,\kappa}_H$ denotes the space of all $\mathbb F^+$-predictable $\mathbb R^d$-valued processes $Z$ with
$$\No{Z}_{\mathbb H^{p,\kappa}_H}^p:=\underset{\mathbb{P} \in \mathcal{P}_H^\kappa}{\sup}\mathbb E^{\mathbb P}\left[\left(\int_0^T|\widehat a_t^{1/2}Z_t|^2dt\right)^{\frac p2}\right]<+\infty.$$

$\mathbb D^{p,\kappa}_H$ denotes the space of all $\mathbb F^+$-progressively measurable $\mathbb R$-valued processes $Y$ with
$$\mathcal P^\kappa_H-q.s. \text{ c\`adl\`ag paths, and }\No{Y}_{\mathbb D^{p,\kappa}_H}^p:=\underset{\mathbb{P} \in \mathcal{P}_H^\kappa}{\sup}\mathbb E^{\mathbb P}\left[\underset{0\leq t\leq T}{\sup}|Y_t|^p\right]<+\infty.$$

$\mathbb J^{p,\kappa}_H$ denotes the space of all $\mathbb F^+$-predictable functions $U$ with
$$\No{U}_{\mathbb J^{p,\kappa}_H}^p:=\underset{\mathbb{P} \in \mathcal{P}_H^\kappa}{\sup}\mathbb E^{\mathbb P}\left[\left(\int_0^T\int_{E}\abs{U_t(x)}^2 \nu^{\mathbb P}_t(dx)dt\right)^{\frac p2}\right]<+\infty.$$

For each $\xi \in L^{1,\kappa}_H$, $\mathbb P\in \mathcal P^\kappa_H$ and $t \in [0,T]$ denote
$$\mathbb E_t^{\mathcal P^\kappa_H,\mathbb P}[\xi]:=\underset{\mathbb P^{'}\in \mathcal P^\kappa_H(t^{+},\mathbb P)}{\esup^{\mathbb P}}\mathbb E^{\mathbb P^{'}}_t[\xi],\ \mathbb P-a.s., \text{ where } \mathcal P^\kappa_H(t^{+},\mathbb P):=\left\{\mathbb P^{'}\in\mathcal P^\kappa_H:\mathbb P^{'}=\mathbb P \text{ on }\mathcal F_t^+\right\}.$$

Then we define for each $p\geq \kappa$,
$$\mathbb L_H^{p,\kappa}:=\left\{\xi\in L^{p,\kappa}_H:\No{\xi}_{\mathbb L_H^{p,\kappa}}<+\infty\right\} \text{ where } \No{\xi}_{\mathbb L_H^{p,\kappa}}^p:=\underset{\mathbb P\in\mathcal P^\kappa_H}{\sup}\mathbb E^{\mathbb P}\left[\underset{0\leq t\leq T}{\esup}^{\mathbb P}\left(\mathbb E_t^{\mathcal P^\kappa_H,\mathbb P}[|\xi|^\kappa]\right)^{\frac{p}{\kappa}}\right].$$

Next, we denote by $\mbox{UC}_b(\Omega)$ the collection of all bounded and uniformly continuous maps $\xi:\Omega\rightarrow \mathbb R$ with respect to the Skorohod topology, and we let
$$\mathcal L^{p,\kappa}_H:=\text{the closure of $\mbox{UC}_b(\Omega)$ under the norm $\No{\cdot}_{\mathbb L^{p,\kappa}_H}$, for every $1<\kappa \leq p$}.$$

For a given probability measure $\mathbb P\in\mathcal P^\kappa_H$, the spaces $L^p(\mathbb P)$, $\mathbb D^p(\mathbb P)$, $\mathbb H^p(\mathbb P)$ and $\mathbb J^p(\mathbb P)$ correspond to the above spaces when the set of probability measures is reduced to the singleton $\left\{\mathbb P\right\}$. Finally,
$\mathbb H^{p}_{loc}(\mathbb P)$ (resp. $\mathbb J^{p}_{loc}(\mathbb P)$) denotes the space of all $\mathbb F^+$-predictable $\mathbb R^d$-valued processes $Z$ (resp.  $\mathbb F^+$-predictable functions $U$) with
$$\left(\int_0^T\abs{\widehat a_t^{1/2}Z_t}^2dt\right)^{\frac p2}<+\infty,\ (\text{resp. }\ \mathbb P-a.s.$$

$\mathbb J^{p}_{loc}(\mathbb P)$ denotes the space of all $\mathbb F^+$-predictable functions $U$ with
$$\left(\int_0^T\int_{E}\abs{U_t(x)}^2 \nu^{\mathbb P}_t(dx)dt\right)^{\frac p2}<+\infty,\ \mathbb P-a.s.$$

\subsection{Formulation}\label{formulation:2bsdej}

We shall consider the following 2BSDEJ, which must hold for $0\leq t\leq T$ and $\mathcal{P}^\kappa_H\text{-q.s.}$
\begin{equation}
Y_t=\xi +\int_t^T\widehat{F}^{\mathbb P}_s(Y_s,Z_s,U_s)ds -\int_t^T Z_sdB^{\mathbb P,c}_s-\int_t^T \int_{E} U_s(x) {\mu}^{\mathbb P}_B(dx,ds) + K^{\mathbb P}_T-K^{\mathbb P}_t.
\label{2bsdej}
\end{equation}

\begin{Definition}
We say $(Y,Z,U)\in \mathbb D^{2,\kappa}_H \times \mathbb H^{2,\kappa}_H \times \mathbb J^{2,\kappa}_H$ is a solution to the $2$BSDEJ \reff{2bsdej} if

\begin{itemize}
\item[$\bullet$] $Y_T=\xi$, $\mathcal{P}^\kappa_H$-q.s.
\item[$\bullet$] For all $\mathbb P \in \mathcal{P}^\kappa_H$ and $0\leq t\leq T$, the process $K^{\mathbb P}$ defined below is predictable and has non-decreasing paths $\mathbb P-a.s.$ 
\begin{equation}
K_t^{\mathbb P}:=Y_0-Y_t - \int_0^t\widehat{F}^{\mathbb P}_s(Y_s,Z_s,U_s)ds+\int_0^tZ_sdB^{\mathbb P,c}_s + \int_0^t \int_{E} U_s(x) {\mu}^{\mathbb P}_B(dx,ds).
\label{2bsdej_K}
\end{equation}

\item[$\bullet$] The family $\left\{K^{\mathbb P}, \mathbb P \in \mathcal P_H^\kappa\right\}$ satisfies the minimum condition
\begin{equation}
K_t^{\mathbb P}=\underset{ \mathbb{P}^{'} \in \mathcal{P}_H^\kappa(t^+,\mathbb{P}) }{ \einf^{\mathbb P} }\mathbb{E}_t^{\mathbb P^{'}}\left[K_T^{\mathbb{P}^{'}}\right], \text{ } 0\leq t\leq T, \text{  } \mathbb P-a.s., \text{ } \forall \mathbb P \in \mathcal P_H^\kappa.
\label{2bsdej.minK}
\end{equation}
\end{itemize}
If the family $\left\{K^{\mathbb P}, \mathbb P \in \mathcal P_H^\kappa\right\}$ can be aggregated into a universal process $K$, we call $(Y,Z,U,K)$ a solution of the $2$BSDEJ \reff{2bsdej}.
\end{Definition}


Following \cite{stz}, in addition to Assumption \ref{assump.href}, we will always assume
\begin{Assumption}\label{assump.h2ref}

\vspace{0.4em}
\rm{(i)} $\mathcal P_H^\kappa$ is not empty.

\vspace{0.4em}
\rm{(ii)} The process $F $ satisfies the following integrability condition 
\begin{equation}
\phi^{2,\kappa}_H:=\underset{\mathbb P\in\mathcal P^\kappa_H}{\sup}\mathbb E^{\mathbb P}\left[\underset{0\leq t\leq T}{\esup}^{\mathbb P}\left(\mathbb E_t^{\mathcal P^\kappa_H,\mathbb P}\left[\int^T_0|\Fh^{\mathbb P,0}_s|^\kappa ds\right]\right)^{\frac{2}{\kappa}}\right]<+\infty
\end{equation}
\end{Assumption}
\vspace{0.4em}

We recall the uniqueness result proved in \cite{kpz3}.

\begin{Theorem}\label{uniqueref}
Let Assumptions \ref{assump.href} and \ref{assump.h2ref} hold. Assume $\xi \in \mathbb{L}^{2,\kappa}_H$ and that $(Y,Z,U)$ is a solution to the $2$BSDEJ \reff{2bsdej}. Then, for any $\mathbb{P}\in\mathcal{P}^\kappa_H$ and $0\leq t_1< t_2\leq T$,
\begin{align}
\label{representationref}
Y_{t_1}&=\underset{\mathbb{P}^{'}\in\mathcal{P}^\kappa_H(t_1^+,\mathbb{P})}{\esup^\mathbb{P}}y_{t_1}^{\mathbb{P}^{'}}(t_2,Y_{t_2}), \text{ }\mathbb{P}-a.s.,
\end{align}
where, for any $\mathbb{P}\in\mathcal{P}^\kappa_H$, $\mathbb{F}^+$-stopping time $\tau$, and $\mathcal{F}^+_{\tau} $-measurable r.v. $\xi\in\mathbb{L}^2({\mathbb P})$, we denote by $(y^{\mathbb{P}},z^{\mathbb{P}},u^\P):=(y^{\mathbb{P}}(\tau,\xi),z^{\mathbb{P}}(\tau,\xi),u^\P(\tau,\xi))$ the solution to the following standard BSDEJ on $0\leq t\leq \tau$
\begin{equation}
\label{bsdej}
y^{\mathbb{P}}_t=\xi + \int_t^{\tau}\widehat{F}^{\mathbb P}_s(y^{\mathbb{P}}_s,z^{\mathbb{P}}_s,u^{\mathbb{P}}_s)ds-\int_t^{\tau}z^{\mathbb{P}}_sdB^{\mathbb P,c}_s - \int_t^{\tau} \int_{E} u^{\mathbb{P}}_s(x) {\mu}^{\mathbb P}_B(dx,ds), \text{ } \mathbb P-a.s.
\end{equation}
\end{Theorem}

\begin{Remark}\label{rem:bsde}
 We first emphasize that existence and uniqueness results for the standard BSDEJs \reff{bsdej} are not given directly by the existing literature, since the compensator of the counting measure associated to the jumps of $B$ is not deterministic. However, since all the probability measures we consider satisfy the martingale representation property and the Blumenthal $0-1$ law, it is clear that we can straightforwardly generalize the proof of existence and uniqueness of Tang and Li \cite{tangli} $($see also \cite{bech} and \cite{crep} for related results$)$. Furthermore, the usual \textit{a priori} estimates and comparison theorems will also hold.
 \end{Remark}
 
\section{A direct existence argument}\label{section.3}
The aim of this section is to prove the following result, which is our first main theorem.

\begin{Theorem}\label{mainref}
Let $\xi\in\mathcal L^{2,\kappa}_H$. Under Assumptions \ref{assump.href} and \ref{assump.h2ref}, there exists a unique solution $(Y,Z,U)\in\mathbb D^{2,\kappa}_H\times\mathbb H^{2,\kappa}_H\times\mathbb J^{2,\kappa}_H$ of the $2\rm{BSDEJ}$ \reff{2bsdej}.
\end{Theorem}

In the article \cite{stz}, the main tool to prove existence of a solution is the so-called regular conditional probability distributions of Stroock and Varadhan \cite{str}. Indeed, these tools allow to give a pathwise construction for conditional expectations. Since, at least when the generator is null, the $y$ component of the solution of a BSDE can be written as a conditional expectation, the r.c.p.d. allows us to construct solutions of BSDEs pathwise. Earlier in the paper, we have identified a candidate solution to the 2BSDEJ as an essential supremum of solutions of classical BSDEJs (see \reff{representationref}). However, those BSDEJs are written under mutually singular probability measures. Hence, being able to construct them pathwise allows us to avoid the problems related to negligible-sets. In this section we will generalize the approach of \cite{stz} to the jump case.

\subsection{Notations}\label{parag.notations}

For the convenience of the reader, we recall below some of the notations introduced in \cite{stz}. Remember that we are working in the Skorohod space $\Omega=\mathbb{D}\left([0,T],\mathbb R^d\right)$ endowed with the Skorohod metric, denoted $d_S$, which makes it a complete and separable space.

\vspace{0.4em}
$\bullet$ For $0\leq t\leq T$, we denote by $\Omega^t:=\left\{\omega\in \mathbb D\left([t,T],\mathbb R^d\right)\right\}$ the shifted canonical space of c\`adl\`ag paths on $[t,T]$ which are null at $t$, $B^t$ the shifted canonical process. $\mathbb F^t$ is the filtration generated by $B^t$. For any local martingale measure $\P$ on $(\Omega^t,\mathcal B(\Omega^t))$, we let $B^{t,\P,c}$, be the continuous local martingale part of $B^t$ and $B^{t,\P,d}$ its discontinuous martingale part. We again associate tot he jumps of $B^t$ a counting measure $\mu_{B^t}$, and we restrict ourselves to the set $\overline{P}^t_W$ of local martingale measures $\P$ such that
\begin{itemize}
	  \item[$\rm{(i)}$] The quadratic variation of $B^{t,\P,c}$ is absolutely continuous with respect to the Lebesgue measure $dt$ and its density takes values in $\mathbb S^{>0}_d$. Using Bichteler integration theory \cite{bich}, we can once more define this density pathwise and we denote it by $\widehat a^t$	\item[$\rm{(ii)}$] The compensator $\lambda^{t,\mathbb P}_s(dx,ds)$ of the jump measure $\mu_{B^t}$ exists under $\mathbb P$ and can be decomposed as follows
	$$\lambda^{t,\P}_s(dx,ds)=\nu^{t,\P}_s(dx)ds,$$
	for some $\F^t$-predictable random measure $\nu^{t,\P}$ on $E$.
\end{itemize}

We will denote by $\widetilde\mu_{B^t}^{t,\mathbb P}(dx,ds)$ the corresponding compensated measure. Let $\mathcal{N}^t$ be the set of $\mathbb F^t$-predictable random measures $\nu$ on $\mathcal{B}(E)$ satisfying 
\begin{equation}
\int_t^T \int_{E}(1\wedge \abs{x}^2)\nu_s(\widetilde\omega,dx)ds <+\infty \text{  and  } \int_t^T\int_{ \abs{x}>1 } x \nu_s(\widetilde\omega,dx)ds <+\infty, \; \forall \ \widetilde\omega \in \Omega^t,
\end{equation}

and let $\mathcal D^t$ be the set of $\mathbb F^t$-predictable processes $\alpha$ taking values in $\mathbb S_d^{>0}$ with $\int_t^T|\alpha_s(\widetilde\omega)|ds<+\infty$, for every $\widetilde\omega \in \Omega^t$. Exactly as in Section \ref{section.1}, we can define semimartingale problems and the corresponding probability measures. Define then
$$\mathcal A^t_{\rm det}:=\left\{(I_d,F),\ F\in\mathcal N^t \text{ and $F$ is deterministic}\right\},$$
let $\widetilde{\Ac}^t_{\rm det}$ be the separable class of coefficients generated by $\Ac_{\rm det}$, let $\mathcal V^t$ designate the measures $F\in\mathcal N^t$ such that $(I_d,F)\in\widetilde{\mathcal A}^t_{\rm det}$ and denote, for any $F\in\mathcal V^t$, by $\mathbb P_{t,F}$ the unique solution to the martingale problem associated to the couple $(I_d,F)$. Next, we define exactly as in Section \ref{section.1} a set $\mathcal R_F^t$ of $\F^t$-predictable functions $\beta:E\longmapsto \R$, and we define for each $F\in\mathcal V^t$ and for each $(\alpha,\beta)\in\mathcal D^t\times\mathcal R_F^t$
$$\mathbb P^{t,\alpha,\beta}_F:=\mathbb P_{t,F}\circ\left(X^{\alpha,\beta}_.\right)^{-1}.$$

Finally, we let
$$\overline{\mathcal P}^t_S:=\underset{F\in\mathcal V^t}{\bigcup}\left\{\mathbb P^{t,\alpha,\beta}_F, \ (\alpha,\beta)\in\mathcal D^t\times\mathcal R_F^t\right\},$$
and we emphasize that this set enjoys the same properties as $\overline{\Pc}_S$. We next define important operations on the shifted spaces and their paths.

\vspace{0.4em}
$\bullet$ For $0\leq s\leq t\leq T$ and $\omega\in \Omega^s$, we define the shifted path $\omega^t\in \Omega^t$ by
$$\omega^t_r:=\omega_r-\omega_t,\text{ }\forall r\in [t,T].$$

$\bullet$ For $0\leq s\leq t\leq T$ and $\omega\in \Omega^s$, $\widetilde \omega\in\Omega^t$ we define the concatenation path $\omega\otimes_t\widetilde \omega\in\Omega^s$ by 
$$(\omega\otimes_t\widetilde \omega)(r):=\omega_r1_{[s,t)}(r)+(\omega_{t}+\widetilde\omega_r)1_{[t,T]}(r),\text{ }\forall r\in[s,T].$$

$\bullet$ For $0\leq s\leq t\leq T$ and a $\mathcal F^s_T$-measurable random variable $\xi$ on $\Omega^s$, for each $\omega \in\Omega^s$, we define the shifted $\mathcal F^t_T$-measurable random variable $\xi^{t,\omega}$ on $\Omega^t$ by
$$\xi^{t,\omega}(\widetilde\omega):=\xi(\omega\otimes_t\widetilde \omega),\text{ }\forall \ \widetilde\omega\in\Omega^t.$$
Similarly, for an $\mathbb F^s$-progressively measurable process $X$ on $[s,T]$ and $(t,\omega)\in[s,T]\times\Omega^s$, we can define the shifted process $\left\{X_r^{t,\omega},r\in[t,T]\right\}$, which is $\mathbb F^t$-progressively measurable.

\vspace{0.4em}
$\bullet$ For a $\mathbb F$-stopping time $\tau$, we use the same simplification as \cite{stz}
$$\omega\otimes_\tau\widetilde \omega:=\omega\otimes_{\tau(\omega)}\widetilde \omega,\text{ }\xi^{\tau,\omega}:=\xi^{\tau(\omega),\omega},\text{ }X^{\tau,\omega}:=X^{\tau(\omega),\omega}.$$

$\bullet$ We define the "shifted" generator by
$$\widehat F^{t,\omega,{\mathbb P}}_s(\widetilde\omega,y,z,u):=F_s(\omega\otimes_t\widetilde\omega,y,z,u,\widehat a^t_s(\widetilde\omega), \nu{t,{\mathbb P}}_s(\widetilde\omega)), \text{ }\forall (s,\widetilde\omega)\in[t,T]\times\Omega^t.$$
Then note that since $F$ is assumed to be uniformly continuous in $\omega$ for the Skorohod topology, then so is $\widehat F^{t,\omega}$. Notice that this implies that for any $\mathbb P\in\overline{\mathcal P}^t_{S}$
$$\mathbb E^\mathbb P\left[\left(\int_t^T\abs{\widehat F^{t,\omega,{\mathbb P}}_s(0,0,0)}^{\kappa}ds\right)^{\frac2\kappa}\right]<+\infty,$$
for some $\omega$ if and only if it holds for all $\omega\in\Omega$.

\vspace{0.4em}
$\bullet$ We also extend Definition \ref{def} in the shifted spaces
\begin{Definition}\label{set_proba_shift}
$\mathcal P^{t,\kappa}_H$ consists of all $\mathbb P:=\mathbb P^{t,\alpha,\beta}_F\in\overline{\mathcal P}^t_{S}$ such that
\begin{itemize}
\item[\rm{(i)}] $\underline{a}^\mathbb P\leq\widehat a^t_s\leq\overline{a}^\mathbb P, \ ds\times d\mathbb P-a.s.\text{ for some }\underline{a}^\mathbb P,\overline{a}^\mathbb P\in\mathbb S^{>0}_d$ and
$\displaystyle\mathbb E^\mathbb P\left[\int_t^T\int_E \abs{x}^2{\nu}^{t,\mathbb P}_s(dx)ds\right]<+\infty.$
\item[\rm{(ii)}] The following integrability condition holds
$$\mathbb E^\mathbb P\left[\left(\int_t^T\abs{\widehat F^{t,\omega,{\mathbb P}}_s(0,0,0)}^{\kappa}ds\right)^{\frac2\kappa}\right]<+\infty,\text{ for all $\omega\in\Omega$}.$$
\end{itemize}
\end{Definition}

\vspace{0.4em}
$\bullet$ Finally, we define the so-called regular conditional probability distributions (r.c.p.d. in the sequel). For given $\omega\in \Omega$, $\mathbb F$-stopping time $\tau$ and $\mathbb P\in\mathcal P^\kappa_H$, the r.c.p.d. of $\mathbb P$ is a probability measure $\mathbb P^\omega_\tau$ on $\mathcal F_T$ such that for every bounded $\mathcal F_T$-measurable random variable $\xi$
$$\mathbb E^\mathbb P_\tau\left[\xi\right](\omega)=\mathbb E^{\mathbb P^\omega_\tau}[\xi],\text{ for $\mathbb P$-a.e. $\omega$.}$$

Besides, $\mathbb P^\omega_\tau$ naturally induces a probability measure $\mathbb P^{\tau,\omega}$ on $\mathcal F_T^{\tau(\omega)}$ such that the $\mathbb P^{\tau,\omega}$-distribution of $B^{\tau(\omega)}$ is equal to the $\mathbb P^\omega_\tau$-distribution of $\left\{B_t-B_{\tau(\omega)},\ t\in[\tau(\omega),T]\right\}.$ Besides, we have
$$\mathbb E^{\mathbb P^\omega_\tau}[\xi]=\mathbb E^{\mathbb P^{\tau,\omega}}[\xi^{\tau,\omega}].$$

\begin{Remark}
We emphasize that the above notations correspond to the ones used in \cite{stz} when we consider the subset of $\Omega$ consisting of all continuous paths from $[0,T]$ to $\mathbb R^d$ whose value at time $0$ is $0$.
\end{Remark}

We now prove that there exists a relation between $(\widehat a^{t,\omega},\left(\nu^{\mathbb P}\right)^{t,\omega})$ and $(\widehat a^t,\nu^{t,\mathbb P^{t,\omega}})$.

\begin{Proposition}\label{relationhata}
Let $\mathbb P\in\mathcal P^\kappa_H$ and $\tau$ be an $\mathbb F$-stopping time. Then, for $\mathbb P$-a.e. $\omega\in\Omega$, we have for $ds\times d\mathbb P^{\tau,\omega}$-a.e. $(s,\widetilde\omega)\in[\tau(\omega),T]\times\Omega^{\tau(\omega)}$
\begin{align*}
&\widehat a_s^{\tau,\omega}(\widetilde\omega)=\widehat a_s^{\tau(\omega)}(\widetilde\omega)\text{, and } (\nu_ s^{{\mathbb P}})^{\tau,\omega}(\widetilde\omega,A)=\nu_s^{\tau(\omega),\mathbb P^{\tau,\omega}}(\widetilde\omega,A)\text{ for every $A\in\mathcal B(E)$.}
\end{align*}
\end{Proposition}

This result is important for us, because it implies that for $\mathbb P$-a.e. $\omega\in\Omega$ and for $ds\times d\mathbb P^{t,\omega}-a.e. \ (s,\widetilde\omega)\in[t,T]\times\Omega^t$
$$F_s\left(\omega\otimes_t\widetilde\omega,y,z,u,\widehat a_s(\omega\otimes_t\widetilde\omega),\nu^{\mathbb P}_s(\omega\otimes_t\widetilde\omega)\right)=F_s\left(\omega\otimes_t\widetilde\omega,y,z,u,\widehat a_s^t(\widetilde\omega),\nu^{t,\mathbb P^{t,\omega}}_s(\widetilde\omega)\right),$$
which justifies the choice we made for the "shifted" generator.

\vspace{0.4em}
\proof
The proof of the equality for $\widehat a$ is the same as in Lemma $4.1$ of \cite{stz2}, so we omit it. Now, for $s\geq \tau$ and for any $A\in\mathcal B(E)$, we know by the Doob-Meyer decomposition that there exist a $\mathbb P$-local martingale $M$ and a $\mathbb P^{\tau,\omega}$-martingale $N$ such that
\begin{align*}
\mu_B([0,s],A)&=M_s+\int_0^s\nu^\mathbb P_r(A)dr,\ \mathbb P-a.s.\\
\mu_{B^{\tau(\omega)}}([\tau(\omega),s],A)&=N_s+\int_{\tau}^s\nu^{\tau(\omega),\mathbb P^{\tau,\omega}}_r(A)dr,\ \mathbb P^{\tau,\omega}-a.s.
\end{align*}

Then, we can rewrite the first equation above for $\mathbb P$-a.e. $\omega\in\Omega$ and for $\mathbb P^{\tau,\omega}$-a.e. $\widetilde\omega\in\Omega^{\tau(\omega)}$
\begin{equation}
\mu_B(\omega\otimes_\tau\widetilde\omega,[0,s],A)=M_s^{\tau,\omega}(\widetilde\omega)+\int_0^s\nu_r^{\mathbb P,\tau,\omega}(\widetilde\omega,A)dr.
\label{eq:trilili}
\end{equation}

Now, by definition of the measures $\mu_B$ and $\mu_{B^{\tau(\omega)}}$, we have
$$\mu_B(\omega\otimes_\tau\widetilde\omega,[0,s],A)=\mu_B(\omega,[0,\tau],A)+\mu_{B^{\tau(\omega)}}(\widetilde\omega,[\tau,s],A).$$

Hence, we obtain from \reff{eq:trilili} that for $\mathbb P$-a.e. $\omega\in\Omega$ and for $\mathbb P^{\tau,\omega}$-a.e. $\widetilde\omega\in\Omega^{\tau(\omega)}$
\begin{align*}
\mu_B(\omega,[0,\tau],A)-\int_0^\tau\nu^\mathbb P_r(\omega,A)dr+N_s(\widetilde\omega)-M_s^{\tau,\omega}(\widetilde\omega)=\int_\tau^s\left(\nu_r^{\mathbb P,\tau,\omega}(\widetilde\omega,A)-\nu^{\tau(\omega),\mathbb P^{\tau,\omega}}_r(\widetilde\omega,A)\right)dr
\end{align*}

On the left-hand side above, the terms which are $\mathcal F_\tau$-measurable are constants in $\Omega^{\tau(\omega)}$ and using the same arguments as in Step $1$ of the proof of Lemma \ref{lemme.technique}, we can show that $M^{\tau,\omega}$ is a $\mathbb P^{\tau,\omega}$-local martingale for $\mathbb P$-a.e. $\omega\in\Omega$. This means that the left-hand side is a $\mathbb P^{\tau,\omega}$-local martingale while the right-hand side is a predictable finite variation process. By the martingale representation property which still holds in the shifted canonical spaces, we deduce that for $\mathbb P$-a.e. $\omega\in\Omega$ and for $ds\times d\mathbb P^{\tau,\omega}$-a.e. $(s,\widetilde\omega)\in[\tau(\omega),T]\times\Omega^{\tau(\omega)}$
$$\int_\tau^s\left(\nu_r^{\mathbb P,\tau,\omega}(\widetilde\omega,A)-\nu^{\tau(\omega),\mathbb P^{\tau,\omega}}_r(\widetilde\omega,A)\right)dr=0,$$
which is the desired result.
\ep

\subsection{Existence when $\xi$ is in $\rm{UC_b}(\Omega)$ }\label{sec.existtt}

When $\xi$ is in $\rm{UC_b}(\Omega)$, we know that there exists a modulus of continuity function $\rho$ for $\xi$ and $F$ in $\omega$. Then, for any $0\leq t \leq s \leq T,\ (y,z,u)\in  \mathbb R \times \mathbb{R}^d\times\hat L^2$ and $\omega,\omega'\in \Omega,\ \widetilde{\omega}\in\Omega^t,\ {\mathbb P\in\mathcal P^{t,\kappa}_H}$,
\begin{align*}
\left|\xi^{t,\omega}\left(\widetilde{\omega}\right)-\xi^{t,\omega'}\left(\widetilde{\omega}\right)\right| \leq \rho\left(d_{S,t}(\omega,\omega')\right) \text{, } \left|\widehat{F}_s^{t,\omega,\mathbb P}\left(\widetilde{\omega},y,z,u\right)-\widehat{F}_s^{t,\omega',\mathbb P}\left(\widetilde{\omega},y,z,u\right)\right| \leq \rho\left(d_{S,t}(\omega,\omega')\right),
\end{align*}
where $d_{S,t}(\omega,\omega'):=d_S(\omega_{.\wedge t},\omega'_{.\wedge t}),$ $d_S$ being the Skorohod distance. We then define for all $\omega\in\Omega$
\begin{equation}
\Lambda\left(\omega\right):=\underset{0\leq s\leq t}{\sup}\Lambda_t\left(\omega\right):=\underset{0\leq s\leq t}{\sup}\ \underset{\mathbb P\in\mathcal P^{t,\kappa}_H}{\sup}\left(\mathbb E^\mathbb P\left[\abs{\xi^{t,\omega}}^2+\int_t^T|\widehat F^{t,\omega,\mathbb P}_s(0,0,0)|^2ds\right]\right)^{1/2}.
\end{equation}

Now since $\widehat F^{t,\omega,\mathbb P}$ is also uniformly continuous in $\omega$ for the Skorohod topology, it is easily verified that $\Lambda\left(\omega\right)<\infty \text{ for some } \omega\in\Omega \text{ iff it holds for all } \omega\in\Omega.$ Moreover, when $\Lambda$ is finite, it is uniformly continuous in $\omega$ for the Skorohod topology and is therefore $\mathcal F_T$-measurable. Now, by Assumption \ref{assump.h2ref}, we have $\Lambda_t\left(\omega\right)<\infty \text{ for all } \left(t,\omega\right)\in\left[0,T\right]\times\Omega.$ To prove existence, we define the following value process $V_t$ pathwise
\begin{equation}\label{sol}
V_t(\omega):=\underset{\mathbb P\in\mathcal P^{t,\kappa}_H}{\sup}\mathcal Y^{\mathbb P,t,\omega}_t\left(T,\xi\right), \text{ for all } \left(t,\omega\right)\in\left[0,T\right]\times\Omega,
\end{equation}
where, for any $\left(t_1,\omega\right)\in\left[0,T\right]\times\Omega,\ \mathbb P\in\mathcal P^{t_1,\kappa}_H, t_2\in\left[t_1,T\right]$, and any $\mathcal F_{t_2}$-measurable $\eta\in\mathbb L^{2}\left(\mathbb P\right) $, we denote $\mathcal Y^{\mathbb P,t_1,\omega}_{t_1}\left(t_2,\eta\right):= y^{\mathbb P,t_1,\omega}_{t_1}$, where $\left(y^{\mathbb P,t_1,\omega},z^{\mathbb P,t_1,\omega},u^{\mathbb P,t_1,\omega}\right) $ is the solution of the following BSDEJ on the shifted space $\Omega^{t_1} $ under $\mathbb P$
\begin{align}\label{eq.bsdeeeeref}
\nonumber y^{\mathbb P,t_1,\omega}_{s}&=\eta^{t_1,\omega}+\int^{t_2}_{s}\widehat{F}^{t_1,\omega,\mathbb P}_{r}\left(y^{\mathbb P,t_1,\omega}_{r},z^{\mathbb P,t_1,\omega}_{r},u_r^{\mathbb P,t_1,\omega}\right)dr-\int^{t_2}_{s}z^{\mathbb P,t_1,\omega}_{r}dB^{t_1,\mathbb P,c}_{r}\\
&\hspace{0.9em}-\int_s^{t_2}\int_{E}u_r^{\mathbb P,t_1,\omega}(x)\widetilde\mu^\mathbb P_{B^{t_1}}(dx,dr), \ \mathbb P-a.s.,\ s\in[t_1,t_2],
\end{align}
where as usual $\widetilde\mu^\mathbb P_{B^{t_1}}(dx,ds):=\mu_{B^{t_1}}(dx,ds)-\nu_s^{t_1,\mathbb P}(dx)ds$. In view of the Blumenthal $0-1$ law, $y^{\mathbb P,t,\omega}_{t}$ is constant for any given $\left(t,\omega\right)$ and $\mathbb P\in\mathcal P^{t,\kappa}_H$, and therefore the value process $V$ is well defined. Let us now show that $V$ inherits some properties from $\xi$ and $F$.

\begin{Lemma}\label{unifcont}
Let Assumptions \ref{assump.href} and \ref{assump.h2ref} hold and consider some $\xi$ in $\rm{UC_b}(\Omega)$. Then for all $\left(t,\omega\right)\in\left[0,T\right]\times\Omega$ we have $\left|V_t\left(\omega\right)\right|\leq C\Lambda_t\left(\omega\right) $. Moreover, for all $\left(t,\omega,\omega'\right)\in\left[0,T\right]\times\Omega^2$, 
$\left|V_t\left(\omega\right)-V_t\left(\omega'\right)\right|\leq C\rho\left(d_{S,t}(\omega,\omega')\right) $. Consequently, $V_t$ is $\mathcal F_t$-measurable for every $t\in\left[0,T\right]$.
\end{Lemma}

\proof
$\rm{(i)}$ For each $\left(t,\omega\right)\in\left[0,T\right]\times\Omega $ and $\mathbb P\in\mathcal P^{t,\kappa}_H $, let $\alpha$ be some positive constant which will be fixed later and let $\eta\in(0,1)$. Since $F$ is uniformly Lipschitz in $(y,z)$ and satisfies Assumption \ref{assump.href}$\rm{(iv)}$, we have
$$\abs{\widehat F_s^{t,\omega,\mathbb P}(y,z,u)}\leq \abs{\widehat F_s^{t,\omega,\mathbb P}(0,0,0)}+C\left(\abs{y}+|\left(\widehat a^t_s\right)^{1/2}z|+\left(\int_{E}\abs{u(x)}^2\nu^{t,\mathbb P}_s(dx)\right)^{1/2}\right).$$

Now apply It\^o's formula. We obtain
\begin{align*}
&e^{\alpha t}\abs{y_t^{\mathbb P,t,\omega}}^2+\int_t^Te^{\alpha s}\abs{(\widehat a^t_s)^{1/2}z_s^{\mathbb P,t,\omega}}^2ds+\int_t^T\int_{E}e^{\alpha s}\abs{u_s^{\mathbb P,t,\omega}(x)}^2\nu^{t,\mathbb P}_s(dx)ds\\
&= e^{\alpha T}\abs{\xi^{t,\omega}}^2+2\int_t^Te^{\alpha s}y_s^{\mathbb P,t,\omega}\widehat F_s^{t,\omega,\mathbb P}(y_s^{\mathbb P,t,\omega},z_s^{\mathbb P,t,\omega},u_s^{\mathbb P,t,\omega})ds\\
&\hspace{0.9em}-\alpha\int_t^Te^{\alpha s}\abs{y_s^{\mathbb P,t,\omega}}^2ds-2\int_t^Te^{\alpha s}y^{\mathbb P,t,\omega}_{s^-}z_s^{\mathbb P,t,\omega}dB^{t,\mathbb P,c}_s\\
&\hspace{0.9em}-\int_t^T\int_{E}e^{\alpha s}\left(2y^{\mathbb P,t,\omega}_{s^-}u_s^{\mathbb P,t,\omega}(x)+\abs{u_s^{\mathbb P,t,\omega}(x)}^2\right)\widetilde\mu^{\mathbb P}_{B^{t}}(dx,ds)\\
&\leq e^{\alpha T}\abs{\xi^{t,\omega}}^2+\int_t^Te^{\alpha s}\abs{\widehat F_s^{t,\omega,\mathbb P}(0,0,0)}^2ds+\left(1+2C+\frac{2C^2}{\eta}-\alpha\right)\int_t^Te^{\alpha s}\abs{y_s^{\mathbb P,t,\omega}}^2ds\\
&\hspace{0.9em}+\eta\int_t^Te^{\alpha s}\abs{(\widehat a^t_s)^{1/2}z_s^{\mathbb P,t,\omega}}^2ds+\eta\int_t^T\int_{E}e^{\alpha s}\abs{u_s^{\mathbb P,t,\omega}(x)}^2\nu^{t,\mathbb P}_s(dx)ds\\
&\hspace{0.9em}-2\int_t^Te^{\alpha s}y^{\mathbb P,t,\omega}_{s^-}z_s^{\mathbb P,t,\omega}dB^{t,\mathbb P,c}_s-\int_t^T\int_{E}e^{\alpha s}\left(2y^{\mathbb P,t,\omega}_{s^-}u_s^{\mathbb P,t,\omega}(x)+\abs{u_s^{\mathbb P,t,\omega}(x)}^2\right)\widetilde\mu^{\mathbb P}_{B^{t}}(dx,ds).
\end{align*}

Now choose $\eta=1/2$ and $\alpha$ large enough. By taking expectation we obtain easily $\abs{y_t^{\mathbb P,t,\omega}}^2\leq C\abs{\Lambda_t(\omega)}^2.$ The result then follows from the arbitrariness of $\mathbb P$.

\vspace{0.4em}
$\rm{(ii)}$ The proof is exactly the same as above, except that one has to use uniform continuity in $\omega$ of $\xi^{t,\omega}$ and $ F^{t,\omega}$. Indeed, for each $\left(t,\omega\right)\in\left[0,T\right]\times\Omega $ and $\mathbb P\in\mathcal P^{t,\kappa}_H $, let $\alpha$ be some positive constant which will be fixed later and let $\eta\in(0,1)$. By It\^o's formula we have, since $ F$ is uniformly Lipschitz
\begin{align*}
&e^{\alpha t}\abs{y_t^{\mathbb P,t,\omega}-y_t^{\mathbb P,t,\omega'}}^2+\int_t^T\scriptstyle e^{\alpha s}\left(\abs{(\widehat a^t_s)^{1/2}(z_s^{\mathbb P,t,\omega}-z_s^{\mathbb P,t,\omega'})}^2+\int_{E}\scriptstyle e^{\alpha s}(u_s^{\mathbb P,t,\omega}-u_s^{\mathbb P,t,\omega'})^2(x)\nu^{t,\mathbb P}_s(dx)\right)ds\\
&\leq e^{\alpha T}\abs{\xi^{t,\omega}-\xi^{t,\omega'}}^2+2C\int_t^Te^{\alpha s}\abs{y_s^{\mathbb P,t,\omega}-y_s^{\mathbb P,t,\omega'}}^2ds\\
&\hspace{0.9em}+2C\int_t^T\abs{y_s^{\mathbb P,t,\omega}-y_s^{\mathbb P,t,\omega'}}\abs{(\widehat a_s^t)^{1/2}(z_s^{\mathbb P,t,\omega}-z_s^{\mathbb P,t,\omega'})}ds\\
&\hspace{0.9em}+2C\int_t^Te^{\alpha s}\abs{y_s^{\mathbb P,t,\omega}-y_s^{\mathbb P,t,\omega'}}\left(\int_{E}\abs{u_s^{\mathbb P,t,\omega}(x)-u_s^{\mathbb P,t,\omega'}(x)}^2\nu^{t,\mathbb P}_s(dx)\right)^{1/2}ds\\
&\hspace{0.9em}+2C\int_t^Te^{\alpha s}\abs{y_s^{\mathbb P,t,\omega}-y_s^{\mathbb P,t,\omega'}}\abs{\widehat F^{t,\omega,\mathbb P}_s(y_s^{\mathbb P,t,\omega},z_s^{\mathbb P,t,\omega},u_s^{\mathbb P,t,\omega})-\widehat F^{t,\omega',\mathbb P}_s(y_s^{\mathbb P,t,\omega},z_s^{\mathbb P,t,\omega},u_s^{\mathbb P,t,\omega})}ds\\
&\hspace{0.9em}-\alpha\int_t^Te^{\alpha s}\abs{y_s^{\mathbb P,t,\omega}-y_s^{\mathbb P,t,\omega'}}^2ds-2\int_t^Te^{\alpha s}(y^{\mathbb P,t,\omega}_{s^-}-y_{s^-}^{\mathbb P,t,\omega'})(z_s^{\mathbb P,t,\omega}-z_s^{\mathbb P,t,\omega'})dB^{t,\mathbb P,c}_s\\
&\hspace{0.9em}-\int_t^T\int_{E}e^{\alpha s}\left(2(y^{\mathbb P,t,\omega}_{s^-}-y^{\mathbb P,t,\omega'}_{s^-})(u_s^{\mathbb P,t,\omega}-u_s^{\mathbb P,t,\omega'})+(u_s^{\mathbb P,t,\omega}-u_s^{\mathbb P,t,\omega'})^2\right)(x)\widetilde\mu^{\mathbb P}_{B^{t}}(dx,ds).
\end{align*}

We then deduce
\begin{align*}
&e^{\alpha t}\abs{y_t^{\mathbb P,t,\omega}-y_t^{\mathbb P,t,\omega'}}^2+\int_t^T\scriptstyle e^{\alpha s}\left(\abs{(\widehat a^t_s)^{\1/2}(z_s^{\mathbb P,t,\omega}-z_s^{\mathbb P,t,\omega'})}^2+\int_{E}\scriptstyle e^{\alpha s}(u_s^{\mathbb P,t,\omega}-u_s^{\mathbb P,t,\omega'})^2(x)\nu^{t,\mathbb P}_s(dx)\right)ds\\
&\leq e^{\alpha T}\abs{\xi^{t,\omega}-\xi^{t,\omega'}}^2+\int_t^Te^{\alpha s}\abs{\widehat F^{t,\omega,\mathbb P}_s(y_s^{\mathbb P,t,\omega},z_s^{\mathbb P,t,\omega},u_s^{\mathbb P,t,\omega})-\widehat F^{t,\omega',\mathbb P}_s(y_s^{\mathbb P,t,\omega},z_s^{\mathbb P,t,\omega},u_s^{\mathbb P,t,\omega})}^2ds\\[0.3em]
&\hspace{0.9em}+\eta\int_t^Te^{\alpha s}\abs{(\widehat a^t_s)^{1/2}(z_s^{\mathbb P,t,\omega}-z_s^{\mathbb P,t,\omega'})}^2ds+\eta\int_t^T\int_{E}e^{\alpha s}\abs{u_s^{\mathbb P,t,\omega}(x)-u_s^{\mathbb P,t,\omega'}(x)}^2\nu^{t,\mathbb P}_s(dx)ds\\
&\hspace{0.9em}+\left(2C+C^2+\frac{2C^2}{\eta}-\alpha\right)\int_t^Te^{\alpha s}\abs{y_s^{\mathbb P,t,\omega}-y_s^{\mathbb P,t,\omega'}}^2ds\\
&\hspace{0.9em}-2\int_t^Te^{\alpha s}(y^{\mathbb P,t,\omega}_{s^-}-y_{s^-}^{\mathbb P,t,\omega'})(z_s^{\mathbb P,t,\omega}-z_s^{\mathbb P,t,\omega'})dB^{t,{\mathbb P},c}_s\\
&\hspace{0.9em}-\int_t^T\int_{E}e^{\alpha s}\left(2(y^{\mathbb P,t,\omega}_{s^-}-y^{\mathbb P,t,\omega'}_{s^-})(u_s^{\mathbb P,t,\omega}-u_s^{\mathbb P,t,\omega'})+(u_s^{\mathbb P,t,\omega}-u_s^{\mathbb P,t,\omega'})^2\right)(x)\widetilde\mu^{\mathbb P}_{B^{t}}(dx,ds).
\end{align*}

Now choose $\eta=1/2$ and $\alpha$ such that $\iota:=\alpha -2C-C^2-\frac{2C^2}{\eta}\geq 0$. We obtain the desired result by taking expectation and using the uniform continuity in $\omega$ of $\xi$ and $F$.
\ep

\vspace{0.4em}
The next proposition is a dynamic programming property verified by the value process, which will be crucial when proving that $V$ provides a solution to the 2BSDEJ with generator $F$ and terminal condition $\xi$. The result and its proof are intimately connected to Proposition $4.7$ in \cite{stz2} and use the same type of arguments.

\begin{Proposition}\label{progdyn}
Under Assumptions \ref{assump.href}, \ref{assump.h2ref} and for $\xi\in \rm{UC_b}(\Omega)$, we have for all $0\leq t_1<t_2\leq T$ and for all $\omega \in \Omega$
$$V_{t_1}(\omega)=\underset{\mathbb P\in \mathcal P^{t_1,\kappa}_H}{\sup}\mathcal Y_{t_1}^{\mathbb P,t_1,\omega}(t_2,V_{t_2}^{t_1,\omega}).$$
\end{Proposition}

\begin{Remark}
Let us emphasize here that all the regularity in $\omega$ we assumed so far is because it is a sufficient condition in order to obtain the measurability and the dynamic programming property for \reff{sol}. It is however clear that such assumptions are too restrictive and we hope to be able to weaken them in a future work. In fact, in a recent paper, Nutz and van Handel \cite{nvh} showed the required regularity and dynamic programming for conditional non-linear expectations $($corresponding roughly to 2BSDEs with a generator equal to $0)$ for terminal conditions which were only upper semi-analytic. An extension of their result to our framework would allow us to get rid off our regularity assumptions, and therefore off the limitations induced by the continuity with respect to the Skorohod topology. Indeed, even for fairly regular functions $f$, the random variable $f(B_t)$ is continuous for the Skorohod distance only for almost every $t\in[0,T]$\footnote{As mentioned in the introduction, the results of \cite{nvh} have actually been extended very recently to a jump setting in \cite{nn2,nn3}. The authors do manage, in the case of a null generator, to construct what is actually a solution to the corresponding 2BSDEJ, without any regularity assumptions on the terminal condition.}  
\end{Remark}

The proof is almost the same as the proof in \cite{stz2}, with minor modifications due to the introduction of jumps. We therefore relegate it to the appendix.

\vspace{0.5em}
Now we are facing the problem of the regularity in $t$ of $V$. Indeed, if we want to obtain a solution to the 2BSDEJ, then it has to be at least c\`adl\`ag, $\mathcal P^\kappa_H-q.s.$
To this end, we define now for all $(t,\omega)$, the $\mathbb F^+$-progressively measurable process
$$V_t^+:=\underset{r\in\mathbb Q\cap(t,T],r\downarrow t}{\overline \lim}V_r.$$

\begin{Lemma}\label{lem.cadlag}
Under the conditions of the previous Proposition, we have
$$V_t^+=\underset{r\in\mathbb Q\cap(t,T],r\downarrow t}{\lim}V_r,\text{ }\mathcal P^\kappa_H-q.s.$$
and thus $V^+$ is c\`adl\`ag, $\mathcal P^\kappa_H-q.s.$
\end{Lemma}
The proof is relegated to the appendix.

\vspace{0.4em}
We follow now Remark $4.9$ in \cite{stz2}, and for a fixed $\mathbb P\in\mathcal P^\kappa_H$, we introduce the following reflected BSDE with jumps (RBSDEJ for short) and with lower obstacle $V^+$ under $\mathbb P$ 
\begin{align*}
&\widetilde{Y}_t^{\mathbb P}=\xi+\int_t^T\widehat{F}^\mathbb P_s(\widetilde{Y}_s^{\mathbb P},\widetilde{Z}_s^{\mathbb P},\widetilde{U}_s^{\mathbb P},\nu)ds-\int_t^T\widetilde{Z}_s^{\mathbb P}dB_s^{\mathbb P,c}-\int_t^T\int_{E}\widetilde{U}_s^{\mathbb P}(x)\widetilde{\mu}^{\mathbb P}_B(dx,ds)+\widetilde{K}^{\mathbb P}_T-\widetilde{K}^{\mathbb P}_t\\
&\widetilde{Y}_t^{\mathbb P}\geq V^+_t,\ 0\leq t\leq T,\ \mathbb P-a.s.\\
&\int_0^T\left(\widetilde{Y}_{s^-}^{\mathbb P}-V^{+}_{s^{-}}\right)d\widetilde{K}^{\mathbb P}_s=0,\text{ } \mathbb P-a.s.,
\end{align*}
where we emphasize that the process $\widetilde{K}^\mathbb P$ is predictable.

\begin{Remark}
Existence and uniqueness of the above RBSDEJ under our Assumptions, with the restrictions that the compensator is not random, have been proved by Hamad\`ene and Ouknine \cite{hamaou} or Essaky \cite{ess}. However, their proofs can be easily generalized to our context.
\end{Remark}

Let us now argue by contradiction and suppose that $\widetilde Y^\mathbb P$ is not equal $\mathbb P-a.s.$ to $V^+$. Then we can assume without loss of generality that $\widetilde Y^\mathbb P_0>V^+_0$, $\mathbb P-a.s.$. Fix now some $\eps>0$ and define the following stopping time
$$\tau^\eps:=\inf\left\{t\geq 0,\ \widetilde Y^\mathbb P_t\leq V^+_t+\eps\right\}.$$

Then $\widetilde Y^\mathbb P$ is strictly above the obstacle before $\tau^\eps$, and therefore $\widetilde K^\mathbb P$ is identically equal to $0$ in $[0,\tau^\eps]$. Hence, we have
$$\widetilde{Y}_t^{\mathbb P}=\widetilde Y^\mathbb P_{\tau^\eps}+\int_t^{\tau^\eps}\widehat{F}^{\mathbb P}_s(\widetilde{Y}_s^{\mathbb P},\widetilde{Z}_s^{\mathbb P},\widetilde{U}_s^{\mathbb P})ds-\int_t^{\tau^\eps}\widetilde{Z}_s^{\mathbb P}dB_s^{\mathbb P,c}-\int_t^{\tau^\eps}\int_{E}\widetilde{U}_s^{\mathbb P}(x)\widetilde{\mu}^{\mathbb P}_B(dx,ds),\ \mathbb P-a.s.$$

Let us now define the following BSDEJ on $[0,\tau^\eps]$
$$y^{+,\mathbb P}_t=V^+_{\tau^\eps}+\int_t^{\tau^\eps}\widehat{F}^{\mathbb P}_s(y_s^{+,\mathbb P},z_s^{+,\mathbb P},u_s^{+,\mathbb P})ds-\int_t^{\tau^\eps}z_s^{+,\mathbb P}dB_s^{\mathbb P,c}-\int_t^{\tau^\eps}\int_{E}u_s^{+,\mathbb P}(x)\widetilde{\mu}^{\mathbb P}_B(dx,ds),\ \mathbb P-a.s.$$

By the standard \textit{a priori} estimates already used in this paper, we obtain that
$$\widetilde Y_0^\mathbb P\leq y^{+,\mathbb P}_0+C\abs{V^+_{\tau^\eps}-\widetilde Y^\mathbb P_{\tau^\eps}}\leq y^{+,\mathbb P}_0+C\eps,$$
by definition of $\tau^\eps$. Following the arguments in Step $1$ of the proof of Theorem $4.5$ in \cite{stz2}, we can show that $y^{+,\mathbb P}_0\leq V^+_0$ which in turn implies $\widetilde Y_0^\mathbb P\leq V^+_0+C\eps,$
hence a contradiction by arbitrariness of $\eps$. Therefore, we have obtained the following decomposition
$$V_t^+=\xi+\int_t^T\widehat{F}^{\mathbb P}_s(V_s^+,\widetilde{Z}_s^{\mathbb P},\widetilde{U}_s^{\mathbb P})ds-\int_t^T\widetilde{Z}_s^{\mathbb P}dB_s^{\mathbb P,c}-\int_t^T\int_{E}\widetilde{U}_s^{\mathbb P}(x)\widetilde{\mu}^{\mathbb P}_B(dx,ds)+\widetilde{K}^{\mathbb P}_T-\widetilde{K}^{\mathbb P}_t,\ \mathbb P-a.s.$$

Finally, since $V^+$ and $B$ are c\`adl\`ag, we can use the result of Karandikar \cite{kar} to give a pathwise definition of $[B,V^+]$. We then have
$$d[Y,B]^c_t=d\langle Y^{\mathbb P,c},B^{\mathbb P,c}\rangle_t=\widetilde Z^\P_td\langle B^{\mathbb P,c},B^{\mathbb P,c}\rangle_t=\widetilde Z_t^\P d[B,B]^c_t=\widehat a_t\widetilde Z_t^\P dt, \text{ } \mathbb P-a.s., \ \forall{\mathbb P}\in\mathcal{P}^\kappa_H,$$
$$\Delta[Y,B]_t=\widetilde U_t^\P(\Delta B_t)\Delta B_t,\text{ } \mathcal{P}^\kappa_H-q.s.,$$
so that we can define aggregators $Z$ and $U$ for the the families $\{\widetilde Z^\mathbb P,\ \mathbb P\in\mathcal P^\kappa_H\},$ and $\{\widetilde U^\mathbb P,\ \mathbb P\in\mathcal P^\kappa_H\}.$

\vspace{0.4em}
We next prove the representation \reff{representationref} for $V$ and $V^+$, and that, as shown in Proposition $4.11$ of \cite{stz2}, we actually have $V=V^+$, $\mathcal P^{\kappa}_H-q.s.$, which shows that in the case of a terminal condition in ${\rm UC_b}(\Omega)$, the solution of the $2$BSDEJ is actually $\mathbb F$-progressively measurable. 

\begin{Proposition}\label{prop.repref}
Assume that $\xi\in {\rm UC_b}(\Omega)$ and that Assumptions \ref{assump.href} and \ref{assump.h2ref} hold. Then we have
$$V_t=\underset{\mathbb P^{'}\in\mathcal P_H(t,\mathbb P)}{\esup^\mathbb P}\mathcal Y_t^{\mathbb P^{'}}(T,\xi)\text{ and } V_t^+=\underset{\mathbb P^{'}\in\mathcal P_H^\kappa(t^+,\mathbb P)}{\esup^\mathbb P}\mathcal Y_t^{\mathbb P^{'}}(T,\xi), \text{ }\mathbb P-a.s., \text{ }\forall \mathbb P\in\mathcal P_H^\kappa.$$

Besides, we also have for all $t$, $V_t=V_t^+, \text{ }\mathcal P_H^\kappa-q.s.$
\end{Proposition}

\proof
The proof for the representations is the same as the proof of proposition $4.10$ in \cite{stz2}, since we also have a stability result for BSDEJs under our assumptions. For the equality between $V$ and $V^+$, we also refer to the proof of Proposition $4.11$ in \cite{stz2}.
\ep

\vspace{0.4em}
Therefore, in the sequel we will use $V$ instead of $V^+$. Finally, we have to check that the minimum condition \reff{2bsdej.minK} holds. Fix $\mathbb P$ in $\mathcal P^\kappa_H$ and $\mathbb P^{'}\in\mathcal P^\kappa_H(t^+,\mathbb P)$. Then, proceeding exactly as in Step $2$ of the proof of Theorem $4.1$ in \cite{kpz3}, but introducing the process $\gamma'$ of Assumption \ref{assump.href}$\rm{(iv)}$ instead of $\gamma$, we can similarly obtain
\begin{align}\label{eq:M'}
V_t-y_t^{\mathbb P^{'}}&\geq\mathbb E^{\mathbb P^{'}}_t\left[\int_t^TM^{' \mathbb P^{'}}_sd\widetilde K_s^{\mathbb P^{'}}\right]
\geq \mathbb E^{\mathbb P^{'}}_t\left[\underset{t\leq s\leq T}{\inf}M^{' \mathbb P^{'}}_s\left(\widetilde K_T^{\mathbb P^{'}}-\widetilde K_t^{\mathbb P^{'}}\right)\right],
\end{align}
where $M^{' \mathbb P^{'}}$ is defined as $M^{\mathbb P^{'}}$ but with $\gamma'$ instead of $\gamma$. Now let us prove that for any $n>1$
\begin{equation}
\mathbb E^{\mathbb P^{'}}_t\left[\left(\underset{t\leq s\leq T}{\inf}M^{' \mathbb P^{'}}_s\right)^{-n}\right]<+\infty,\ \mathbb P^{'}-a.s.
\label{eq:momentsexp}
\end{equation}

First we have
\begin{align*}
M^{' \mathbb P^{'}}_s&=\exp\left(\int_t^s\lambda_rdr+\int_t^s\eta_r\widehat a_r^{-1/2}dB_s^{\mathbb P^{'},c}-\frac12\int_t^s\abs{\eta_r}^2dr+\int_t^s\int_E\gamma'_r(x)\widetilde\mu^{\mathbb P^{'}}_B(dx,dr)\right)\\
&\hspace{0.9em}\times\prod_{t\leq r\leq s}(1+\gamma'_r(\Delta B_r))e^{-\gamma'_r(\Delta B_r)}.
\end{align*}

Define, then $A^{\mathbb P^{'}}_s=\mathcal E\left(\int_t^s\eta_r\widehat a_r^{-1/2}dB_s^{\mathbb P^{'},c}\right)$ and $C^{\mathbb P^{'}}_s=\mathcal E\left(\int_t^s\int_E\gamma'_r(x)\widetilde\mu^{\mathbb P^{'}}_B(dx,dr)\right)$. Notice that both these processes are strictly positive martingales, since $\eta$ and $\gamma'$ are bounded and we have assumed that $\gamma'$ is strictly greater than $-1$. We have
$$M^{' \mathbb P^{'}}_s=\exp\left(\int_t^s\lambda_rdr\right)A^{\mathbb P^{'}}_s C^{\mathbb P^{'}}_s.$$

Since the process $\lambda$ is bounded, we have
\begin{align*}
\left(\underset{t\leq s\leq T}{\inf}M^{' \mathbb P^{'}}_s\right)^{-n}&\leq C\left(\underset{t\leq s\leq T}{\inf}A^{\mathbb P^{'}}_sC^{\mathbb P^{'}}_s\right)^{-n}=C\left(\underset{t\leq s\leq T}{\sup}\left\{(A^{\mathbb P^{'}}_sC^{\mathbb P^{'}}_s)^{-1}\right\}\right)^{n}.
\end{align*}

\vspace{0.4em}
Using the Doob inequality for the submartingale $(A_sC_s)^{-1}$, we obtain
\begin{align*}
\mathbb E^{\mathbb P^{'}}_t\left[\left(\underset{t\leq s\leq T}{\inf}M^{' \mathbb P^{'}}_s\right)^{-n}\right]&\leq C\mathbb E^{\mathbb P^{'}}_t\left[(C^{\mathbb P^{'}}_T A^{\mathbb P^{'}}_T)^{-n}\right]\leq C\left(\mathbb E^{\mathbb P^{'}}_t\left[(C^{\mathbb P^{'}}_T)^{-2n}\right]\mathbb E^{\mathbb P^{'}}_t\left[(A^{\mathbb P^{'}}_T)^{-2n}\right]\right)^{1/2}<+\infty,
\end{align*}
where we used the fact that since $\eta$ is bounded, the continuous stochastic exponential $A^{\mathbb P^{'}}$ has negative moments of any order, and where the same result holds for the purely discontinuous stochastic exponential $C^{\mathbb P^{'}}$ by Lemma A.$4$ in \cite{kpz3}.

\vspace{0.4em}
Then, we have for any $p>1$
\begin{align*}
&\mathbb E^{\mathbb P^{'}}_t\left[\widetilde K_T^{\mathbb P^{'}}-\widetilde K_t^{\mathbb P^{'}}\right]\\
&=\mathbb E^{\mathbb P^{'}}_t\left[\left(\underset{t\leq s\leq T}{\inf}M^{' \mathbb P^{'}}_s\right)^{1/p}\left(\widetilde K_T^{\mathbb P^{'}}-\widetilde K_t^{\mathbb P^{'}}\right)\left(\underset{t\leq s\leq T}{\inf}M^{' \mathbb P^{'}}_s\right)^{-1/p}\right]\\
&\leq \left(\mathbb E^{\mathbb P^{'}}_t\left[\underset{t\leq s\leq T}{\inf}M^{' \mathbb P^{'}}_s\left(\widetilde K_T^{\mathbb P^{'}}-\widetilde K_t^{\mathbb P^{'}}\right)\right]\right)^{1/p}\left(\mathbb E^{\mathbb P^{'}}_t\left[\underset{t\leq s\leq T}{\inf}(M^{' \mathbb P^{'}}_{s})^{-\frac{2}{p-1}}\right]\mathbb E^{\mathbb P^{'}}_t\left[\left(\widetilde K_T^{\mathbb P^{'}}-\widetilde K_t^{\mathbb P^{'}}\right)^2\right]\right)^{\frac{p-1}{2p}}\\
&\leq C\left(\underset{\mathbb P^{'}\in\mathcal P^{\kappa}_H(t^+,\mathbb P)}{\esup^\mathbb P}\mathbb E^{\mathbb P^{'}}\left[\left(\widetilde K_T^{\mathbb P^{'}}-\widetilde K_t^{\mathbb P^{'}}\right)^2\right]\right)^{\frac{p-1}{2p}}\left(V_t-y_t^{\mathbb P^{'}}\right)^{1/p},\ \mathbb P-a.s.,
\end{align*}
where we used \reff{eq:momentsexp}. Arguing as in Step $\rm{(iii)}$ of the proof of Theorem \ref{uniqueref}, the above inequality along with Proposition \ref{prop.repref} shows that we have
$$\underset{\mathbb P^{'}\in\mathcal P^{\kappa}_H(t^+,\mathbb P)}{\einf^\mathbb P}\mathbb E^{\mathbb P^{'}}\left[\widetilde K_T^{\mathbb P^{'}}-\widetilde K_t^{\mathbb P^{'}}\right]=0,\ \mathbb P-a.s.,$$
that is to say that the minimum condition \ref{2bsdej.minK} is satisfied. 

\vspace{0.3em}
Finally, when the terminal condition is in $\Lc^{2,\kappa}_H$, it suffices to approximate it by a sequence $(\xi_n)_{n\geq 0}\subset {\rm UC_b}(\Omega)$, and to pass to the limit using the a priori estimates obtained in \cite{kpz3}. The proof is similar to step (ii) of the proof of Theorem 4.6 in \cite{stz}, we therefore omit it.
\ep

\section{Fully non-linear PIDEs}\label{sec.PIDE}

\subsection{Markovian 2BSDEJs}
In this section, we specialize our discussion on 2BSDEJs by considering the so-called Markovian case for which the generator $F$ and the function $H$ have a specified deterministic form
$$H_t(\omega, y,z,u,\gamma,v)=:h(t,B_t(\omega),y,z,u,\gamma,v),\text{ and } F_t(\omega, y,z,u,a,\nu)=: f(t,B_t(\omega),y,z,u,a,\nu),$$
for some function $f:[0,T]\times\R^d\times\R\times\R^d\times \hat L^2\times\mathbb S^{>0}_d\times\Nc\rightarrow \R$ and some function $h:[0,T]\times\R^d\times\R\times\R^d\times \hat L^2\times D_1\times D_2$. We also denote by $D^1_{f(t,x,y,z,u)}$ the domain of $F$ in $a$ and by $D^2_{f(t,x,y,z,u)}$ the domain of $f$ in $\nu$, for a fixed $(t,x,y,z,u)$

\vspace{0.5em}
Moreover, we also need to define Markovian counterparts of the set $\Vc$ of predictable compensators, and of the sets $\Rc_F$ of predictable functions. More precisely, we let $\Vc^m$ be the set of measures in $\Nc$ which does not depend on $t$ and $\omega$ (i.e. the L\'evy measures in $\Nc$), and for any $F\in\Vc^m$, let $\Rc_F^m$ to be the set of functions $\beta:E\rightarrow E$ which verify
$$\abs{\beta}(x)\leq C(1\wedge\abs{x}), \ F(dx)-a.e.$$
and such that $x\longmapsto \beta(x)$ is strictly monotone on the support of $F$. Finally, we let $\mathfrak V^m$ to be the set of measures $\widetilde F:=F\circ\left(\beta\right)^{-1}$ for some $F\in\Vc^m$ and some $\beta\in\Rc^m_F$.

\vspace{0.5em}
We can now define the Legendre-Fenchel transform of the generator $f$ as follows
$$\hat h(t,x,y ,z,u,\gamma,v):=\underset{(a,\nu)\in\mathbb S^{>0}_d\times\mathfrak V^m}{\sup}\left\{\frac12{\rm Tr}(a\gamma)+\int_E(Av)(x,e)\nu(de)-f(t,x,y,z,u,a,\nu)\right\},$$
for any $(t,x,y,z,u,\gamma,v)\in[0,T]\times \R^d\times\R\times\R^d\times \hat L^2\times C^2_b(E)$, where $C^2_b(E)$ denotes the set of functions from $E$ to $E$ which are $C^2$ with a bounded gradient and Hessian, which we endow with the topology of uniform convergence on compact sets, and where we recall that the non-local operator $A$ is defined in \reff{def_v_nu}.

\vspace{0.5em}
For simplicity, we abuse notations and let $\mathcal P^\kappa_h:=\mathcal P^\kappa_H$, as well as $\mathcal P^{\kappa,t}_h:=\mathcal P^{\kappa,t}_H$ for any $t\in[0,T]$. Assumptions \ref{assump.href} and \ref{assump.h2ref} now take the following (stronger) form
\begin{Assumption} \label{assump.href3}
\rm{(i)} $\Pc^\kappa_h$ is not empty and the domains $D^1_{f(t,x,y,z,u)}=D^1_{f(t)}$ and $D^2_{f(t,y,z,u)}=D^2_{f(t)}$ are independent of $(x,y,z,u)$.

\vspace{0.4em}
\rm{(ii)} The following uniform Lipschitz-type property holds. For all $(y,y',z,z',u,t,a,\nu,x)$
\begin{align*}
 &\abs{ f(t,x,y,z,u,a,\nu)- f(t,x,y',z',u,a,\nu)}\leq C\left(\abs{y-y'}+\abs{ a^{1/2}\left(z-z'\right)}\right).
\end{align*}

\vspace{0.4em}
(iii) The map $t\longmapsto f(t,x,y,z,u,a,\nu)$ has left-limits and is uniformly continuous from the right, uniformly in $(a,\nu)\in D^1_{f(t)}\times D^2_{f(t)}$.

\vspace{0.4em}
\rm{(iv)} For all $(t,x,y,z,u^1,u^{2},a,\nu)$, there exist two functions $\gamma$ and $\gamma'$ such that
\begin{align*}
 \int_{E}\delta^{1,2} u(x)\gamma'_t(x)\nu(dx)\leq f(t,x,y,z,u^1,a,\nu)- f(t,x,y,z,u^2,a,\nu)&\leq \int_{E}\delta^{1,2} u(x)\gamma_t(x)\nu(dx),
\end{align*}
where $\delta^{1,2} u:=u^1-u^2$ and $c_1(1\wedge \abs{x}) \leq \gamma_t(x) \leq c_2(1\wedge \abs{x})$ with $-1+\delta\leq c_1\leq0, \; c_2\geq 0,$
and $c_1'(1\wedge \abs{x}) \leq \gamma'_t(x) \leq c_2'(1\wedge \abs{x})$ with $-1+\delta\leq c_1'\leq0, \; c_2'\geq 0,$
for some $\delta >0$.

\vspace{0.4em}
\rm{(v)} $F$ is uniformly continuous in $x$, that is to say that there exists some modulus of continuity $\rho$ such that for all $(t,x,x',y,z,u,a,\nu)$
$$\abs{f(t,x,y,z,u,a,\nu)-f(t,x',y,z,u,a,\nu)}\leq \rho\left(\abs{x-x'}\right).$$
\end{Assumption}

\begin{Remark}\label{rem.linear}
With the exception of $(iii)$, the above assumptions are simple restatements of Assumptions \ref{assump.href} and \ref{assump.h2ref}. We emphasize that $(iii)$ above will be important in the proof of the Feynman-Kac representation formula. Notice also that since $\R^d$ is a convex space, it is always possible to choose the modulus $\rho$ to be both concave and with linear growth.
\end{Remark}

We consider a given Borel-measurable function $g:\R^d\longrightarrow\R$. The rest of this section will be devoted to relationships existing between the following 2BSDEJ
\begin{equation}\label{2bsdejmark}
Y_t=g(B_T)-\int_t^Tf\left(s,B_s,Y_s,Z_s,U_s,\widehat a_s,\nu^\P_s\right)ds-\int_t^TZ_sdB_s+K_T^\P-K_t^\P,\ \Pc^\kappa_h-a.s.,
\end{equation}
and the following fully non-linear PIDE
\begin{equation}\label{pide}
\begin{cases}
-\partial_tu(t,x)-\hat h(t,x,u(t,x),Du(t,x),\Kc u(t,x,\cdot),D^2u(t,x),u(t,x+\cdot))=0,\ (t,x)\in[0,T)\times\R^d\\
u(T,x)=g(x),\ x\in\R^d,
\end{cases}
\end{equation}
where $y\longmapsto \Kc v(t,x,y)$ is a function from $E$ to $E$ defined by
$$\Kc v(t,x,y):=v(t,x+y)-v(t^-,x).$$

\subsection{Smooth solutions of \reff{pide} and Feynman-Kac formula}
We start by showing that a smooth solution to \reff{pide} provides a solution to the 2BSDEJ \reff{2bsdejmark}. As was already showed in \cite{stz}, and unlike what happens for classical Markovian BSDEJs, the proof is not a simple application of It\^o's formula. Indeed, verifying the minimal condition \reff{2bsdej.minK} creates unavoidable technical difficulties.
\begin{Theorem}
Let Assumption \ref{assump.href3} hold and assume in addition that $h$ is continuous $($where continuity with respect to its fifth and seventh argument are to be understood w.r.t. the topology of uniform convergence on compact sets$)$, that $D^1_{f}:=D^1_{f(t)}$ and $D^2_{f}:=D^2_{f(t)}$ are actually independent of $t$, that $D^1_{f}$ is bounded from above and away from $0$, that $D^2_f$ is such that
\begin{equation}\label{condition}\underset{\nu\in D^2_{f}\cap\mathfrak V^m}{\sup}\int_E\left(\abs{x}^2{\bf 1}_{\abs{x}<1}+\abs{x}{\bf 1}_{\abs{x}\geq 1}\right)\nu(dx)<+\infty,\end{equation}
 and that $g$ is bounded. Let $u\in C^{1,2}([0,T),\R^d)$ be a classical solution of \reff{pide} with bounded gradient and Hessian, such that in addition $$\{(u,Du,\Kc u)(t,B_{t^-}), \ t\in[0,T]\}\in\D^{2,\kappa}_H\times\H^{2,\kappa}_H\times\mathbb J^{2,\kappa}_H.$$ Then, if we define
$$Y_t:=u(t,B_t), \;ÊZ_t:=Du(t,B_{t^-}),\ U_t(\cdot):=\Kc u(t,B_{t^-},\cdot),\ \Gamma_t:=D^2u(t,B_{t^-}),\ K_t^\P:=\int_0^tk_s^\P ds,$$
\vspace{-0.3em}
$$k_t^\P:=\hat h (t, B_{t^-},Y_t,Z_t, U_t,\Gamma_t,u(t,B_{t^-}+\cdot))-\frac12{\rm Tr}[\widehat a_t\Gamma_t]-\int_E(Au)(B_{t^-},x)\nu^\P_t(dx)+f(t,B_t,Y_t,Z_t,U_t,\widehat a_t,\nu^\P_t),$$
$(Y,Z,U)$ is the unique solution of \reff{2bsdejmark}.
\end{Theorem}
\begin{Remark}
The condition \reff{condition} above seems difficult to avoid with our approach here, and basically demands uniform moments for the small and large jumps of the canonical process over the whole uncertainty set. Moreover, we would like to point out that this condition also appears in the recent work \cite{nn3}, when the authors look at viscosity solution to PIDE \reff{pide} when $f=0$.
\end{Remark}
\proof
A simple application of It\^o's formula shows that $(Y,Z,U)$ does satisfy the equation \reff{2bsdej}. Since in addition $Y_T=g(B_T)\in\L^{2,\kappa}_H$ (because $g$ is bounded), it only remains to verify that for any $\P\in\Pc^\kappa_h$ and for any $t\in[0,T]$
$$\underset{\P^{'}\in\Pc^\kappa_h(t^+,\P)}{\einf^\P}\E^{\P^{'}}_t\left[\int_t^Tk_s^{\P^{'}}ds\right]=0,\ \P-a.s.$$   

Towards this goal, we follow the proof of Theorem $5.3$ in \cite{stz} and we adapt it to our jump framework. Let us outline the proof for the sake of clarity. The main idea is, as is common in stochastic control problems, to find an $\eps$-optimal control in the definition of $\hat h$, which means here both a volatility process and a jump measure. The main difficulty after that is to be able to find a probability measure in $\Pc^\kappa_h$ such that the characteristics of the canonical process $B$ under this measure coincide with the $\eps$-optimal controls. We emphasize that even though it may be possible to find such a measure in the larger set $\overline{\Pc}_W$ (and even in this case it may prove impossible, see Remark 2.3 in \cite{stz}), it is not clear at all that this measure will be in $\overline{\Pc}_S$, and thus in $\Pc^\kappa_h$.

\vspace{0.4em}
We now start the proof. By a classical measurable selection argument, for any $\eps>0$, we can find a predictable process $a^\eps$ taking values in $D^1_f$ and a predictable random measure $\nu^\eps$, taking values in $D^2_f\cap\mathfrak V^m$ (i.e. there exist $(F^\eps,\beta^\eps)\in\Vc^m\times\Rc^m_F$ such that $\nu^\eps=F^\eps\circ (\beta^\eps)^{-1}$) with
$$ \hat h (t, B_{t^-},Y_t,Z_t, U_t,\Gamma_t,u(t,B_{t^-}+\cdot))\leq \frac12{\rm Tr}[a^\eps_t\Gamma_t]+\int_E(Au)(B_{t^-},x)\nu^\eps_t(dx)-f(t,B_t,Y_t,Z_t,U_t,a^\eps_t,\nu^\eps_t)+\eps.$$

Fix now some $\P:=\P^{\alpha,\beta}_F\in\Pc^\kappa_H$ and some $t\in[0,T]$. We will now show that we can find some $(\alpha^\eps,b^\eps,\widetilde F^\eps)$ such that $\P^{\alpha^\eps,b^\eps}_{\widetilde F^\eps}\in\Pc^\kappa_H(t^+,\P)$ and for $s\in[t,T]$
$$\widehat a_s=a^\eps,\ \nu^{\P^{\alpha^\eps,b^\eps}_{\widetilde F^\eps}}_s=\nu^\eps,\ ds\times d\P^{\alpha^\eps,b^\eps}_{\widetilde F^\eps}-a.e.$$

For notational simplicity, let us define for $s\geq r\geq t$
\begin{align*}
X_s^r:=&\ h(s,B_{s^-},Y_s,Z_s, U_s,\Gamma_s,u(s,B_{s^-}+\cdot))- \frac12{\rm Tr}[a^\eps_r\Gamma_s]-\int_E(Au)(B_{s^-},x)\nu^\eps_r(dx)\\
&+f(s,B_s,Y_s,Z_s,U_s,a^\eps_r,\nu^\eps_r).
\end{align*}
We next define a sequence of $\F$-stopping times. Let
\begin{align*}
\tau_0^\eps:&=\inf\Big\{s\geq t,\ X_s^t\geq 2\eps\text{ or }X_{s^-}^t\geq 2\eps\Big\}\wedge T\\
\tau_{n+1}^\eps:&=\inf\Big\{s\geq \tau_n^\eps,\ X_s^{\tau_n^\eps}\geq 2\eps\text{ or }X_{s^-}^{\tau_n^\eps}\geq 2\eps\Big\}\wedge T,\ n\geq 0.
\end{align*}
Notice that since by definition $X^s_s\leq \eps$ for any $s\geq t$, we always have $\tau_{n+1}^\eps>\tau_n^\eps$, for any $n\geq 0$, and $\tau_0^\eps>t$. Besides, since $B,Y,Z,U,\Gamma, u$ are all c\`adl\`ag, it is a classical result that the $\tau^\eps_n$ are indeed $\F$-stopping times.

\vspace{0.5em}
Next, for any $\omega\in\Omega$, the maps $t\mapsto \hat h(t,x,y,z,u,\gamma,v)$ and $t\mapsto f(t,x,y,z,u,a,\nu)$ are respectively uniformly continuous from the right (see \cite{applebaum}, Appendix $2.8$ for more details) and uniformly continuous from the right uniformly in $(a,\nu)\in D^1_f\times D^2_f$ (since they are c\`adl\`ag on the compact $[0,T]$). Moreover, since we assumed that $D^1_f$ was bounded from above and away from $0$, that 
$$\underset{\nu\in D^2_{f}\cap\mathfrak V^m}{\sup}\int_E\left((1\wedge \abs{x}^2)+\abs{x}{\bf 1}_{\abs{x}\geq 1}\right)\nu(dx)<+\infty,$$
and that $Du$ and $D^2u$ were bounded, we can deduce that for any omega, the function
 \begin{align*}
 h(s,B_{s^-}(\omega),Y_s(\omega),Z_s(\omega), U_s(\omega),\Gamma_s(\omega),u(s,B_{s^-}(\omega)+\cdot))- \frac12{\rm Tr}[a\Gamma_s(\omega)]-\int_E(Au)(B_{s^-}(\omega),x)\nu(dx)\\+f(s,B_s(\omega),Y_s(\omega),Z_s(\omega),U_s(\omega),a,\nu),
 \end{align*}
 is uniformly continuous from the right in $s$, uniformly in $(a,\nu)$. This implies that the $\tau_n^\eps$ cannot accumulate and that it is possible to find some $\delta(\eps,\omega)>$, independent of $n$, such that $\tau_{n+1}^\eps(\omega)\geq\tau_n^\eps(\omega)+\delta(\eps,\omega)$. In particular, this also implies that there exists some finite $N\in\N$ such that $\tau_n^\eps=T$ for any $n\geq N$.
 
 \vspace{0.5em}
 Let us now define for $s\geq \tau_0^\eps$
 $$\tilde a^\eps_s:=\sum_{n=0}^{+\infty}a^\eps_{\tau^\eps_n}{\bf 1}_{s\in[\tau_n^\eps,\tau_{n+1}^\eps)},\ \tilde F^\eps_s:=F_t{\bf 1}_{0\leq s\leq \tau_0^\eps}+{\bf 1}_{s\geq \tau_0^\eps}\sum_{n=0}^{+\infty}F^\eps_{\tau^\eps_n}{\bf1}_{s\in[\tau_n^\eps,\tau_{n+1}^\eps)},\  \tilde \beta^\eps_s:=\sum_{n=0}^{+\infty}\beta^\eps_{\tau^\eps_n}{\bf 1}_{s\in[\tau_n^\eps,\tau_{n+1}^\eps)}.$$
 Consider next the following SDE on $[\tau_0^\eps,T]$
 \begin{equation}\label{sde}
 dZ_s= \left(a^\eps_s(Z_\cdot)\right)^{1/2}dB^{\P_{0,\tilde F^\eps},c}_s+\int_E\tilde\beta_s(Z_\cdot,x)\left(\mu_B(dx,ds)-\tilde F^\eps_s(dx)ds\right),\ \P_{0,\tilde F^\eps}-a.s.
 \end{equation}
 
 It is proved in Lemma \ref{lemma.sde} that the above SDE has a unique strong solution $Z^\eps$ on $[\tau_0^\eps,T]$ such that $Z^\eps_{\tau_0^\eps}=0.$ Define then
 $$\alpha^\eps_t:=\alpha_t{\bf 1}_{0\leq t\leq \tau_0^\eps}+{\bf 1}_{t\geq \tau_0^\eps}\tilde a^\eps(Z^\eps_\cdot),\ b^\eps_t(x):=\beta_t(x){\bf 1}_{0\leq t\leq \tau_0^\eps}+{\bf 1}_{t\geq \tau_0^\eps}\tilde \beta^\eps(Z^\eps_\cdot,x).$$
 It is immediate to verify that $\alpha^\eps\in\mathcal D$, $\tilde F^\eps\in\nu$ and $b^\eps\in\mathcal R_{\tilde F^\eps}$ (see the proof of Lemma A.3 in \cite{kpz3} for similar arguments). We can therefore define the probability measure $\P_{\tilde F^\eps}^{\alpha^\eps,b^\eps}\in\mathcal P^\kappa_h$, which by definition, coincides with $\P$ on $\mathcal F_{t^+}$ (since $\tau_0^\eps>t$). Using $(2.6)$ and $(2.7)$ in \cite{kpz3}, we then deduce that
 $$\widehat a_s= \tilde a^\eps_s,\ \nu^{\P^{\alpha^\eps,b^\eps}_{\tilde F^\eps}}_s=\nu_s^{\tilde F^\eps,\tilde\beta^\eps},\ \P^{\alpha^\eps,b^\eps}_{\tilde F^\eps}-a.s.\text{ on }[\tau_0^\eps,T].$$
 
 This implies that, $\P^{\alpha^\eps,b^\eps}_{\tilde F^\eps}-a.s$,
\begin{align*}
 \hat h (s, B_{s^-},Y_s,Z_s, U_s,\Gamma_s,u(s,B_{s^-}+\cdot))\leq &\ \frac12{\rm Tr}[\widehat a_s\Gamma_s]+\int_E(Au)(B_{s^-},x)\nu^{\P^{\alpha^\eps,b^\eps}_{\tilde F^\eps}}_s(dx)\\
 &-f\left(s,B_s,Y_s,Z_s,U_s,\widehat a_s,\nu^{\P^{\alpha^\eps,b^\eps}_{\tilde F^\eps}}_s\right)+\eps,\text{ for a.e. $s\in[\tau_0^\eps,T]$}.
 \end{align*}
 Finally,
 
 $$\underset{\P^{'}\in\Pc^\kappa_h(t^+,\P)}{\einf^\P}\E^{\P^{'}}_t\left[\int_t^Tk_s^{\P^{'}}ds\right]\leq 2\eps(T-t)+\E^{\mathbb P^{\alpha^\eps,b^\eps}_{\tilde F^\eps}}\left[\int_{\tau_0^\eps}^Tk_s^{\mathbb P^{\alpha^\eps,b^\eps}_{\tilde F^\eps}}ds\right]\leq 4\eps(T-t).$$
Since $\eps>0$ was arbitrary, this ends the proof. 
\ep

\subsection{2BSDEJs and viscosity solutions to fully non-linear PIDEs}
\subsubsection{Time-space regularity of Markovian solutions to 2BSDEJs}
In this section, we specialize the discussion and notations of Section \ref{sec.existtt} to the Markovian framework and obtain additional regularity results. 

\vspace{0.3em}
For simplicity, let us denote
$$B_s^{t,x}:=x+B^t_s,\ \text{for all $(t,s,x)\in[0,T]\times[t,T]\times\R^d,$}$$
and for any $(t,x)\in[0,T]\times\R^d$, any $\F^t$-stopping time $\tau$, any $\P\in\Pc^{\kappa,t}_h$, and r.v. $\eta\in L^2(\P)$ which is $\Fc^t_\tau$-measurable, we let $(\Yc^{\P,t,x},\Zc^{\P,t,x},\Uc^{\P,t,x}):=(\Yc^{\P,t,x}(\tau,\eta),\Zc^{\P,t,x}(\tau,\eta),\Uc^{\P,t,x}(\tau,\eta))$ be the unique solution to the following BSDEJ
\begin{align}\label{bsdejmark}
\nonumber \Yc^{\mathbb P,t,x}_{s}&=\eta+\int^{\tau}_{s}f\left(r,B_r^{t,x},\Yc^{\mathbb P,t,x}_{r},\Zc^{\mathbb P,t,x}_{r},\Uc_r^{\mathbb P,t,x},\widehat a^t_r,\nu^{t,\P}_r\right)dr-\int^{\tau}_{s}\Zc^{\mathbb P,t,x}_{r}dB^{t,\mathbb P,c}_{r}\\
&\hspace{0.9em}-\int_s^{\tau}\int_{E}\Uc_r^{\mathbb P,t,x}(e)\widetilde\mu^\mathbb P_{B^{t}}(de,dr), \ \mathbb P-a.s.,\ s\in[t,\tau].
\end{align}

Then, exactly as the process $V$ defined in \reff{sol}, we consider the value function
$$u(t,x):=\underset{\P\in\Pc^{t,\kappa}_h}{\sup}\Yc^{\P,t,x}_t\left(T,g(B_T^{t,x})\right),\ (t,x)\in[0,T]\times \R^d,$$
which is indeed deterministic due to the Blumenthal $0-1$ law, which holds true for any $\P\in\Pc^{t,\kappa}_h$.

Before stating the next result, we need to consider in addition the following assumption.
\begin{Assumption}\label{assump.markovian}
The map $x\longmapsto g(x)$ has linear growth and 
\begin{align}
\Lambda(t,x):=\underset{\P\in\Pc^{t,\kappa}_h}{\sup}\left(\E^\P\left[\abs{g\left(B^{t,x}_T\right)}+\int_t^T\abs{f\left(s,B_s^{t,x},0,0,0,\widehat a^t_s,\nu^{t,\P}_s\right)}^\kappa\right]\right)^{\frac1\kappa}, \label{grandlambda.def}
\end{align}
is such that
$$\underset{\P\in\Pc^{t,\kappa}_h}{\sup}\E^\P\left[\underset{t\leq s\leq T}{\sup}\left(\Lambda\left(s,B^{t,x}_s\right)\right)^2\right]<+\infty,\text{ for any $(t,x)\in[0,T]\times\R^d$}.$$
\end{Assumption}

\begin{Remark}
We emphasize that we could make the weaker assumption that $g$ has polynomial growth, say of integer order $p$ but, to compensate this, we would need the additional assumption that the compensators $\nu$ that we consider have moments of order $p$, thus reducing the set of probability measures we allow for. The wellposedness of our 2BSDEJ \ref{2bsdej} would still hold. For the sake of simplicity we will directly ask that $g$ has linear growth.
\end{Remark}

We refer to Remark 5.8 in \cite{stz} for sufficient conditions ensuring that Assumption \ref{assump.markovian} holds true. The next result generalizes Theorem 5.9 and Proposition 5.10 of \cite{stz}.
\begin{Proposition}\label{prop.reg}
We have the following results
\begin{itemize}
\item[{\rm(i)}] Let Assumptions \ref{assump.href3} and \ref{assump.markovian} hold, and assume furthermore that $g$ is uniformly continuous. Then the 2BSDEJ \reff{2bsdej} has a unique solution $(Y,Z,U)\in\D^{2,\kappa}_H\times\H^{2,\kappa}_H\times\mathbb J^{2,\kappa}_H$. Moreover, we have the identity $Y_t=u(t,B_t)$ and the function $u$ is uniformly continuous in $x$, uniformly in $t$, and right-continuous in $t$.

\item[{\rm(ii)}] Let Assumptions \ref{assump.href3} and \ref{assump.markovian} hold, and assume furthermore that $g$ is lower semi-continuous. Then $u$ is lower semicontinuous in $(t,x)$, from the right in $t$, that is to say that for any $(t,x)\in[0,T]\times \R^d$ and any sequence $(t_n,x_n)_{n\geq 0}$ such that
$$(t_n,x_n)\underset{n\rightarrow +\infty}{\longrightarrow}(t,x)\text{ and } t_n> t,\text{ for all $n\geq 0$},$$
we have $\underset{n\rightarrow +\infty}{\underline{\lim}}u(t_n,x_n)\geq u(t,x).$
\end{itemize}
\end{Proposition}

In order to prove this Proposition, we will need the following weak dynamic programming property, in the spirit of the work by Bouchard and Touzi \cite{bt}. It is very closely related to the proof of the dynamic programming property of Proposition \ref{progdyn}, with the additional difficulty that less regularity is assumed on $g$. Moreover, its proof is very close to the proofs of Proposition 5.14 and Lemmas 6.2 and 6.4 in \cite{stz}. Hence, we will only sketch some parts of its proof, which is relegated to the appendix.

\begin{Lemma}\label{partial.dpp}
Under assumptions \ref{assump.href3} and \ref{assump.markovian}, for any family of $\F^t$-stopping times $\{\tau^{\P}, \mathbb P\in\mathcal P^{t,\kappa}_h\}$:
\begin{align}
u(t,x) \leq \underset{\mathbb P\in\mathcal P^{t,\kappa}_h}{\sup} \mathcal Y^{\mathbb P,t,x}_t\left(\tau^{\P},X\right), \text{ for all } \left(t,x\right)\in\left[0,T\right]\times\R^d \label{partial.dpp1}
\end{align}
and for any $\mathcal{F}^t_{\tau^{\P}}$-measurable r.v. $X$ such that $X\geq u(\tau^{\P},B_{\tau^{\P}}^{t,x})$, $\P$-a.s. Moreover, when the function $g$ is lower semi-continuous, 
\begin{align}
u(t,x) = \underset{\mathbb P\in\mathcal P^{t,\kappa}_h}{\sup} \mathcal Y^{\mathbb P,t,x}_t\left(\tau^{\P},u(\tau^{\P},B_{\tau^{\P}}^{t,x})\right), \text{ for all } \left(t,x\right)\in\left[0,T\right]\times\R^d, \label{partial.dpp2}
\end{align}
\end{Lemma}

\vspace{0.3em}
\proof [of Proposition \ref{prop.reg}.]
Proposition \ref{prop.reg} (i) can be proved exactly as in \cite{stz}, using Theorem \ref{mainref}, Lemmas \ref{unifcont} and \ref{lem.cadlag} and Proposition \ref{prop.repref}.

\vspace{0.3em}
(ii) We follow \cite{stz}. We start by introducing the functional
$$\mathfrak J(\P,t,x):=\mathbb E^\P\left[\mathfrak y_t^\P\left(t,x\right))\right],\ (\P,t,x)\in \Pc^\kappa_h\times[0,T]\times\R^d,$$
where $\mathfrak y ^{\P}(t,x)$ is the first component of the solution to the BSDEJ under $\P$ with terminal condition $g(x+B_T-B_t)$ and generator $f(s,x+B_s-B_t,y,z,u,a,\nu)$. The first step of the proof is to show the following identity
\begin{equation}\label{identity}
u(t,x)=\underset{\P\in\Pc^\kappa_h}{\sup}\mathfrak J(\P,t,x).
\end{equation}

First, by \reff{eq.picard}, we have for any $\P\in\mathcal P^\kappa_h$ and for $\P-a.e.$ $\omega\in\Omega$
$$\mathfrak y_t^\P(t,x)(\omega)=\Yc_t^{\P^{t,\omega},t,x}\left(T,g\left(B_T^{t,x}\right)\right)\leq u(t,x).$$
This implies, that $\mathfrak J(\P,t,x)\leq u(t,x)$. The other inequality can be proved exactly as in the proof of Proposition 5.10 in \cite{stz}. Then, it is clearly sufficient to show that the map $(t,x)\longmapsto J(\P,t,x)$ is lower semi-continuous, from the right in $t$, for any $\P\in\mathcal P^\kappa_h$. 

\vspace{0.3em}
Consider thus some $(t,x)\in[0,T]\times \R^d$, some $\P\in\Pc^\kappa_h$ and a sequence $(t_n,x_n)_{n\geq 0}$ such that $(t_n,x_n)\rightarrow (t,x)$ and $t_n>t$ for any $n\geq 0$. Consider the following $\underline{\lim}$
$$\xi_n:=\underset{k\geq n}{\inf} g(x_k+B_T-B_{t_k}),\ f^{n,\P}(s,y,z,u):=\underset{k\geq n}{\inf}  f(s,x_k+B_s-B_{t_k},y,z,u,\widehat a_s,\nu^\P_s),$$
$$\overline\xi:=\underset{n\rightarrow+\infty}{\lim}\xi_n,\ \overline{f}^\P:=\underset{n\rightarrow+\infty}{\lim}f^{n,\P},$$
and let $(\Yc^n,\Zc^n,\Uc^n)$ denote the solution to the BSDEJ under $\P$ with terminal condition $\xi_n$ and generator $f^{n,\P}$. Since $g$ and the modulus of uniform continuity of $f$ have linear growth in $x$ (remember Remark \reff{rem.linear}), since $f$ is uniformly Lipschitz continuous, and since under any of the measure considered the canonical process $B$ is a square-integrable martingale (see Definition \ref{def}), it can be checked directly that this BSDEJ has a unique solution, and by stability for BSDEJ that 
$$\underset{n\rightarrow +\infty}{\lim}\E^\P[\Yc^n_t]=\E^\P[\overline{\Yc}_t],$$
where $\overline{\Yc}$ denotes the first component of the solution of the BSDEJ with terminal condition $\overline\xi$ and generator $\overline{f}^\P$.

\vspace{0.3em}
Since $g$ is lower semi-continuous, $f$ is uniformly continuous in $x$, $B$ is c\`adl\`ag and $t_n>t$ for any $n\geq 0$, we deduce that
$$\overline\xi\geq g(x+B_T-B_t),\text{and }\overline f^\P(s,y,z,u)=f(s,x+B_s-B_t,y,z,u,\widehat a_s,v_s^\P).$$
Hence
\begin{align*}
\underset{n\rightarrow +\infty}{\underline \lim} \mathfrak J(\P,t_n,x_n)\geq \underset{n\rightarrow +\infty}{\lim}\E^\P[\Yc^n_t]=\E^\P[\overline{\Yc}_t]\geq \E^\P[\mathfrak y^\P_t(t,x)]=\mathfrak J(\P,t,x),
\end{align*}
which ends the proof.
\ep

\subsubsection{Viscosity solution of PIDE \reff{pide}}
\begin{Definition}
$(i)$ A bounded lower-semicontinuous function $v$ is called a viscosity super-solu-tion of the PIDE \reff{pide} if $v(T,\cdot)\geq g(\cdot)$ and if for every function $\varphi\in C^{3}_b([0,T)\times\R^d)$ $($i.e. the space of functions from $[0,T)\times\R^d$ which are thrice continuously differentiable with bounded derivatives$)$ such that 
$$0=v(t_0,x_0)-\varphi(t_0,x_0)=\underset{(t,x)\in[0,T)\times \R^d}{\min}(v-\varphi)(t,x),$$
we have
$$-\partial_tv(t_0,x_0)-\hat h(t_0,x_0,\varphi(t_0,x_0),D\varphi(t_0,x_0),\Kc \varphi(t_0,x_0,\cdot),D^2\varphi(t_0,x_0),\varphi(t_0,x_0+\cdot))\geq 0.$$
$(ii)$ A bounded upper-semicontinuous function $v$ is called a viscosity sub-solution of the PIDE \reff{pide} if $v(T,\cdot)\leq g(\cdot)$ and if for every function $\varphi\in C^{3}_b([0,T)\times\R^d)$ such that 
$$0=v(t_0,x_0)-\varphi(t_0,x_0)=\underset{(t,x)\in[0,T)\times \R^d}{\max}(v-\varphi)(t,x),$$
we have
$$-\partial_tv(t_0,x_0)-\hat h(t_0,x_0,\varphi(t_0,x_0),D\varphi(t_0,x_0),\Kc \varphi(t_0,x_0,\cdot),D^2\varphi(t_0,x_0),\varphi(t_0,x_0+\cdot))\leq 0.$$
$(iii)$ A continuous function $v$ is a viscosity solution of \reff{pide} if it is both a viscosity sub and supersolution.
\end{Definition}

This section is devoted to the proof of the following result, which generalizes to the case of 2BSDEJs Proposition 5.4 of \cite{nn3}, which considers the case $f=0$. The arguments are classical, as soon as one has at disposition a dynamic programming principle. We however give a detailed proof, since the presence of the non-linearity $f$ complicates the estimates. 
\begin{Theorem}\label{th:visco}
Let Assumption \ref{assump.href3} hold and assume in addition that $h$ is continuous, that $t\longmapsto f(t,x,y,z,a,\nu)$ is uniformly continuous from the right, uniformly in all the other variables, that $D^1_{f}:=D^1_{f(t)}$ and $D^2_{f}:=D^2_{f(t)}$ are actually independent of $t$, that $D^1_{f}$ is bounded from above and away from $0$, that $D^2_f$ is such that
$$\underset{\nu\in D^2_{f}\cap\mathfrak V^m}{\sup}\int_E\left(\abs{x}^2+\abs{x}{\bf 1}_{\abs{x}\geq 1}\right)\nu(dx)<+\infty.$$
 and that $g$ is uniformly continuous and bounded. Then, $u$ is a viscosity solution of the PIDE \reff{pide}.
\end{Theorem}
\begin{Remark}
We would to point out two differences with \cite{nn3}. First, the set of test functions that we consider is not the same, since we assume more regularity. However, as is well-known in viscosity theory, this is actually without loss of generality by simple density arguments. Then, the integrability assumptions on the jumps of the canonical process under the measures considered is not the same. Indeed, we assume a uniform control for the second moment of both the small and large jumps, while \cite{nn3} only does it for the small jumps. This seems unavoidable in our setting since we want the canonical process to be square-integrable under every measures, and we want to have a uniform control on its norm. Nonetheless, the added assumption allows us to get rid off the assumption of \cite{nn3} on the limit as $\epsilon$ goes to $0$ of the first order moment of small jumps $($see their condition $(5.2))$. 
\end{Remark}

\begin{Remark}
Of course, Theorem \ref{th:visco} should be complemented with a comparison theorem which would then imply uniqueness of viscosity solutions to \reff{pide}. Given the length of this paper, we will refrain from studying this problem here, and we refer instead the reader to Proposition $5.5$ in \cite{nn3} and the references therein for examples of assumptions under which such a result holds.
\end{Remark}
Before proving this theorem, we will need the following lemmas, which notably insure that the function $u$ is jointly continuous, which is needed if we want to prove that it is a (continuous) viscosity solution of the PIDE \reff{pide}.
\begin{Lemma}\label{lemma:estimates}
Let the assumptions of Theorem \reff{th:visco} hold. Then, for some constant $C>0$, we have that for any $(t,t')\in[0,T]^2$
$$\underset{\P\in\Pc^{t',\kappa}_h}{\sup}\E^\P\left[\underset{t'\leq s\leq t}{\sup}\abs{B_s^{t'}}^2\right]\leq C(t-t'),\ \underset{\P\in\Pc^{t',\kappa}_h}{\sup}\E^\P\left[\underset{t'\leq s\leq t}{\sup}\abs{B_s^{t'}}\right]\leq C\sqrt{\abs{t-t'}}.$$
\end{Lemma}
\proof
First of all, by BDG inequality, we have for any $\P\in\Pc^{t',\kappa}_h$ and for some constant $C>0$ which may vary from line to line
$$\E^\P\left[\underset{t'\leq s\leq t}{\sup}\abs{B_s^{t'}}^2\right]\leq C\E^\P\left[\int_{t'}^t\abs{\widehat a^{t'}_s}ds+\int_{t'}^t\int_E\abs{x}^2\nu^{t',\P}_s(dx)ds\right]\leq C(t-t'),$$
where we have used the fact that $D^1_f$ is bounded and that
$\underset{\nu\in D^2_{f}\cap\mathfrak V^m}{\sup}\int_E\abs{x}^2\nu(dx)$ is finite. The second term can be treated similarly.
\ep
\begin{Lemma}\label{lemma:unifcontt}
Let the assumptions of Theorem \reff{th:visco} hold. Then, the map $t\longmapsto u(t,x)$ is uniformly continuous. More precisely, if $\rho$ denotes the modulus of continuity in $x$ of $g$ and $f$, we have for some constant $C>0$
$$\abs{u(t,x)-u(t',x)}\leq C\left(\rho\left(C\sqrt{\abs{t-t'}}\right)^{\frac12}+(1+\abs{x})\sqrt{\abs{t-t'}}\right), \ (t,t',x)\in[0,T]^2\times \R^d.$$
\end{Lemma}
\proof
By \reff{identity}, we have
\begin{align*}
\abs{u(t,x)-u(t',x)}&\leq \underset{\P\in\Pc^\kappa_h}{\sup}\E^\P\left[\abs{\mathfrak y_t^\P(t,x)-\mathfrak y_{t'}^\P(t',x) }\right]\\
&\leq \underset{\P\in\Pc^\kappa_h}{\sup}\E^\P\left[\abs{\mathfrak y_t^\P(t,x)-\mathfrak y_{t}^\P(t',x) }\right]+\underset{\P\in\Pc^\kappa_h}{\sup}\E^\P\left[\abs{\mathfrak y_t^\P(t',x)-\mathfrak y_{t'}^\P(t',x) }\right].
\end{align*}
For the first term on the right-hand side, we have by classical linearization arguments for Lipschitz BSDEJs, following the same line as Step $2$ of the proof of Theorem $4.1$ in \cite{kpz3} that
$$\underset{\P\in\Pc^\kappa_h}{\sup}\E^\P\left[\abs{\mathfrak y_t^\P(t,x)-\mathfrak y_{t}^\P(t',x) }\right]\leq  C\underset{\P\in\Pc^\kappa_h}{\sup}\E^\P\left[\underset{0\leq s\leq T}{\sup}\abs{M_s}\rho\left(\abs{B_t-B_{t'}}\right)\right],$$
where $M$ is a process such that
$$\underset{\P\in\Pc^\kappa_h}{\sup}\E^\P\left[\underset{0\leq s\leq T}{\sup}\abs{M_s}^p\right]<+\infty, \ \text{for any $p\geq 1$.}$$
Using twice Cauchy-Schwarz inequality and the fact that $\rho$ has linear growth, we deduce that
\begin{align*}
\underset{\P\in\Pc^\kappa_h}{\sup}\E^\P\left[\abs{\mathfrak y_t^\P(t,x)-\mathfrak y_{t}^\P(t',x) }\right]\leq &\ \underset{\P\in\Pc^\kappa_h}{\sup}\E^\P\left[\underset{0\leq s\leq T}{\sup}\abs{M_s}^4\right]^{\frac14}C\left(1+\underset{\P\in\Pc^{\kappa}_h}{\sup}\E^\P\left[\underset{0\leq s\leq T}{\sup}\abs{B_s}^2\right]^{\frac14}\right)\\
&\times\underset{\P\in\Pc^\kappa_h}{\sup}\E^\P\left[\rho\left(\abs{B_t-B_{t'}}\right)\right]^{\frac12}.\\
\leq &\ C\rho\left(C\sqrt{\abs{t-t'}}\right)^{\frac12},
\end{align*}
where we used in the last line the fact the supremum of $M$ has moments of any order, Lemma \ref{lemma:estimates} as well as Jensen's inequality (remember that $\rho$ is concave).

\vspace{0.3em}
We then have, assuming w.l.o.g. that $t\leq t'$ and denoting by $(\mathfrak z^\P(t',x),\mathfrak u^\P(t',x))$ the second and third components of the solution of the BSDEJ associated to $\mathfrak y^\P(t',x)$
\begin{align*}
&\underset{\P\in\Pc^\kappa_h}{\sup}\E^\P\left[\abs{\mathfrak y_t^\P(t',x)-\mathfrak y_{t'}^\P(t',x) }\right]\\
&\leq \underset{\P\in\Pc^\kappa_h}{\sup}\E^\P\left[\int_t^{t'}\abs{f\left(s,x+B_s-B_{t'},\mathfrak y^\P(t',x),\mathfrak z^\P(t',x),\mathfrak u^\P(t',x),\widehat a_s,\nu^\P_s\right)}ds\right]\\
&\leq C\sqrt{t'-t}\left(1+\abs{x}+\underset{\P\in\Pc^\kappa_h}{\sup}\left\{\No{B}_{\D^2(\P)}^2+\No{\mathfrak y_s^\P(t',x)}_{\D^2(\P)}^2+\No{\mathfrak z_s^\P(t',x)}^2_{\H^2(\P)}+\No{\mathfrak u^\P(t',x)}_{\mathbb J^2(\P)}\right\}\right.\\
&\left.\hspace{0.9em}+\underset{\P\in\Pc^\kappa_h}{\sup}\mathbb E^\P\left[\int_0^T\abs{f(s,0,0,0,0,\widehat a_s,\nu^\P_s)}^2ds\right]\right)\\
&\leq C(1+\abs{x})\sqrt{t'-t},
\end{align*}
where we used the fact that $f$ is uniformly continuous in $x$, uniformly Lipschitz in $(y,z,u)$ and that $g$ is bounded and $f$ is sufficiently integrable. Hence the desired result.
\ep

\vspace{0.3em}
We can now proceed to the

\vspace{0.3em}
\proof[Proof of Theorem \ref{th:visco}]

{\bf Viscosity super-solution}: First of all, by Proposition \ref{prop.reg} and Lemma \ref{lemma:unifcontt}, we know that the map $(t,x)\longmapsto u(t,x)$ is continuous. Let us now prove the viscosity super-solution property. Let $\varphi\in C^{3}_b([0,T)\times\R^d)$ and $(t_0,x_0)\in[0,T)\times\R^d$ be such that
$$0=u(t_0,x_0)-\varphi(t_0,x_0)=\underset{(t,x)\in[0,T)\times \R^d}{\min}(u-\varphi)(t,x).$$
Fix some $\eta>0$ such that $t_0+\eta<T$. By \reff{partial.dpp2} with the constant family of stopping times $t_0+\eta$, we know that
\begin{align}\label{ineq:visco}
\nonumber\varphi(t_0,x_0)=u(t_0,x_0) &= \underset{\mathbb P\in\mathcal P^{t_0,\kappa}_h}{\sup} \mathcal Y^{\mathbb P,t_0,x_0}_{t_0}\left(t_0+\eta,u(t_0+\eta,B_{t_0+\eta}^{t_0,x_0})\right)\\
&\geq \underset{\mathbb P\in\mathcal P^{t_0,\kappa}_h}{\sup} \mathcal Y^{\mathbb P,t_0,x_0}_{t_0}\left(t_0+\eta,\varphi(t_0+\eta,B_{t_0+\eta}^{t_0,x_0})\right),
\end{align}
where we used the comparison theorem for BSDEJs in the last inequality.

\vspace{0.3em}
For notational simplicity, we denote by $(y^{\P,t_0,x_0,\varphi},z^{\P,t_0,x_0,\varphi},u^{\P,t_0,x_0,\varphi})$ the solution to the BSDEJ associated to $\mathcal Y^{\mathbb P,t_0,x_0}_{t_0}\left(t_0+\eta,\varphi(t_0+\eta,B_{t_0+\eta}^{t_0,x_0})\right)$. We also define for any $(\P,\omega,s,a,\nu)\in\Pc^{t_0,\kappa}_h\times\Omega\times[t_0,t_0+\eta]\times D^1_f\times (D^2_f\cap\mathfrak V^m)$
\begin{align*}
&\mathcal L^{t_0,x_0}\varphi(\omega,s,a,\nu):= \varphi_t(s,B_s^{t_0,x_0}(\omega))+\frac12\Tr{aD^2\varphi(s,B_s^{t_0,x_0}(\omega)}+\int_E(A\varphi)(B_{s^-}^{t_0,x_0}(\omega),e)\nu(de)\\
&-f\left(s,B_s^{t_0,x_0}(\omega),\varphi(s,B_s^{t_0,x_0}(\omega)),D\varphi(s,B_s^{t_0,x_0}(\omega)),\mathcal K\varphi(s,B_{s^-}^{t_0,x_0}(\omega),\cdot),a,\nu\right).
\end{align*}
Let us then define
$$\delta y_s^{\P,t_0,x_0,\varphi}:=y_s^{\P,t_0,x_0,\varphi}-\varphi(s,B_s^{t_0,x_0}),\ \delta z_s^{\P,t_0,x_0,\varphi}:=z_s^{\P,t_0,x_0,\varphi}-D\varphi(s,B_s^{t_0,x_0})$$
$$ \delta u_s^{\P,t_0,x_0,\varphi}(\cdot):=u_s^{\P,t_0,x_0,\varphi}(\cdot)-\mathcal K\varphi(s,B_s^{t_0,x_0},\cdot)$$
\begin{align*}
\delta f^{\P,t_0,x_0,\varphi}_s:=&\ f\left(s,B_s^{t_0,x_0},y_s^{\P,t_0,x_0,\varphi},z_s^{\P,t_0,x_0,\varphi},u_s^{\P,t_0,x_0,\varphi},\widehat a^{t_0}_s,\nu^{\P,t_0}_s\right)
\\
&-f\left(s,B_s^{t_0,x_0},\varphi(s,B_s^{t_0,x_0}),D\varphi(s,B_s^{t_0,x_0}),\mathcal K\varphi(s,B_{s^-}^{t_0,x_0},\cdot),\widehat a^{t_0}_s,\nu^{\P,t_0}_s\right).
\end{align*}

By a simple application of It\^o's formula under $\P$, we deduce that
\begin{align*}
\nonumber \delta y^{\P,t_0,x_0,\varphi}_{t_0}=&\ \int_{t_0}^{t_0+\eta}\mathcal L^{t_0,x_0}\varphi(s,\widehat a^{t_0}_s,\nu^{\P,t_0}_s)ds-\int_{t_0}^{t_0+\eta}\delta f^{\P,t_0,x_0,\varphi}_sds-\int_{t_0}^{t_0+\eta}\delta z^{\P,t_0,x_0,\varphi}_sd(B^{t_0})^{c,\P}_s\\
&-\int_{t_0}^{t_0+\eta}\int_E\delta u^{\P,t_0,x_0,\varphi}_s(e)\left(\mu_{B^{t_0}}(de,ds)-\nu^{\P,t_0}_s(de)ds\right).
\end{align*}
Using the same arguments that lead us to \reff{eq:M'}, and in particular using the Lipschitz properties of $f$, we can define a positive c\`adl\`ag process $M'$, whose supremum has finite moments of any order such that
\begin{align}\label{eq:visco1}
\nonumber\delta y^{\P,t_0,x_0,\varphi}_{t_0}&\geq \E^\P\left[\int_{t_0}^{t_0+\eta}M'_s\mathcal L^{t_0,x_0}\varphi(s,\widehat a^{t_0}_s,\nu^{\P,t_0}_s)ds\right]\\
\nonumber&\geq \E^\P\left[\int_{t_0}^{t_0+\eta}M'_s\left(\mathcal L^{t_0,x_0}\varphi(s,\widehat a^{t_0}_s,\nu^{\P,t_0}_s)-\mathcal L^{t_0,x_0}\varphi(t_0,\widehat a^{t_0}_s,\nu^{\P,t_0}_s)\right)ds\right]\\
&\hspace{0.9em}+\E^\P\left[\int_{t_0}^{t_0+\eta}M'_s\mathcal L^{t_0,x_0}\varphi(t_0,\widehat a^{t_0}_s,\nu^{\P,t_0}_s)ds\right].
\end{align}
By \reff{ineq:visco}, we know that the left-hand side of \reff{eq:visco1} is non-positive for every $\P\in\Pc^{t_0,\kappa}_h.$ Let us therefore fix some $(a,\nu)\in D^1_f\times (D^2_f\cap\mathfrak V^m)$ and consider the measure $\P(a,\nu)$ under which $\widehat a^{t_0}$ and $\nu^{\P(a,\nu),t_0}$ are equal to $a$ and $\nu$ respectively\footnote{Notice that since $a$ and $\nu$ are deterministic, such a measure does exist and is in $\Pc^{t_0,\kappa}_h$.}. We deduce that
\begin{align}\label{eq:visco2}
\nonumber&\Lc^{\P(a,\nu),t_0,x_0}\varphi(t_0,a,\nu)\int_{t_0}^{t_0+\eta}\E^{\P(a,\nu)}[M'_s]ds\\
&\leq \E^{\P(a,\nu)}\left[\underset{t_0\leq s\leq T}{\sup}M'_s\int_{t_0}^{t_0+\eta}\abs{\mathcal L^{t_0,x_0}\varphi(s,a,\nu)-\mathcal L^{t_0,x_0}\varphi(t_0,a,\nu)}ds\right]
\end{align}
Let us now estimate the right-hand side of \reff{eq:visco2}. We first have, using the fact that $\varphi\in C^{3}_b([0,T)\times\R^d)$ and Lemma \reff{lemma:estimates}
\begin{align}\label{phit}
&\nonumber\E^{\P(a,\nu)}\left[\underset{t_0\leq s\leq T}{\sup}M'_s\int_{t_0}^{t_0+\eta}\abs{\varphi_t(s, B_s^{t_0,x_0})-\varphi_t(t_0,x_0)}ds\right]\\
\nonumber&\leq C \E^{\P(a,\nu)}\left[\underset{t_0\leq s\leq T}{\sup}M'_s\int_{t_0}^{t_0+\eta}(s-t_0+\abs{B_s^{t_0}})ds\right]\\
\nonumber &\leq C \left(\eta^2+\int_{t_0}^{t_0+\eta}\left(\E^{\P(a,\nu)}\left[\underset{t_0\leq t\leq T}{\sup}(M'_t)^2\right]\right)^{1/2}\left(\E^{\P(a,\nu)}\left[\abs{B_s^{t_0}}^2\right]\right)^{1/2}ds\right)\\
&\leq C(\eta^2+\eta^{3/2}).
\end{align}
Similarly, we obtain that
\begin{align}\label{phixx}
\E^{\P(a,\nu)}\left[\underset{t_0\leq s\leq T}{\sup}M'_s\int_{t_0}^{t_0+\eta}\abs{\frac12\Tr{aD^2\varphi(s,B_s^{t_0,x_0})}-\frac12\Tr{aD^2\varphi(t_0,x_0)}}ds\right]\leq C(\eta^2+\eta^{3/2}).
\end{align}
We also have 
\begin{align}\label{phisaut}
\nonumber&\mathbb E^{\P(a,\nu)}\left[\underset{t_0\leq s\leq T}{\sup}M'_s\int_{t_0}^{t_0+\eta}\abs{\int_E(A\varphi)(B_{s^-}^{t_0,x_0}(\omega),e)\nu(de)-\int_E(A\varphi)(x_0,e)\nu(de)}ds\right]\\
&\leq C\underset{\nu\in D^2_{f}\cap\mathfrak V^m}{\sup}\left\{\int_E\abs{e}^2\nu(de)\right\}(\eta^2+\eta^{3/2}).
\end{align}
Finally, for the term involving the generator $f$, we have, using that $f$ is uniformly continuous in $(t,x)$, and uniformly Lipschitz in $(y,z,u)$
\begin{align}\label{phif}
\nonumber&\mathbb E^{\P(a,\nu)}\Big[\underset{t_0\leq s\leq T}{\sup}M'_s\int_{t_0}^{t_0+\eta}\Big|f(s,B_s^{t_0,x_0},\varphi(s,B_s^{t_0,x_0}),D\varphi(s,B_s^{t_0,x_0}),\Kc\varphi(s,B_s^{t_0,x_0},\cdot),a,\nu)\\
\nonumber &\hspace{3.5em}-f(t_0,x_0,\varphi(t_0,x_0),D\varphi(t_0,x_0),\Kc\varphi(t_0,x_0,\cdot),a,\nu)\Big|ds\Big]\\
\nonumber&\leq C\int_{t_0}^{t_0+\eta}\left(\rho(s-t_0)+\mathbb E^{\P(a,\nu)}\left[\underset{t_0\leq s\leq T}{\sup}M'_s\rho(\abs{B_s^{t_0}})\right]+\abs{s-t_0}+\mathbb E^{\P(a,\nu)}\left[\underset{t_0\leq s\leq T}{\sup}M'_s\abs{B_s^{t_0}}\right]\right)ds\\
&\hspace{0.9em}+C\int_{t_0}^{t_0+\eta}\int_E\E^{\P(a,\nu)}\left[\underset{t_0\leq s\leq T}{\sup}M'_s\abs{\Kc\varphi(s,B_s^{t_0,x_0},e)-\Kc\varphi(t_0,x_0,e)}\right](1\wedge\abs{e})\nu(de)ds
\end{align}
Since $\rho$ is concave, by Jensen's inequality, the first term on the right-side verifies
\begin{align*}
&\int_{t_0}^{t_0+\eta}\left(\rho(s-t_0)+\mathbb E^{\P(a,\nu)}\left[\underset{t_0\leq s\leq T}{\sup}M'_s\rho(\abs{B_s^{t_0}})\right]+\abs{s-t_0}+\mathbb E^{\P(a,\nu)}\left[\underset{t_0\leq s\leq T}{\sup}M'_s\abs{B_s^{t_0}}\right]\right)ds\\
&\leq C\left(\eta\rho(\eta)+\eta^2+\eta^{3/2}+\int_{t_0}^{t_0+\eta}\mathbb E^{\P(a,\nu)}\left[\underset{t_0\leq s\leq T}{\sup}M'_s\rho\left(\abs{B_s^{t_0}}\right)\right]ds\right).
\end{align*}
Then, we have by using twice the Cauchy-Schwarz inequality and using the fact that $\rho$ is concave and has linear growth
\begin{align*}
\mathbb E^{\P(a,\nu)}\left[\underset{t_0\leq s\leq T}{\sup}M'_s\rho\left(\abs{B_s^{t_0}}\right)\right]\leq &\ C\left(\E^{\P(a,\nu)}\left[\underset{t_0\leq s\leq T}{\sup}(M'_s)^4\right]\right)^{\frac14}\left(1+\E^{\P(a,\nu)}\left[\abs{B_s^{t_0}}^2\right]\right)^{\frac14}\\
&\times\left(\E^{\P(a,\nu)}\left[\rho\left(\abs{B_s^{t_0}}\right)\right]\right)^{\frac12}\\
\leq &\ C\rho^{1/2}\left(\E^{\P(a,\nu)}\left[\abs{B_s^{t_0}}\right]\right).
\end{align*}
Since $\rho^{1/2}$ is also concave, we obtain by Jensen's inequality that
\begin{align*}
&\int_{t_0}^{t_0+\eta}\left(\rho(s-t_0)+\mathbb E^{\P(a,\nu)}\left[\underset{t_0\leq s\leq T}{\sup}M'_s\rho(\abs{B_s^{t_0}})\right]+\abs{s-t_0}+\mathbb E^{\P(a,\nu)}\left[\underset{t_0\leq s\leq T}{\sup}M'_s\abs{B_s^{t_0}}\right]\right)ds\\
&\leq C\left(\eta\rho(\eta)+\eta^2+\eta^{3/2}+\eta\rho^{1/2}(\sqrt{\eta})\right).
\end{align*}
For the second term, we have similarly
\begin{align*}
&\int_{t_0}^{t_0+\eta}\E^{\P(a,\nu)}\left[\underset{t_0\leq s\leq T}{\sup}M'_s\int_E\abs{\Kc\varphi(s,B_s^{t_0,x_0},e)-\Kc\varphi(t_0,x_0,e)}\right](1\wedge\abs{e})\nu(de)ds\\
&\leq C\underset{\nu\in D^2_{f}\cap\mathfrak V^m}{\sup}\left\{\int_E(\abs{e}\wedge\abs{e}^2)\nu(de)\right\}(\eta^2+\eta^{3/2}),
\end{align*}
so that \reff{phif} becomes
\begin{align}\label{phif2}
\nonumber&\mathbb E^{\P(a,\nu)}\Big[\underset{t_0\leq s\leq T}{\sup}M'_s\int_{t_0}^{t_0+\eta}\Big|f(s,B_s^{t_0,x_0},\varphi(s,B_s^{t_0,x_0}),D\varphi(s,B_s^{t_0,x_0}),\Kc\varphi(s,B_s^{t_0,x_0},\cdot),a,\nu)\\
\nonumber &\hspace{3.5em}-f(t_0,x_0,\varphi(t_0,x_0),D\varphi(t_0,x_0),\Kc\varphi(t_0,x_0,\cdot),a,\nu)\Big|ds\Big]\\
&\leq C(\eta\rho(\eta)+\eta\rho^{1/2}(\sqrt{\eta})+\eta^2+\eta^{3/2}).
\end{align}
Hence, using \reff{phit}, \reff{phixx}, \reff{phisaut} and \reff{phif2} in \reff{eq:visco2}, we obtain
\begin{align*}
&\Lc^{\P(a,\nu),t_0,x_0}\varphi(t_0,a,\nu)\int_{t_0}^{t_0+\eta}\E^{\P(a,\nu)}[M'_s]ds\\
&\leq C\underset{\nu\in D^2_{f}\cap\mathfrak V^m}{\sup}\left\{\int_E(\abs{e}{\bf 1}_{\abs{e}>1}+\abs{e}^2)\nu(de)\right\}\eta\left(\rho(\eta)+\rho^{1/2}(\sqrt{\eta})+\eta+\eta^{1/2}\right).
\end{align*}
Dividing both sides by $\eta$, letting $\eta$ go to $0$ and using the fact that $M'$ is c\`adl\`ag, we deduce that
$$\E^{\P(a,\nu)}\left[M'_{t_0}\right]\Lc^{\P(a,\nu),t_0,x_0}\varphi(t_0,a,\nu)\leq 0,$$
which is the desired result since $M'$ is positive.

\vspace{0.3em}
{\bf Viscosity sub-solution}: Sub-solution property can be treated similarly, so we only detail the steps different from the proof for super-solution property. Let $\varphi\in C^{3}_b([0,T)\times\R^d)$ and $(t_0,x_0)\in[0,T)\times\R^d$ be such that
$$0=u(t_0,x_0)-\varphi(t_0,x_0)=\underset{(t,x)\in[0,T)\times \R^d}{\max}(u-\varphi)(t,x).$$
Fix some $\eta>0$ such that $t_0+\eta<T$. By \reff{partial.dpp2} with the constant family of stopping times $t_0+\eta$, we know that
\begin{align}\label{ineq:viscosub}
\nonumber\varphi(t_0,x_0)=u(t_0,x_0) &= \underset{\mathbb P\in\mathcal P^{t_0,\kappa}_h}{\sup} \mathcal Y^{\mathbb P,t_0,x_0}_{t_0}\left(t_0+\eta,u(t_0+\eta,B_{t_0+\eta}^{t_0,x_0})\right)\\
&\leq \underset{\mathbb P\in\mathcal P^{t_0,\kappa}_h}{\sup} \mathcal Y^{\mathbb P,t_0,x_0}_{t_0}\left(t_0+\eta,\varphi(t_0+\eta,B_{t_0+\eta}^{t_0,x_0})\right),
\end{align}
where we used the comparison theorem for BSDEJs in the last inequality.

\vspace{0.3em}
By the same arguments that lead us to \reff{eq:visco1}, and in particular using the Lipschitz properties of $f$ as Step $2$ of the proof of Theorem $4.1$ in \cite{kpz3}, we can define a positive c\`adl\`ag process $M$, whose supremum has finite moments of any order such that
\begin{align*}
\nonumber\delta y^{\P,t_0,x_0,\varphi}_{t_0}&\leq \E^\P\left[\int_{t_0}^{t_0+\eta}M_s\mathcal L^{t_0,x_0}\varphi(s,\widehat a^{t_0}_s,\nu^{\P,t_0}_s)ds\right]\\
\nonumber&\leq \E^\P\left[\int_{t_0}^{t_0+\eta}M_s\left(\mathcal L^{t_0,x_0}\varphi(s,\widehat a^{t_0}_s,\nu^{\P,t_0}_s)-\mathcal L^{t_0,x_0}\varphi(t_0,\widehat a^{t_0}_s,\nu^{\P,t_0}_s)\right)ds\right]\\
\nonumber&\hspace{0.9em}+\E^\P\left[\int_{t_0}^{t_0+\eta}M_s\mathcal L^{t_0,x_0}\varphi(t_0,\widehat a^{t_0}_s,\nu^{\P,t_0}_s)ds\right].\\
\nonumber&\leq \E^\P\left[\int_{t_0}^{t_0+\eta}M_s\left(\mathcal L^{t_0,x_0}\varphi(s,\widehat a^{t_0}_s,\nu^{\P,t_0}_s)-\mathcal L^{t_0,x_0}\varphi(t_0,\widehat a^{t_0}_s,\nu^{\P,t_0}_s)\right)ds\right]\\
&\hspace{0.9em}+\underset{(a,\nu)\in\mathbb S^{>0}_d\times\mathfrak V^m}{\sup}\mathcal L^{t_0,x_0}\varphi(t_0,a,\nu)\eta\E^\P\left[\underset{t_0\leq s\leq T}{\sup}M_s\right].
\end{align*}

Taking supremum on both sides of the above equation and by \reff{ineq:viscosub}, we deduce that
\begin{align*}
&-\frac1\eta\underset{\mathbb P\in\mathcal P^{t_0,\kappa}_h}{\sup} \E^{\P}\left[\underset{t_0\leq s\leq T}{\sup}M_s\int_{t_0}^{t_0+\eta}\abs{\mathcal L^{t_0,x_0}\varphi(s,\widehat a^{t_0}_s,\nu^{\P,t_0}_s)-\mathcal L^{t_0,x_0}\varphi(t_0,\widehat a^{t_0}_s,\nu^{\P,t_0}_s)}ds\right]\\ 
&\leq \underset{(a,\nu)\in\mathbb S^{>0}_d\times\mathfrak V^m}{\sup}\mathcal L^{t_0,x_0}\varphi(t_0,a,\nu)\underset{\mathbb P\in\mathcal P^{t_0,\kappa}_h}{\sup}\E^\P\left[\underset{t_0\leq s\leq T}{\sup}M_s\right].
\end{align*}
Exactly as in the proof of the super-solution property, we can show that the term on the l.h.s. of the above inequality tends to $0$ as $\eta$ goes to $0$, which finishes the proof by the positivity of $M$. 
\ep
\begin{appendix}

\section{Appendix}

\subsection{Technical proofs}

\proof[Proof of Proposition \ref{progdyn}]
W.l.o.g.,  we assume that $t_1=0$ and $t_2=t$. Thus, we have to prove
$$V_0(\omega)=\underset{\mathbb P\in \mathcal P^\kappa_H}{\sup}\mathcal Y_0^\mathbb P(t,V_t).$$

Denote $(y^\mathbb P,z^\mathbb P,u^\mathbb P):=(\mathcal  Y^\mathbb P(T,\xi),\mathcal Z^\mathbb P(T,\xi),\mathcal U^\mathbb P(T,\xi))$

\vspace{0.4em}
\rm{(i)} For any $\mathbb P\in \mathcal P^\kappa_H$, we know by Lemma \ref{lemme.technique} in the Appendix, that for $\mathbb P-a.e.$ $\omega\in\Omega$, the r.c.p.d. $\mathbb P^{t,\omega}\in\mathcal P^{t,\kappa}_H$. Now thanks to the paper of Tang and Li \cite{tangli}, we know that the solution of BSDEJs on the Wiener-Poisson space with Lipschitz generator can be constructed via Picard iteration. Thus, it means that at each step of the iteration, the solution can be formulated as a conditional expectation under $\mathbb P$. By the properties of the r.p.c.d. and Proposition \ref{relationhata}, this entails that
\begin{equation}\label{eq.picard}
y_t^\mathbb P(\omega)=\mathcal Y_t^{\mathbb P^{t,\omega},t,\omega}(T,\xi), \text{ for } \mathbb P-a.e.\text{ } \omega\in\Omega.
\end{equation}

Hence, by definition of $V_t$ and the comparison principle for BSDEJs, we get that $y_0^\mathbb P\leq \mathcal Y_0^\mathbb P(t,V_t)$. By arbitrariness of $\mathbb P$, this leads to 
$$V_0(\omega)\leq\underset{\mathbb P\in \mathcal P^\kappa_H}{\sup}\mathcal Y_0^\mathbb P(t,V_t).$$

\vspace{0.4em}
\rm{(ii)} For the other inequality, we proceed as in \cite{stz2}. Let $\mathbb P\in\mathcal P_H^\kappa$ and $\epsilon>0$. By separability of $\Omega$, there exists a partition $(E_t^i)_{i\geq 1}\subset \mathcal F_t$ such that $d_{S,t}(\omega,\omega')\leq \epsilon$ for any $i$ and any $\omega,\omega'\in E_t^i$.

\vspace{0.4em}
Now for each $i$, fix a $\widehat \omega_i\in E_t^i$ and let $\mathbb P^i_t$ be an $\epsilon-$optimizer of $V_t(\widehat\omega_i)$. If we define for each $n\geq 1$, $\mathbb P^n:=\mathbb P^{n,\epsilon}$ by
\begin{align}
 \mathbb P^n(E):=\mathbb E^\mathbb P\left[\sum_{i=1}^n\mathbb E^{\mathbb P^i_t}\left[1_E^{t,\omega}\right]1_{E_t^i}\right]+\mathbb P(E\cap\widehat E^n_t),\text{ where } \widehat E^n_t:=\bigcup_{i>n}E^i_t, \label{multi_proba}
\end{align}

then, by Lemma \ref{lemma_multi_proba}, we know that $\mathbb P^n\in \mathcal P^\kappa_H$. Besides, by Lemma \ref{unifcont} and its proof, we have for any $i$ and any $\omega\in E_t^i$
\begin{align*}
V_t(\omega)&\leq V_t(\widehat\omega_i)+C\rho(\epsilon)\leq\mathcal Y_t^{\mathbb P^i_t,t,\widehat\omega_i}(T,\xi)+\epsilon+C\rho(\epsilon)\\
&\leq \mathcal Y_t^{\mathbb P^i_t,t,\omega}(T,\xi)+\epsilon+C\rho(\epsilon)=\mathcal Y_t^{(\mathbb P^n)^{t,\omega},t,\omega}(T,\xi)+\epsilon+C\rho(\epsilon),
\end{align*}
where we used successively the uniform continuity of $V$ in $\omega$, the definition of $\mathbb P^i_t$, the uniform continuity of $\mathcal Y_t^{\mathbb P,t,\omega}$ in $\omega$ and finally the definition of $\mathbb P^n$.

\vspace{0.4em}
Then, it follows from \reff{eq.picard} that
\begin{equation}\label{eq.ggggg}
V_t\leq y_t^{\mathbb P^n}+\epsilon+C\rho(\epsilon),\text{ }\mathbb P^n-a.s. \text{ on } \bigcup_{i=1}^nE_t^i.
\end{equation}

Let now $(y^n,z^n,u^n):=(y^{n,\epsilon},z^{n,\epsilon},u^{n,\epsilon})$ be the solution of the following BSDEJ on $[0,t]$
\begin{align}
\label{grougrouref}
\nonumber y_s^n&=\left[y_t^{\mathbb P^n}+\epsilon+C\rho(\epsilon)\right]1_{\cup_{i=1}^nE_t^i}+V_t1_{\widehat E_t^n}+\int_s^t\widehat F^{\mathbb P^n}_r(y^n_r,z^n_r,u^n_r)dr-\int_s^tz^n_rdB^{\mathbb P^n,c}_r\\
&\hspace{0.9em}-\int_s^t\int_{E}u_r^n(x)\widetilde\mu^{\mathbb P^n}_B(dx,dr),\text{ }\mathbb P^n-a.s.
\end{align}

By the comparison principle for BSDEJs, we know that $\mathcal Y^\mathbb P_0(t,V_t)\leq y_0^n$. Then since $\mathbb P^n=\mathbb P$ on $\mathcal F_t$, the equality \reff{grougrouref} also holds $\mathbb P-a.s.$ Using the same arguments and notations as in the proof of Lemma \ref{unifcont}, we obtain
$$\abs{y_0^n-y_0^{\mathbb P^n}}^2\leq C\mathbb E^\mathbb P\left[\epsilon^2 +\rho(\epsilon)^2+\abs{V_t-y_t^{\mathbb P^n}}^21_{\widehat E_t^n}\right].$$

Then, by Lemma \ref{unifcont}, we have
\begin{align*}
\mathcal Y_0^\mathbb P(t,V_t)\leq y_0^n&\leq y_0^{\mathbb P^n}+C\left(\epsilon+\rho(\epsilon)+\scriptstyle\left(\mathbb E^\mathbb P\left[\Lambda_t^2 1_{\widehat E_t^n}\right]\right)^{\frac12}\right)\leq V_0+C\left(\epsilon+\rho(\epsilon)+\scriptstyle\left(\mathbb E^\mathbb P\left[\Lambda_t^2 1_{\widehat E_t^n}\right]\right)^{\frac12}\right).
\end{align*}

\vspace{0.4em}
Then it suffices to let $n$ go to $+\infty$, use the dominated convergence theorem, and let $\epsilon$ go to $0$.
\ep

\vspace{0.5em}
\proof[Proof of Lemma \ref{lem.cadlag}]
For each $\mathbb P$, we define $\widetilde V^\mathbb P:=V-\mathcal Y^\mathbb P(T,\xi).$ Then, we recall that we have $$\widetilde V^\mathbb P\geq 0,\ \mathbb P-a.s.$$ 

Now for any $0\leq t_1< t_2\leq T$, let $(y^{\mathbb P,t_2},z^{\mathbb P,t_2},u^{\mathbb P,t_2}):=(\mathcal Y^\mathbb P(t_2,V_{t_2}),\mathcal Z^\mathbb P(t_2,V_{t_2}),\mathcal U^\mathbb P(t_2,V_{t_2}))$. Once more, we remind that since solutions of BSDEJs can be defined by Picard iterations, we have by the properties of the r.p.c.d. that
$$\mathcal Y^\mathbb P_{t_1}(t_2,V_{t_2})(\omega)=\mathcal Y_{t_1}^{\mathbb P^{t_1,\omega},t_1,\omega}(t_2,V_{t_2}^{t_1,\omega}),\text{ for $\mathbb P-a.e.$ $\omega$}.$$

Hence, we conclude from Proposition \ref{progdyn} that $V_{t_1}\geq y_{t_1}^{\mathbb P,t_2},\text{ }\mathbb P-a.s.$ Denote
$$\widetilde y_t^{\mathbb P,t_2}:=y_t^{\mathbb P,t_2}-\mathcal Y^{\mathbb P}_t(T,\xi),\text{ }\widetilde z_t^{\mathbb P,t_2}:=\widehat a_t^{-1/2}(z_t^{\mathbb P,t_2}-\mathcal Z^{\mathbb P}_t(T,\xi)),\text{ }\widetilde u_t^{\mathbb P,t_2}:=u_t^{\mathbb P,t_2}-\mathcal U^{\mathbb P}_t(T,\xi).$$

Then $\widetilde V_{t_1}^\mathbb P\geq \widetilde y_{t_1}^{\mathbb P,t_2}$ and $(\widetilde y^{\mathbb P,t_2},\widetilde z^{\mathbb P,t_2},\widetilde u^{\mathbb P,t_2})$ satisfies the following BSDEJ on $[0,t_2]$
$$\widetilde y^{\mathbb P,t_2}_t=\widetilde V_{t_2}^\mathbb P+\int_t^{t_2}f_s^\mathbb P(\widetilde y^{\mathbb P,t_2}_s,\widetilde z^{\mathbb P,t_2}_s,\widetilde u_s^{\mathbb P,t_2})ds-\int_t^{t_2}\widetilde z^{\mathbb P,t_2}_sdW_s^\mathbb P-\int_t^{t_2}\int_{\mathbb R^d}\widetilde u^{\mathbb P,t_2}_s(x)\widetilde\mu^\mathbb P_B(dx,ds),$$
where
\begin{align*}
f_t^\mathbb P(\omega,y,z,u):&=\widehat F^\mathbb P_t(\omega,y+\mathcal Y^{\mathbb P}_t(\omega),\widehat a_t^{-1/2}(\omega)(z+\mathcal Z^{\mathbb P}_t(\omega)),u+{\mathcal U}^{\mathbb P}_t(\omega))\\
&\hspace{0.9em}-\widehat{F}^\mathbb P_t(\omega,\mathcal Y^{\mathbb P}_t(\omega),\mathcal Z^{\mathbb P}_t(\omega),\mathcal U^{\mathbb P}_t(\omega)).
\end{align*}

By the definition given in Royer \cite{roy}, we conclude that $\widetilde V^\mathbb P$ is a positive $f^\mathbb P$-supermartingale under $\mathbb P$. Since $f^\mathbb P(0,0,0)=0$, we can apply the downcrossing inequality proved in \cite{roy} to obtain classically that for $\mathbb P-a.e.$ $\omega$, the limit
$$\underset{r\in\mathbb Q\cup(t,T],r\downarrow t}{\lim}\widetilde V^\mathbb P_r(\omega)$$
exists for all $t$. Finally, since ${\mathcal Y}^{\mathbb P}$ is c\`adl\`ag, we obtain the desired result.
\ep

\begin{Lemma}\label{lemma.sde}
The SDE \reff{sde} with initial condition $Z_{\tau^0_\eps}=0$ has a unique solution on $[\tau_0^\eps,T]$.
\end{Lemma}

\proof
The proof follows the line of the proof of Example 4.5 in \cite{stz3}, and we provide it for comprehensiveness. For simplicity, we only prove the result for $\tau_0^\eps=0$. This does not pertain any loss of generality, since the general result can proved similarly by working on shifted spaces instead. We proceed by induction and let $Z^{0,\eps}$ be the solution of the SDE
$$Z^{\eps,0}_t=\int_0^t\left(a_0^\eps(Z^{\eps,0}_\cdot)\right)^{1/2}dB_s^{\P_{0,\tilde F^\eps},c}+\int_E\beta^\eps_0(Z^{\eps,0}_\cdot,x)\left(\mu_B(dx,ds)-\tilde F^\eps_s(dx)ds\right),\ \P_{0,\tilde F^\eps}-a.s.$$
Since $a_0^\eps$ and $\beta_0^\eps$ are actually $\mathcal F_0$-measurable, they are deterministic and thus $Z^{\eps,0}$ is indeed well-defined. Let then $\tilde \tau_0^\eps:=0$ and $\tilde \tau_1^\eps:=\tau_1^\eps(Z^{\eps,0})$. By Lemma $9.4$ in \cite{stz3}, $\tilde \tau_1^\eps$ is still an $\F$-stopping time. We pursue the construction by setting $Z^{\eps,1}_t:=Z^{\eps,0}_t$ for $t\in[0,\tilde\tau_1^\eps]$ as well as, for $t\geq \tilde\tau_1^\eps$
$$Z^{\eps,1}_t=Z^{\eps,0}_{\tilde\tau_1^\eps}+\int_{\tilde\tau_1^\eps}^t\left(a_1^\eps(Z^{\eps,1}_\cdot)\right)^{1/2}dB_s^{\P_{0,\tilde F^\eps},c}+\int_E\beta^\eps_1(Z^{\eps,1}_\cdot,x)\left(\mu_B(dx,ds)-\tilde F^\eps_s(dx)ds\right),\ \P_{0,\tilde F^\eps}-a.s.$$

Using the fact that $a_1^\eps$ and $\beta_1^\eps$ are $\Fc_{\tau_1}$-measurable, we can then argue as in \cite{stz3} to obtain that $a_1^\eps(Z^{\eps,1})=a_1^\eps(Z^{\eps,0})$ and $\beta_1^\eps(Z^{\eps,1},x)=\beta_1^\eps(Z^{\eps,0},x)$. Therefore $Z^{\eps,1}$ is also well defined. By repeating the procedure for $n\geq 2$, and since we know that there exists $N\in\N$ such that $\tau_n^\eps=T$ for $n\geq N$, after a finite number of steps we have constructed the unique strong solution $Z^\eps$ to the SDE on $[0,T)$. Since it is c\`adl\`ag, we extend it at time $T$ by setting $Z_T^\eps:=\underset{t\uparrow T}{\lim}\ Z_t^\eps$, which finishes the construction.
\ep

\subsection{The measures $\P^{\alpha,\beta}_F$}

\begin{Lemma}\label{stopping_time}
Let $\tau$ be an $\F$-stopping time. Let $\omega \in \Omega$, $s\geq \tau(\omega)$ and $H$ be a $\Fc^{\tau(\omega)}_s$-measurable random variable. There exists a $\Fc_s$-measurable random variable $\widetilde H$, such that 
\begin{align}
 H = \widetilde H^{\tau,\omega}. \label{def_S_tilde}
\end{align}
\end{Lemma}

\proof
\reff{def_S_tilde} means that for all $\widetilde \omega \in \Omega^{\tau(\omega)},$ $H(\widetilde \omega) =\widetilde H^{\tau,\omega}(\widetilde \omega) = \widetilde H(\omega\otimes_{\tau(\omega)}\widetilde \omega).$ We set 
\begin{align*}
\forall \ \omega_1 \in \Omega, \; \widetilde H(\omega_1):= H(\omega_1^{\tau(\omega)}). 
\end{align*}
Using the fact that for $s \geq \tau(\omega)$ 
\begin{align*}
 (\omega\otimes_{\tau(\omega)}\widetilde \omega)^{\tau(\omega)}(s)=\omega(\tau(\omega))+\widetilde \omega(s)-\omega(\tau(\omega))-\widetilde\omega(\tau(\omega))=\widetilde \omega(s),
\end{align*}
we have that $\widetilde H$ satisfies \reff{def_S_tilde} by construction. Indeed
\begin{align*}
\widetilde H^{\tau,\omega}(\widetilde \omega)= H\left[ \left( \omega\otimes_{\tau(\omega)}\widetilde \omega\right)^{\tau(\omega)}\right]=H(\widetilde \omega).
\end{align*}

Notice now that for any $\omega \in \Omega$, $\widetilde H: \Omega \rightarrow [\tau(\omega),T]$ is (Borel) measurable as a composition of measurable mappings. Finally, let us prove that $\widetilde H$ is $\Fc_s$-measurable. Since $H$ is $\Fc^{\tau(\omega)}$-measurable, there exists some measurable function $\phi$ such that for any $\widetilde \omega\in\Omega^{\tau(\omega)}$
$$H(\widetilde\omega)=\phi\left(B_t^{\tau(\omega)}(\widetilde\omega),\ \tau(\omega)\leq t\leq s\right).$$
Therefore, we have for any $\omega^1\in\Omega$
$$\widetilde H(\omega^1)=\phi\left(B_t(\omega^1)-B_\tau(\omega_1),\ \tau(\omega)\leq t\leq s\right),$$
which clearly implies that $\widetilde H$ is indeed $\Fc_s$-measurable.
%
\ep

\vspace{0.4em}
\begin{Lemma}\label{lemme.technique}
Let $\P \in \overline{\mathcal P}_{S}$ and $\tau$ be an $\F$-stopping time. Then
\begin{equation*}
 \P^{\tau,\omega} \in \overline{\mathcal P}_{S}^{\tau(\omega)} \text{ for } \; \P-a.e. \; \omega \in \Omega.
\end{equation*}
\end{Lemma}

\proof

\underline{{\bf{Step $1$}}}: Let us first prove that
\begin{align}
 (\P_{0,F})^{\tau,\omega} = \P_{\tau(\omega),F^{\tau,\omega}}, \;\; \P_{0,F}\text{-a.s. on } \Omega, \label{egalite.probas}
\end{align}
where $(\P_{0,F})^{\tau,\omega}$ denotes the probability measure on $\Omega^{\tau}$, constructed from the regular conditional probability distribution (r.c.p.d.) of $\P_{0,F}$ for the stopping time $\tau$, evaluated at $\omega$, and $ \P_{\tau(\omega),F^{\tau,\omega}}$ is the unique solution of the martingale problem $(\P^1,\tau(\omega),T,I_d,F^{\tau,\omega})$, where $\P^1$ is such that $\P^1(B^{\tau}_{\tau}=0)=1$.

\vspace{0.4em}
We recall that thanks to Remark 2.2 in \cite{kpz3}, it is enough to show that outside a $\P_{0,F}$-negligible set, the shifted processes $M^{\tau}, J^{\tau},Q^{\tau}$ (which are defined in this Remark)  are $(\P_{0,F})^{\tau,\omega}$-local martingales. In order to show this, for any $\omega$ in $\Omega$, and any $t\geq s\geq \tau(\omega)$, take any $\Fc_s^{\tau(\omega)} $-measurable random variable $H$. By Lemma \ref{stopping_time}, there exists a $\Fc_s $-measurable random variable $\widetilde{H}$ such that $H=\widetilde H^{\tau,\omega}$. Then, following the definitions in Subsection \ref{parag.notations}, we have
\begin{align*}
\Delta B_t^{\tau,\omega}(\widetilde{\omega}) &= \Delta B_t (\omega \otimes_{\tau} \widetilde{\omega}) = \Delta (\omega \otimes_{\tau} \widetilde{\omega})(t) = \Delta \omega_t \mathbf{1}_{\{t \leq \tau \}} + \Delta \widetilde{\omega}_t \mathbf{1}_{\{t > \tau \}},
\end{align*}
and for $t \geq \tau$
\begin{align*}
B_t (\omega \otimes_{\tau} \widetilde{\omega}) &=  (\omega \otimes_{\tau} \widetilde{\omega})(t) = \omega_{\tau} +\widetilde{\omega}_t  = B_{\tau}(\omega)+B_t^{\tau}(\widetilde{\omega}).
\end{align*}

From this we obtain
\begin{align*}
M_t^{\tau,\omega}(\widetilde{\omega}) = M_t (\omega \otimes_{\tau} \widetilde{\omega})=& \ B_t (\omega \otimes_{\tau} \widetilde{\omega}) - \sum_{u \leq t} \mathbf{1}_{\abs{\Delta B_u (\omega \otimes_{\tau} \widetilde{\omega}) }>1} \Delta B_u (\omega \otimes_{\tau} \widetilde{\omega})\\
 &+ \int_{0}^t\int_E x \mathbf{1}_{\abs{x}>1}F_u(\omega \otimes_{\tau} \widetilde{\omega},dx) \,du\\
=&  \ B_t^{\tau}(\widetilde{\omega}) + B_t(\omega) - \sum_{u \leq \tau} \mathbf{1}_{\abs{\Delta \omega_u}>1} \Delta \omega_u - \sum_{\tau < u \leq t} \mathbf{1}_{\abs{\Delta B_u^{\tau}(\widetilde{\omega})}>1} \Delta B_u^{\tau}(\widetilde{\omega})\\
&+ \int_{0}^{\tau}\int_E x \mathbf{1}_{\abs{x}>1}F_u(\omega,dx)\,du +\int_{\tau}^t\int_E x \mathbf{1}_{\abs{x}>1}F_u^{\tau,\omega}(\widetilde{\omega},dx)\,du\\
 =& \ M_t^{\tau}(\widetilde{\omega}) + M_{\tau}(\omega), \;\; \forall \omega \in \Omega.
\end{align*}
Localizing if necessary, we can now compute
\begin{align}
\E^{(\P_{0,F})^{\tau,\omega}} \left[ HM_t^{\tau} \right] &= \E^{(\P_{0,F})^{\tau,\omega}} \left[ H(M_t^{\tau,\omega} - M_{\tau}(\omega)) \right]\nonumber \\
&=\E^{(\P_{0,F})^{\tau,\omega}} \left[ \widetilde H^{\tau,\omega}M_{t}^{\tau,\omega} \right] - \E^{(\P_{0,F})^{\tau,\omega}} \left[ \widetilde H^{\tau,\omega}\right] M_{\tau}(\omega)\nonumber\\
&= \E^{\P_{0,F}}_{\tau}[\widetilde HM_{t}](\omega)- \E^{\P_{0,F}}_{\tau}[\widetilde H](\omega)M_{\tau}(\omega), \text{for $\mathbb P_{0,F}$-a.e. $\omega$}\nonumber\\
 &=\E^{\P_{0,F}}_{\tau}[\widetilde HM_{s}](\omega)- \E^{\P_{0,F}}_{\tau}[\widetilde H](\omega)M_{\tau}(\omega), \text{for $\mathbb P_{0,F}$-a.e. $\omega$}\nonumber\\
 &=\E^{(\P_{0,F})^{\tau,\omega}} \left[ HM_s^{\tau} \right]
, \text{for $\mathbb P_{0,F}$-a.e. $\omega$,} \label{eq1_lemme_technique}
\end{align}

where we use the fact that $M$ is a $\P_{0,F}$-local martingale.  
Since $H$ is arbitrary, we have that $M^{\tau}$ is a $(\P_{0,F})^{\tau,\omega}$-local martingale for $\mathbb P_{0,F}$-a.e. $\omega$. 


\vspace{0.4em}
We treat the case of the process $J^{\tau}$ analogously and write
\begin{align*}
J_t^{\tau,\omega}(\widetilde{\omega}) =& \ \left(M_t^{\tau,\omega}(\widetilde{\omega}) \right)^2 - t - \int_0^{\tau} \int_E x^2 F_u(\omega,dx)du - \int_{\tau}^t \int_E x^2 F_u^{\tau,\omega}(\widetilde{\omega},dx)du \\
=& \ \left(M_t^{\tau}(\widetilde{\omega})\right)^2 + \left(M_{\tau}(\omega)\right)^2 +2 M_t^{\tau}(\widetilde{\omega}) M_{\tau}(\omega) - (t-\tau) - \int_{\tau}^t \int_E x^2 F_u^{\tau,\omega}(\widetilde{\omega},dx)du\\
&- \int_0^{\tau} \int_E x^2 F_u(\omega,dx)du - \tau\\
=& \ J_t^{\tau}(\widetilde{\omega}) +  J_{\tau}(\omega) + 2 M_t^{\tau}(\widetilde{\omega}) M_{\tau}(\omega).
\end{align*}

Then we can compute the expectation, for $\P_{0,F}$-a.e. $\omega$
\begin{align*}
\E^{(\P_{0,F})^{\tau,\omega}} \left[ HJ_t^{\tau} \right] &= \E^{(\P_{0,F})^{\tau,\omega}} \left[ HJ_t^{\tau,\omega} - 2H M_t^{\tau} M_{\tau}(\omega) \right] - \E^{(\P_{0,F})^{\tau,\omega}} \left[ H\right]J_{\tau}(\omega)\\
&=\E^{\P_{0,F}}_\tau \left[\widetilde HJ_t\right](\omega)-\E^{(\P_{0,F})^{\tau,\omega}} \left[ H M_s^\tau M_\tau(\omega)\right]-\E^{(\P_{0,F})^{\tau,\omega}} \left[ H\right]J_{\tau}(\omega)\\
&=\E^{\P_{0,F}}_\tau \left[\widetilde HJ_s\right](\omega)-\E^{(\P_{0,F})^{\tau,\omega}} \left[ H M_s^\tau M_\tau(\omega)\right]-\E^{(\P_{0,F})^{\tau,\omega}} \left[ H\right]J_{\tau}(\omega)\\
&=\E^{(\P_{0,F})^{\tau,\omega}} \left[ HJ_s^{\tau} \right].
\end{align*}

$J^{\tau}$ is then a $(\P_{0,F})^{\tau,\omega}$-local martingale $\text{for $\mathbb P_{0,F}$-a.e. $\omega$}$. Finally, we do the same kind of calculation for $Q^{\tau}$, and we obtain
\begin{align*}
Q_t^{\tau,\omega}(\widetilde{\omega}) =&\ \int_0^t \int_E g(x) \mu_B(\omega \otimes_{\tau} \widetilde{\omega},dx,du) - \int_0^t \int_E g(x) F_u^{\tau,\omega}(\widetilde{\omega},dx)\,du\\
=& \ \int_0^{\tau} \int_E g(x) \mu_B(\omega,dx,du) + \int_{\tau}^t \int_E g(x) \mu_{B^{\tau}}(\widetilde{\omega},dx,du) \\
&- \int_0^{\tau} \int_E g(x) F_u(\omega,dx)du - \int_{\tau}^t \int_E g(x) F_u^{\tau,\omega}(\widetilde{\omega},dx)\,du\\
=& \ Q_t^{\tau}(\widetilde{\omega}) + Q_{\tau}(\omega).
\end{align*}

And again we compute the expectation over the $\widetilde{\omega} \in \Omega^{\tau}$, under the measure $(\P_{0,F})^{\tau,\omega}$
\begin{align*}
\E^{(\P_{0,F})^{\tau,\omega}} \left[ HQ_t^{\tau} \right] &= \E^{(\P_{0,F})^{\tau,\omega}} \left[ \widetilde H^{\tau,\omega}Q_t^{\tau,\omega} - HQ_{\tau}(\omega) \right]\\
&= \E^{\P_{0,F}}_{\tau}[\widetilde HQ_t](\omega)- \E^{(\P_{0,F})^{\tau,\omega}} \left[ H\right] Q_{\tau}(\omega),\text{ for $\mathbb P_{0,F}$-a.e. $\omega$}\\
&= \E^{(\P_{0,F})^{\tau,\omega}} \left[ HQ_s^{\tau} \right], \text{ for $\mathbb P_{0,F}$-a.e. $\omega$.}
\end{align*}

We have the desired result, and conclude that (\ref{egalite.probas}) holds true. We can now deduce that for any $(\alpha,\beta)\in\mathcal D\times \Rc_F$
\begin{equation}
 \P^{\tau(\omega),\alpha^{\tau,\omega},\beta^{\tau,\omega}}_{F^{\tau,\omega}} \in \overline{\mathcal P}_{S}^{\tau(\omega)}, \; \P_{0,F}\text{-a.s. on } \Omega. \label{step1}
\end{equation}

Indeed, if $(\alpha,\beta) \in \Dc \times \Rc_F$, then $(\alpha^{\tau,\omega}, \beta^{\tau,\omega}) \in \Dc^{\tau(\omega)} \times \Rc_{F^{\tau(\omega)}}^{\tau(\omega)}$, because for $\P_{0,F}-a.e.\ \omega$ and for Lebesgue almost every $s\in [\tau(\omega),T]$ 
$$
\int_{\tau(\omega)}^T \abs{\alpha_s^{\tau,\omega}(\widetilde{\omega})}ds < \infty, (\P_{0,F})^{\tau,\omega}-a.s. \ (\text{and thus } \P_{\tau(\omega),F^{\tau,\omega}}-a.s.),$$
$$\abs{\beta^{\tau,\omega}_s}(\widetilde \omega,x)\leq C(1\wedge \abs{x}), F^{\tau,\omega}(\widetilde \omega,dx)-a.e.,\ \text{for } \P_{\tau(\omega),F^{\tau,\omega}}-a.e. \ \widetilde \omega,$$
$$x\longmapsto \beta^{\tau,\omega}_s(\widetilde \omega, x) \text{ is strictly monotone for } F^{\tau,\omega}(\widetilde \omega,dx)-a.e. \  x,\ \text{for } \P_{\tau(\omega),F^{\tau,\omega}}-a.e. \ \widetilde \omega.$$

\vspace{0.4em}
\underline{{\bf Step $2$}}: We define $\widetilde{\tau}:= \tau \circ X^{\alpha,\beta}$, $\widetilde{\alpha}^{\tau,\omega} := \alpha^{\widetilde{\tau},\zeta_{\alpha,\beta}(\omega)}$, $\widetilde{F}^{\tau,\omega} := F^{\widetilde{\tau},\zeta_{\alpha,\beta}(\omega)}$ and $\widetilde{\beta}^{\tau,\omega} := \beta^{\widetilde{\tau},\zeta_{\alpha,\beta}(\omega)}$ where $\zeta_{\alpha,\beta}$ is a measurable map such that $B=\zeta_{\alpha,\beta}(X^{\alpha,\beta})$, $\P_{0,F}$-a.s. We refer to Lemma $2.2$ in \cite{stz2} (which can be proved similarly in our setting) for the existence of $\zeta_{\alpha,\beta}$. Moreover, 
$\widetilde{\tau}$ is an $\mathbb F$-stopping time and we have $\tau = \widetilde{\tau} \circ \zeta_{\alpha,\beta}$, $\P_{F}^{\alpha,\beta}-a.s.$ by definition. Using (\ref{step1}), we deduce 
\begin{equation*}
 \P^{\tau(\omega),\widetilde{\alpha}^{\tau,\omega},\widetilde{\beta}^{\tau,\omega}}_{\widetilde F^{\tau,\omega}} \in \overline{\mathcal P}_{S}^{\tau(\omega)}, \ \P^{\alpha, \beta}_F\text{-a.s. on } \Omega.
\end{equation*}

\vspace{0.4em}
\underline{{\bf Step $3$}}: We show that 
\begin{equation*}
 \E^{\P^{\alpha, \beta}_F} \left[ \phi\left(B_{t_1\wedge \tau}, \dots, B_{t_n\wedge \tau}\right)\psi\left(B_{t_1},\dots,B_{t_n}\right)\right] = \E^{\P^{\alpha, \beta}_F} \left[ \phi\left(B_{t_1\wedge \tau}, \dots, B_{t_n\wedge \tau}\right)\psi_{\tau}\right]
\end{equation*}
for every $0<t_1<\dots<t_n\leq T$, every continuous and bounded functions $\phi$ and $\psi$ and  where
\begin{equation*}
 \psi_{\tau}(\omega) = \E^{\P^{\tau(\omega),\widetilde{\alpha}^{\tau,\omega},\widetilde{\beta}^{\tau,\omega}}_{\widetilde{F}^{\tau,\omega}}} \left[\psi(\omega(t_1), \dots, \omega(t_k), \omega(t)+B_{t_{k+1}}^t, \dots, \omega(t)+B_{t_{n}}^t) \right],
\end{equation*}
for $t:=\tau(\omega) \in [t_k, t_{k+1})$.

\vspace{0.4em}
Recall that $\P^{\tau(\omega),\widetilde{\alpha}^{\tau,\omega},\widetilde{\beta}^{\tau,\omega}}_{\widetilde{F}^{\tau,\omega}}$ is defined by $\P^{\tau(\omega),\widetilde{\alpha}^{\tau,\omega},\widetilde{\beta}^{\tau,\omega}}_{\widetilde{F}^{\tau,\omega}} = \P_{\tau(\omega),\widetilde F^{\tau,\omega}} \circ \left(X^{\widetilde{\alpha}^{\tau,\omega},\widetilde{\beta}^{\tau,\omega}}\right)^{-1}$. For simplicity, we denote $$\widetilde \P:=\P_{\tau(\omega),\widetilde F^{\tau,\omega}}.$$ We then have
\begin{equation*}
\begin{split}
& \psi_{\tau}(\omega) = \E^{\widetilde{\P}} \Big[\psi\Big(\omega(t_1), \dots, \omega(t_k), \omega(t)+\int_t^{t_{k+1}} \left(\alpha_s^{\widetilde{\tau},\zeta_{\alpha,\beta}(\omega)}\right)^{1/2}d(B_s^{\tau(\omega)})^{\widetilde \P, c} \\
 &\hspace{0.5em}+ \int_t^{t_{k+1}}\int_{E} \beta_s^{\widetilde\tau,\zeta_{\alpha,\beta}(\omega)}(x) \left (\mu_{B^{\tau(\omega)}}(dx,ds) - F_s^{\widetilde{\tau},\zeta_{\alpha,\beta}(\omega)}(dx)ds\right) , \dots, \omega(t)\\
 &\hspace{0.5em}+\int_t^{t_{n}} \left(\alpha_s^{\widetilde{\tau},\beta_{\alpha}(\omega)}\right)^{1/2}d(B_s^{\tau(\omega)})^{\widetilde \P,c}+ \int_t^{t_{n}}\int_{E} \beta_s^{\widetilde\tau,\zeta_{\alpha,\beta}(\omega)}(x) \left (\mu_{B^{\tau(\omega)}}(dx,ds) - F_s^{\widetilde{\tau},\zeta_{\alpha,\beta}(\omega)}(dx)ds\right)  \Big) \Big].
\end{split}
\end{equation*}

Denote $$\widehat \P:= \P_{\widetilde \tau,F^{\widetilde{\tau},\omega}}.$$ Then, $\forall \ \omega \in \Omega$, if $t:= \widetilde{\tau}(\omega) = \tau\left(X^{\alpha,\beta}(\omega)\right) \in [t_k, t_{k+1}[$, 
\begin{align}\label{eq_psi_2}
\nonumber &\psi_{\tau}\left(X^{\alpha,\beta}(\omega)\right)=\E^{\widehat \P} \Big[ \psi\Big( X_{t_1}^{\alpha,\beta}(\omega), \dots, X_{t_k}^{\alpha,\beta}(\omega), X_{t}^{\alpha,\beta}(\omega) +\int_t^{t_{k+1}} \left(\alpha_s^{\widetilde{\tau},\omega}\right)^{1/2}d(B_s^{\widetilde \tau(\omega)})^{\widehat \P,c} \\
\nonumber 
&\hspace{0.5em}+ \int_t^{t_{k+1}}\int_{E} \beta_s^{\widetilde\tau,\omega}(x) \left (\mu_{B^{\widetilde\tau(\omega)}}(dx,ds) - F_s^{\widetilde{\tau},\omega}(dx)ds\right) , \dots, X_{t}^{\alpha,\beta}(\omega)+\int_t^{t_{n}} \left(\alpha_s^{\widetilde{\tau},\omega}\right)^{1/2} d(B_s^{\widetilde{\tau}(\omega)})^{\widehat \P, c} \\
 &\hspace{0.5em}+ \int_t^{t_{n}}\int_{E} \beta_s^{\widetilde\tau,\omega}(x)\left (\mu_{B^{\widetilde\tau(\omega)}}(dx,ds) - F_s^{\widetilde{\tau},\omega}(dx)ds\right) \Big) \Big].
\end{align}

We remark that for every $\omega \in \Omega$,
\begin{align*}
 &\alpha_s(\omega) = \alpha_s\left(\omega \otimes_{\widetilde{\tau}(\omega)} \omega^{\widetilde{\tau}(\omega)}\right) = \alpha_s^{\widetilde{\tau},\omega}\left(\omega^{\widetilde{\tau}(\omega)}\right),
\end{align*}
and similar relations hold for both $F$ and $\beta$.

\vspace{0.5em}
By definition, the $(\P_{0,F})^{\widetilde{\tau},\omega}$-distribution of $B^{\widetilde{\tau}(\omega)}$ is equal to the $(\P_{0,F})^{\omega}_{\widetilde{\tau}}$-distribution of $(B_{\cdot} - B_{\widetilde{\tau}(\omega)})$. Therefore since by \reff{egalite.probas}, $(\P_{0,F})^{\widetilde \tau,\omega} = \P_{\widetilde\tau(\omega),F^{\widetilde \tau,\omega}}, \;\; \P_{0,F}\text{-a.s. on } \Omega$, (\ref{eq_psi_2}) then becomes
\begin{equation*}
\begin{split}
 \psi_{\tau}\left(X^{\alpha,\beta}(\omega)\right)=\ &\E^{(\P_{0,F})^{\omega}_{\widetilde{\tau}}} \Big[ \psi\Big( X_{t_1}^{\alpha,\beta}(\omega), \dots, X_{t_k}^{\alpha,\beta}(\omega), X_{t}^{\alpha,\beta}(\omega) +\int_t^{t_{k+1}} \alpha_s^{1/2}(B)d(B_s^{(\P_{0,F})^\omega_{\widetilde \tau},c}) \\
 &\hspace{-2.7em}+ \int_t^{t_{k+1}}\int_{E} \beta_s(x) (\mu_{B}(dx,ds) - F_s(dx)ds) , \dots, X_{t}^{\alpha,\beta}(\omega)+\int_t^{t_{n}} \alpha_s^{1/2}(B)d(B_s^{(\P_{0,F})^\omega_{\widetilde \tau},c}) \\
 &\hspace{-2.7em}+ \int_t^{t_{n}}\int_{E} \beta_s(x) (\mu_{B}(dx,ds) - F_s(dx)ds) \Big) \Big] \\
 =\ & \E^{(\P_{0,F})^{\omega}_{\widetilde{\tau}}} \Big[ \psi\Big( X_{t_1}^{\alpha,\beta}, \dots, X_{t_k}^{\alpha,\beta}, X_{t_{k+1}}^{\alpha,\beta}, \dots, X_{t_n}^{\alpha,\beta} \Big) \Big] \\
 =\ & \E^{\P_{0,F}} \left[ \psi\Big( X_{t_1}^{\alpha,\beta}, \dots, X_{t_k}^{\alpha,\beta}, X_{t_{k+1}}^{\alpha,\beta}, \dots, X_{t_n}^{\alpha,\beta} \Big) | \mathcal{F}_{\widetilde{\tau}} \right](\omega), \; \; \P_{0,F}\text{-a.s. on } \Omega.
\end{split}
\end{equation*}

Then we have 
\begin{equation*}
\begin{split}
& \E^{\P^{\alpha, \beta}_F} \left[ \phi\left(B_{t_1\wedge \tau}, \dots, B_{t_n\wedge \tau}\right)\psi_{\tau}\right] = \E^{\P_{0,F}} \left[ \phi\Big( X_{t_1\wedge \widetilde{\tau}}^{\alpha,\beta}, \dots, X_{t_n\wedge \widetilde{\tau}}^{\alpha,\beta} \Big)\psi_{\widetilde{\tau}}(X^{\alpha,\beta}) \right] \\
& = \E^{\P_{0,F}} \left[ \phi\Big( X_{t_1\wedge \widetilde{\tau}}^{\alpha,\beta}, \dots, X_{t_n\wedge \widetilde{\tau}}^{\alpha,\beta} \Big) \E^{\P_{0,F}} \left[ \psi\Big( X_{t_1}^{\alpha,\beta}, \dots, X_{t_k}^{\alpha,\beta}, X_{t_{k+1}}^{\alpha,\beta}, \dots, X_{t_n}^{\alpha,\beta} \Big) | \mathcal{F}_{\widetilde{\tau}} \right] \right] \\
&= \E^{\P_{0,F}} \left[ \phi\Big( X_{t_1\wedge \widetilde{\tau}}^{\alpha,\beta}, \dots, X_{t_n\wedge \widetilde{\tau}}^{\alpha,\beta} \Big) \psi\Big( X_{t_1}^{\alpha,\beta}, \dots, X_{t_k}^{\alpha,\beta}, X_{t_{k+1}}^{\alpha,\beta}, \dots, X_{t_n}^{\alpha,\beta} \Big) \right] \\
& = \E^{\P^{\alpha, \beta}_F} \left[ \phi\left(B_{t_1\wedge \tau}, \dots, B_{t_n\wedge \tau}\right) \psi\left(B_{t_1}, \dots, B_{t_n}\right)\right].
\end{split}
\end{equation*}

\underline{{\bf Step $4$}}: Now we prove that $\P^{\tau,\omega} = \P^{\widetilde{\alpha}^{\tau,\omega},\widetilde{\beta}^{\tau,\omega}}_{\widetilde{F}^{\tau,\omega}}$, $\P$-a.s. on $\Omega$.

\vspace{0.4em}
By definition of the conditional expectation, 
\begin{equation*}
 \psi_{\tau}(\omega) = \E^{\P^{\tau,\omega}} \left[\psi(\omega(t_1), \dots, \omega(t_k), \omega(t)+B_{t_{k+1}}^t, \dots, \omega(t)+B_{t_{n}}^t) \right], \text{ $\P^{\alpha,\beta}_F$-a.s.},
\end{equation*}
where $t:= \tau(\omega) \in [t_k, t_{k+1}[$, and where the $\P^{\alpha,\beta}_F$-null set can depend on $(t_1, \dots, t_n)$ and $\psi$, but we can choose a common null set by standard approximation arguments. 

\vspace{0.4em}
Then by a density argument we obtain 
\begin{equation*}
 \E^{\P^{\tau,\omega}}\left[\eta\right] = \E^{\P^{\widetilde{\alpha}^{\tau,\omega},\widetilde{\beta}^{\tau,\omega}}_{\widetilde{F}^{\tau,\omega}}} \left[\eta\right], \text{ for $\mathbb P^{\alpha,\beta}_F$-a.e. $\omega$},
\end{equation*}
for every bounded and $\mathcal{F}_T^{\tau(\omega)}$-measurable random variable $\eta$. This implies $\P^{\tau,\omega} = \P^{\widetilde{\alpha}^{\tau,\omega},\widetilde{\beta}^{\tau,\omega}}_{\widetilde{F}^{\tau,\omega}}$, $\P$-a.s. on $\Omega$. And from the Step 1 we deduce that $\P^{\tau,\omega} \in \overline{\mathcal P}_{S}^{\tau(\omega)}$.
\ep

\vspace{0.4em}
\begin{Lemma}
\label{lemma_multi_proba} 
We have $\P^n \in \Pc^{\kappa}_H$, where $\P^n$ is defined by (\ref{multi_proba}).
\end{Lemma}

\proof
Since by definition, $\P^i_t \in \Pc^{t,\kappa}_H$ and $\P \in \Pc^\kappa_H$, we have $\P^i_t = \P^{\alpha^i, \beta^i}_{F^i}$ and $\P = \P^ {\alpha, \beta}_F$, for $F^i\in\Vc^t $, $(\alpha^i, \beta^i) \in \Dc^t\times\Rc_{F^i}^t$ and $(F,\alpha,\beta) \in \Vc\times\Dc\times\Rc_F$, $i=1, \dots, n$. 
Next we define
\begin{align*}
\overline{\alpha}_s &:= \alpha_s \mathbf{1}_{[0,t)}(s) + \left[\sum_{i=1}^n \alpha^i_s \mathbf{1}_{E^i_t}(X^{\alpha,\beta}) + \alpha_s \mathbf{1}_{\Eh^n_t}(X^{\alpha,\beta}) \right] \mathbf{1}_{[t,T]}(s), \text{  and}\\
\overline{F}_s &:= F_s \mathbf{1}_{[0,t)}(s) + \left[\sum_{i=1}^n F^i_s \mathbf{1}_{E^i_t}(X^{\alpha,\beta}) + F_s \mathbf{1}_{\Eh^n_t}(X^{\alpha,\beta}) \right] \mathbf{1}_{[t,T]}(s), \text{  and}\\
\overline{\beta}_s &:= \beta_s \mathbf{1}_{[0,t)}(s) + \left[\sum_{i=1}^n \beta^i_s \mathbf{1}_{E^i_t}(X^{\alpha,\beta}) + \beta_s \mathbf{1}_{\Eh^n_t}(X^{\alpha,\beta}) \right] \mathbf{1}_{[t,T]}(s).
\end{align*}
Now following the arguments in the proof of step 3 of Lemma \ref{lemme.technique}, we can prove that for any $0<t_1< \dots < t_k=t<t_{k+1}<t_n$ and any continuous and bounded functions $\phi$ and $\psi$,
\begin{align*}
&\E^{\P^{\alpha,\beta}_F} \left[ \phi(B_{t_1}, \dots, B_{t_k}) \sum_{i=1}^n \E^{\P^{\alpha^i,\beta^i}_{F^i}} \left[\psi(B_{t_1}, \dots, B_{t_k},B_t + B^t_{t_{k+1}}, \dots, B_t + B^t_{t_n})\right]\mathbf{1}_{E^i_t} \right]\\
 &= \E^{\P^{\overline{\alpha},\overline{\beta}}_{\overline{F}}} \left[ \phi(B_{t_1}, \dots, B_{t_k})\psi(B_{t_1}, \dots, B_{t_n}) \right].
\end{align*}
This implies that $\P^n = \P^{\overline{\alpha},\overline{\beta}}_{\overline{F}} \in \overline{\Pc}_{S}$. And since all the probability measures $\P^i$ satisfy the requirements of Definition \ref{set_proba_shift}, we have $\P^n  \in \Pc^{\kappa}_H$.
\ep

\subsection{A weak dynamic programming principle}

\proof [of Lemma \ref{partial.dpp}.]
The proof follows closely the steps of Proposition $5.14$ and Lemma $6.2$ and $6.4$ in \cite{stz}.

\vspace{0.3em}
Let us first fix $\P$ and $X$ and denote $\tau = \tau^{\P}$ for simplicity. By the a priori estimates and recalling definition \ref{grandlambda.def}, we have $\abs{u(t,x)} \leq C \Lambda(t,x)$ and then we can assume w.l.o.g that $\abs{X} \leq C \Lambda(\tau,B_{\tau}^{t,x})$, $\P$-a.s. Then by assumption \ref{assump.markovian}, $X$ is in $\L^2(\P)$ and hence $\mathcal Y^{\mathbb P,t,x}_t\left(\tau,X\right)$ is well defined. Now by \reff{eq.picard}, 
\begin{align*}
\mathcal Y^{\mathbb P,t,x}_t\left(T,g(B_T^{t,x})\right) = \mathcal Y^{\mathbb P,t,x}_t\left(\tau,\mathcal Y^{\mathbb P^{\tau,\omega},\tau(\omega),B_{\tau}^{t,x}(\omega)}_{\tau}\left( T,g(B_T^{\tau(\omega),B_{\tau}^{t,x}(\omega)})\right)\right).
\end{align*}
By Lemma \ref{lemme.technique}, $\P^{\tau,\omega} \in \mathcal{P}_h^{\kappa,\tau(\omega)}$, for $\P$-a.e. $\omega \in \Omega^t$. Next by definition of $u$ and $X$,
\begin{align*}
\mathcal Y^{\mathbb P^{\tau,\omega},\tau(\omega),B_{\tau}^{t,x}(\omega)}_{\tau}\left( T,g(B_T^{\tau(\omega),B_{\tau}^{t,x}(\omega)})\right) &\leq u\left(\tau(\omega),B_{\tau}^{t,x}(\omega) \right)\leq X(\omega), \; \text{ for $\P$-a.e. $\omega \in \Omega^t$.}
\end{align*}
Finally, the comparison theorem for standard BSDEJs implies that
\begin{align*}
\mathcal Y^{\mathbb P,t,x}_t\left(T,g(B_T^{t,x})\right) \leq \mathcal Y^{\mathbb P,t,x}_t\left(\tau,X\right),
\end{align*}
and taking the supremum over $\P$ yields inequality \ref{partial.dpp1}.

\vspace{0.4em}
Let us now prove equality \reff{partial.dpp2} when we know in addition that $g$ lower semi-continuous . By Proposition \ref{prop.reg}, $u$ is lower semi-continuous in the variables $(t,x)$, from the right in $t$, and therefore $u$ is measurable and $u(\tau,B_{\tau}^{t,x})$ is $\mathcal{F}_{\tau}$-measurable. This implies that inequality \ref{partial.dpp1} holds for the particular choice $X=u(\tau,B_{\tau}^{t,x})$.

\vspace{0.4em}
We now prove the reverse inequality. The first step is to show that 
\begin{align}
\mathcal Y^{\mathbb P,s,x}_s\left(t,\phi(B_t^{s,x})\right) \leq u(s,x),\label{dpp_step1}
\end{align}
for any $s<T$ and fixed $t\in (s,T]$ and for any continuous real-valued function $\phi$ such that $-\Lambda(t,\cdot)\leq \phi(\cdot)\leq \Lambda(t,\cdot)$. Following step 2 of the proof of Theorem $3.5$ in \cite{bt}, for any $\epsilon >0$, we can find sequences $(x_i,r_i)_{i\geq 1} \subset \R^d \times (0,T]$ and $\P_i\in \mathcal{P}_h^{\kappa,t}$ such that 
\begin{align*}
\mathcal Y^{\mathbb \P_i,t,\cdot}_t\left(T,g(B_T^{t,\cdot})\right)\geq \phi(\cdot)-\epsilon \text{ on } B(x_i,r_i),
\end{align*}
where $B(x_i,r_i)$ denote the open balls centered at $x_i$ with radius $r_i$, which form a countable cover of $\R^d$. From this we can build a partition $(A_i)_{i\geq 1}$ of $\R^d$ defined by $A_i:=B(x_i,r_i) \backslash \cup_{j<i}B(x_j,r_j)$ and a partition of $\Omega$ defined in the following way
\begin{align*}
E^i:=\{B_t^{s,x} \in A_i\}, \; i\geq 1, \text{ and } \, \widehat E^n:= \cup_{i>n}E^i, \; n\geq 1.
\end{align*}
Then 
\begin{align*}
\Omega = \left(\cup_{i=1}^n E^i \right)\cup \widehat E^n  \;\text{ and } \, \lim_{n\to +\infty}\P(\widehat E^n)=0.
\end{align*}
Exactly as in \reff{multi_proba}, we define 
\begin{align*}
 \mathbb P^n(E):=\mathbb E^\mathbb P\left[\sum_{i=1}^n\mathbb E^{\mathbb P^i}\left[1_E^{t,\omega}\right]1_{E^i}\right]+\mathbb P(E\cap\widehat E^n),\text{ for all } E \in \mathcal F_T^s.
\end{align*}
It follows from Lemma \ref{lemma_multi_proba} and the definition of $\P^n$ that $\P^n$ coincide with $\P$ on $\mathcal F_t$, $\P^n \in \mathcal P_h^{s,\kappa}$ and $(\P^n)^{s,\omega}=\P^i$, for $\P$-a.e. $\omega$ in $E_i$, $1\leq i\leq n$.
Then we have
\begin{align*}
\mathcal Y^{\P^n,s,x}_t\left(T,g(B_T^{s,x})\right)(\omega) &= \mathcal Y^{\P^i,t,B_t^{s,x}(\omega)}_t\left(T,g(B_T^{t,B_t^{s,x}(\omega)})\right)\\
&\geq \phi(B_t^{s,x}(\omega))-\epsilon \;\;\;\; \text{ for $\P$-a.e. $\omega$ in $E_i$, $1\leq i\leq n$. }
\end{align*}
We use now the comparison theorem for BSDEs:
\begin{align*}
u(s,x) &\geq \mathcal Y^{\P^n,s,x}_s\left(T,g(B_T^{s,x})\right) = \mathcal Y^{\mathbb \P,s,x}_s\left(t,\mathcal Y^{\P^n,s,x}_t(T,g(B_T^{s,x}))\right)\\
&\geq \mathcal Y^{\P,s,x}_s\left(t,\phi(B_t^{s,x})-\epsilon \right) \mathbf{1}_{(\widehat E^n)^c} + \mathcal Y^{\P^n,s,x}_t(T,g(B_T^{s,x}))\mathbf{1}_{(\widehat E^n)}
\end{align*}
Since $\epsilon$ was arbitrary, we can use the stability results for BSDEs to conclude that inequality \reff{dpp_step1} holds true. 

\vspace{0,4em}
The rest of the proof is exactly the same as the end of the proof of Lemma $6.4$ in \cite{stz}: we use the lower semi-continuity of $u(t,\cdot)$ to approximate it by an increasing sequence of continuous functions and we use inequality \reff{dpp_step1} to prove the desired inequality for constant stopping times. Then we have the inequality for stopping times taking finitely many values by a simple backward induction and for stopping times taking countably many values by a limiting argument. Finally we approximate an arbitrary stopping time by a decreasing sequence of stopping times taking countably many values, and then we only need the function $u$ to be lower semi-continuous in $(t,x)$ from the right in $t$, to proceed exactly as in \cite{stz}.
\ep


\end{appendix}

\end{document}